\documentclass[reqno]{amsart}
\usepackage{amscd}
\usepackage[dvips]{graphics}
\usepackage{color,amsfonts,amssymb,amsmath,indentfirst,amsthm,mathrsfs,lmodern}
\usepackage{bm}
\usepackage{mathtools}
\usepackage{verbatim}
\usepackage{hyperref}
\def\beq{\begin{equation}}
\def\eeq{\end{equation}}
\def\ba{\begin{array}}
\def\ea{\end{array}}

\allowdisplaybreaks

\numberwithin{equation}{section}
\newenvironment{abs}{\textbf{Abstract}\mbox{  }}{ }
\newenvironment{key words}{\textbf{Keywords}\mbox{  }}{ }
\newtheorem*{theorem*}{Theorem A}
\newtheorem{theorem}{Theorem}[section]
\newtheorem{lemma}[theorem]{Lemma\hspace*{0.1em}}
\newtheorem{definition}[theorem]{Definition\hspace*{0.1em}}
\newtheorem{proposition}[theorem]{Proposition\hspace*{0.1em}}

\newtheorem{remark}[theorem]{Remark\hspace*{0.1em}}
\renewenvironment{proof}{{\noindent\textbf{Proof.}}}{\hspace*{0.1em}\hfill$\Box$}

\begin{document}
\title[Sharp Onofri trace inequality and $n$-Liouville equation]{\textbf{Sharp Onofri trace inequality on the upper half space and quasi-linear Liouville equation with Neumann boundary}
}
\author[J. Dou]{Jingbo Dou}
\address{Jingbo Dou, School of Mathematics and Statistics, Shaanxi Normal University, Xi'an, Shaanxi, 710119, People's Republic of China}
\email{jbdou@snnu.edu.cn}

\author[Y. Han]{Yazhou Han}
\address{Yazhou Han, School of Mathematics and Statistics, China Jiliang Normal University, Hangzhou, Zhejiang, 310018, People's Republic of China}
\email{yazhou\_han@msn.com}

\author[S. Yuan]{Shuang Yuan}
\address{Shuang Yuan, School of Mathematics and Statistics, Shaanxi Normal University, Xi'an, Shaanxi, 710119, People's Republic of China}
\email{yuanshuang@snnu.edu.cn}

\author[Y. Zhou]{Yang Zhou}
\address{Yang Zhou, School of Mathematical Sciences, University of Science and Technology of China, Hefei, Anhui, 230026, People's Republic of China}
\email{zy19700816@mail.ustc.edu.cn}

\date{}
\maketitle

\noindent
\begin{abs}
In this paper, we establish a sharp Onofri trace inequality on the upper half space $\overline{\mathbb R_+^n} (n\geq 2)$ by considering the limiting case of Sobolev trace inequality and classify its extremal functions on a suitable weighted Sobolev space. For this aim, by the Serrin-Zou type identity and the Pohozaev type identity, we show the classification of the solutions for a quasi-linear Liouville equation with Neumann boundary which is closely related to the Euler-Lagrange equation of the Onofri trace inequality and it has independent research value. The regularity and asymptotic estimates of solutions to the above equation are essential to discuss.
\end{abs}\\
\begin{key words}
Onofri trace inequality, quasi-linear Liouville equation, classification of solution, regularity, asymptotic behavior, Serrin-Zou type identity, Pohozaev type identity 
\end{key words}\\
\textbf{Mathematics Subject Classification(2020).}
26D10, 35B08, 35B40, 35J92, 35B53, 35J66.

\tableofcontents

\indent

\section{\textbf{Introduction}\label{Section 1}}
The classical Onofri inequality and Onofri trace inequality play an important role in the prescribed Gaussian curvature problem (or Nirenberg problem) for surfaces with and without boundary. Recently, del Pino and Dolbeault \cite{PD2013} established Onofri inequality in $\mathbb{R}^n\ (n\geq 3)$ as the limiting case of a family of  interpolation inequalities. 
In this work, we will extend Onofri trace inequality on $\overline{\mathbb{R}_+^2}$
to the higher-dimensional Euclidean upper half-space based on sharp Sobolev trace inequality and study the classification of the extremal functions. To achieve this, we shall investigate the quasi-linear Liouville equation with Neumann boundary. We state the background and motivations as follows. 

\subsection{Onofri inequality and Liouville equation}
The classical Onofri inequality on the unit sphere $\mathbb S^2$ states that
\begin{equation}\label{I-SO-2D} 
\log\int_{\mathbb S^2}e^v\mathrm d\sigma_2- \int_{\mathbb S^2}v\mathrm d\sigma_2\le\frac14 \int_{\mathbb S^2}|\nabla v|^2\mathrm d\sigma_2
\end{equation}
for any $v\in H^1(\mathbb S^2)$, where $\mathrm d\sigma_2$ denotes normalized surface measure. 
Moser \cite{M1970} first established the rough inequality \eqref{I-SO-2D} without a sharp constant, and then Onofri \cite{O1982} proved the inequality \eqref{I-SO-2D} 
based on an earlier result of Aubin \cite{A1979} and showed that the extremal function $u\in H^1(\mathbb S^2)$ is only constant 
by modulo conformal transformations. Hong \cite{H1986} also showed the extremal function $u\in H^1(\mathbb S^2)$ and the best constant. The inequality \eqref{I-SO-2D} is closely related to the prescribed Gaussian curvature (or Nirenberg problem), such as \cite{CC1986, CD1987, CG2023, CY1987, CY1988, GM2018, WX1999, WX2009}.  

Using the stereographic projection from $y=(y_1,y_2,y_3)\in \mathbb S^2$ to $x=(x_1,x_2)\in\mathbb R^2$,
inequality \eqref{I-SO-2D} can be reformulated into the Euclidean Onofri inequality as
\begin{equation}\label{I-EO-2D} 
\log\int_{\mathbb R^2}e^u\mathrm d\nu_2- \int_{\mathbb R^2}u\mathrm d\nu_2\le\frac{1}{16\pi} \int_{\mathbb R^2}|\nabla u|^2\mathrm dx
\end{equation} for any $u(x)\in L^1(\mathbb R^2,\mathrm d\nu_2)$ and $\nabla u\in L^2(\mathbb R^2,\mathrm dx)$, where $\mathrm d\nu_2=\frac{\mathrm dx}{\pi(1+|x|^2)^2}$. The survey article by Dolbeault, Esteban and Jankowiak \cite{DEJ2015} proved the Onofri inequality \eqref{I-EO-2D} through various limiting procedures based on some functional inequalities and rigidity methods (see also \cite{DET2008}).

Let $X, Y, y, z\in\mathbb R^n$ and $1<p<\infty$, define

$$R_{p}(X,Y):=|X+Y|^p-|X|^p-p|X|^{p-2}X\cdot Y, $$
$$H_n(y,z):=R_n\big(-\frac{n|y|^{-\frac{n-2}{n-1}}}{1+|y|^{\frac{n}{n-1}}}y,\frac{n-1}{n}z\big).$$
For any open set $\Omega\subseteq\mathbb R^n$ and $1<p<\infty$, denote Sobolev space as
$$W^{1,p}(\Omega):=\{u\in L^p(\Omega):\nabla u\in L^p(\Omega)\},$$
$$W^{1,p}_{\operatorname{loc}}(\Omega):=\{u\in W^{1,p}(U):\forall U\subset\subset\Omega\},$$
 and
$$W^{1,p}_{\operatorname{loc}}(\overline{\Omega}):=\{
u\in W^{1,p}(V\cap\overline{\Omega}):  V\subset\subset\mathbb R^n,V\cap\overline{\Omega}\neq\emptyset\}.$$
Recently, del Pino and Dolbeault \cite{PD2013} extended inequality \eqref{I-EO-2D} to the higher-dimensional Euclidean space as
\begin{equation}\label{I-EO-inq}
\log \int_{\mathbb R^n}e^{w}\mathrm d\nu_n-\int_{\mathbb R^n}w\mathrm d\nu_n\leq\beta_n\int_{\mathbb R^n}H_n(x,\nabla w)\mathrm dx,
\end{equation}
for any $w\in C^{\infty}_{c}(\mathbb R^n)$, where $\nu_n(x)=\frac n{|\mathbb S^{n-1}|}\frac1 {\big(1+|x|^{\frac n{n-1}}\big)^n}$, $\mathrm d\nu_n:=\nu_n(x)\mathrm dx$ and $
\beta_n=\frac{n^{1-n}\Gamma(\frac{n}{2})}{2(n-1)\pi^{\frac{n}{2}}}$. They pointed out that the equality in \eqref{I-EO-inq} is achieved
by constants, showed that $\beta_n$ is optimal 
and proposed an open question: \textbf{Are there nontrivial extremal functions of \eqref{I-EO-inq}?}

Agueh, Boroushaki, and Ghoussoub \cite{ABG2017} also established inequality \eqref{I-EO-inq} by the mass transport method and showed that the infimum of the inequality \eqref{I-EO-inq} was attained at $w=0$ in the space $ W^{1,n}(\mathbb R^n)$ with $n\ge 2$. We also refer to some improved versions of \eqref{I-EO-inq} involving singular weighted function in \cite{DEJ2015,LL2018}.

Recently, Borgia, Cingolani, Mancini \cite{BCM2025} extended the inequality \eqref{I-EO-inq} to the weighted Sobolev space
\begin{align*} 
W_{\nu_n}(\mathbb R^n):=&\{u\in L^1(\mathbb R^n,\mathrm d\nu_n)\,: \ |\nabla u|\in L^n(\mathbb R^n,\mathrm dx), \\
&\quad|\nabla u|^2|\nabla \log\nu_n|^{n-2}\in L^1(\mathbb R^n,\mathrm dx) \},
\end{align*}
which is not equivalent to the space $ W^{1,n}(\mathbb R^n)$. Moreover, it is proved that $W_{\nu_n}(\mathbb R^n)\backslash$ $W^{1,n}(\mathbb R^n) \neq \emptyset$ in \cite{BCM2025}. 

Recently, Dou, Wu and Yuan \cite{DWY2025}  formulated the Euler-Lagrange equation of \eqref{I-EO-inq}, up to a constant multiple, as
\begin{equation}\label{I-Onofri-ELn}
-\operatorname{div}\big[|X +\nabla w|^{n-2}(X +\nabla w)\big]=\frac{e^{w}}{\left(1+|x|^{\frac{n}{n-1}}\right)^n}, \qquad \text{in}~\mathbb R^n,
\end{equation}
where $X =-\frac{n^2}{n-1}\frac{|x|^{-\frac{n-2}{n-1}}}{1+|x|^{\frac{n}{n-1}}}x=\nabla [\log\frac{1}{\left(1+|x|^{\frac{n}{n-1}}\right)^n}]$. Then, they classified the solutions of \eqref{I-Onofri-ELn} and got a family of nontrivial extremal functions in the space $W_{\nu_n}(\mathbb R^n)$. This gave a positive answer to the open question.

In fact, let $u(x)=w(x)+\log\frac{1}{\left(1+|x|^{\frac{n}{n-1}}\right)^n}$ and equation \eqref{I-Onofri-ELn} can be simplified to the following quasi-linear Liouville equation (also called $n$-Liouville equation)
\begin{equation}\label{I-Liouv-eq}
 -\Delta_n u=-\operatorname{div}(|\nabla u|^{n-2}\nabla u)=e^{u},\quad\text{ in }~\mathbb R^n.
\end{equation}
When $n=2$, equation \eqref{I-Liouv-eq} is closely related to the prescribed Gaussian curvature (see \cite{CC1986, CD1987, CY1987, CY1988} and the references therein).  
Under the finite mass condition  $\int_{\mathbb R^2}e^{u}\mathrm dx<\infty$, Chen and Li \cite{CL1991} classified the $C^2$-smooth solutions of \eqref{I-Liouv-eq} by the method of moving planes 
and proved that $u$ must be of the form
\[
u(x)=\log\frac{8\lambda^2}{(1+\lambda^2|x-x_0|^2)^2}
\]
with $\lambda>0,~x,~x_0\in \mathbb R^2$. See also \cite{CY1997, WX1999, WX2009}. It implies that there exists a family of nontrivial $C^2$-extremal functions of equation \eqref{I-EO-inq} when $n=2$.

Under the assumption of finite mass $\int_{\mathbb R^n}e^{u}\mathrm dx<\infty$ and $u\in W^{1,n}_{\mathrm{loc}}(\mathbb R^n)$, Esposito \cite{E2018} classified entire 
solutions of equation \eqref{I-Liouv-eq} as follows 
\begin{equation}\label{I-Solution-E}
u(x)=\log\frac{n\big(\frac{n^2}{n-1}\big)^{n-1}\lambda^n}{\left(1+\lambda^{\frac{n}{n-1}}|x-x_0|^{\frac{n}{n-1}}\right)^n}
\end{equation}
for any $\lambda>0$ and $x_0,x\in \mathbb R^n$. He exploited conformal invariance of $\Delta_n$ under the Kelvin transformation, Pohozaev identity and isoperimetric argument to classify the solutions of equation \eqref{I-Liouv-eq}.
This technique was early used in the Liouville system by Chanillo and Kiessling \cite{CK1995}. Recently, it has also been applied to classify solutions to the anisotropic $n$-Liouville equation by Ciraolo and Li \cite{CL2024} and the anisotropic $N$-Liouville equation on the convex cones by Dai, Gui, and Luo \cite{DGL2024}. This approach differs from the classical strategy, the method of moving planes/spheres such as \cite{CL1991, CW1994}. What's more, some progress has been made in the study of classification of solutions for the Liouville equation in a Riemannian manifold. By assuming both the finite mass condition and the asymptotic lower bound, Cai and Lai \cite{CL2024a} proved the classification results of Liouville equation on complete surfaces with nonnegative Gauss curvature. After that, removing the condition of finite mass and still assuming the asymptotic lower bound, Ciraolo, Farina and Polvara \cite{CFP2024} gave the proof of classification result of solutions to Liouville equation on a complete, connected, non-compact, boundaryless Riemannian surface with nonnegative Ricci curvature using the method of $P$-function. Very recently, Sun and Wang \cite{SW2025} extended the classification results in \cite{CFP2024}  to higher dimensions. Under only the assumption of asymptotic lower bound
\[
u(x)\ge-\frac{n}{n-1}\log\big[r(x)G^\frac12(r(x))\big]-C,\quad r(x)\ge C^{-1}
\]
for some positive constant $C$, where $G$ is any positive nondecreasing function satisfying $\int_{c}^\infty\frac{\mathrm ds}{sG(s)}=\infty$ for some $c>1$ and $r(x)=\operatorname{dist}(x, o)$ is the distance function from $x$ to some fixed point $o \in M^n$, they classified solutions of the quasi-linear Liouville equation \eqref{I-Liouv-eq}  in a complete $n$-dimensional Riemannian manifold $(M^n,g)$ with nonnegative Ricci curvature by the $P$-function approach.

\subsection{Onofri trace inequality and our motivations}
The Onofri trace inequality on the unit disk $\overline{\mathbb{B}^2}\subset\mathbb R^2$ as stated in \cite{OPS1988} asserts that
\begin{equation}\label{I-Ot-2D}
\log\int_{\mathbb S^{1}}e^{u}\mathrm d\sigma_{1}-\int_{\mathbb S^{1}}u\mathrm d\sigma_{1}\leq\frac{1}{4\pi}\int_{\mathbb{B}^2}|\nabla u|^2\mathrm dx
\end{equation}
for any $u\in W^{1,2}(\mathbb{B}^2)$, where $\mathrm{d}\sigma_{1}$ denotes the normalized 
curve measure. 
The inequality is in fact the first inequality of Lebedev-Milin inequalities \cite{D1983}, which play an important role in the proof of Bieberbach's conjecture. From a more geometric angle of view, Osgood, Phillips and Sarnak \cite{OPS1988} proved the inequality \eqref{I-Ot-2D} and discussed the extremals of functional determinants on manifolds with boundary (such as, the extremal function is achieved on the disc). They also showed that the equality holds in \eqref{I-Ot-2D} if and only if $u=\log|\tau'|+c$, where $\tau:\mathbb{B}^2\rightarrow\mathbb{B}^2$ is a M\"{o}bius transformation, $|\tau'|$ is the Jacobian and $c\in\mathbb R$.

Under the conformal transformation $F:\mathbb{R} _+^2\ni y=(y_1,y_2)\mapsto x=(x_1,x_2)\in \mathbb B^2$ defined as
$$(x_1,x_2)=\frac{2(y_1,y_2+1)}{y_1^2+(y_2+1)^2}-(0,1),$$
inequality \eqref{I-Ot-2D} is equivalent to 
\begin{equation}\label{Ineq-trace-2D}
\log \int_{\mathbb{R}} e^w\mathrm d\mu_2 -\int_{\mathbb{R}} w\mathrm d\mu_2\leq \frac 1{4\pi} \int_{\mathbb{R}_+^2} |\nabla w|^2\mathrm dy,
\end{equation}
where $\mu_2(y):=\frac 1 {\pi(y_1^2+(1+y_2)^2)}$ and $\mathrm d\mu_2:=\mu_2(y_1,0)\mathrm dy_1=\frac{\mathrm dy_1}{\pi(y_1^2+1)}$. 

Consider the Liouville equation with Neumann boundary on the upper half-plane
\begin{equation}\label{I-2hL-eq}
\begin{cases}
\Delta u+be^{u}=0,&\text{ in }\mathbb R_+^2,\\
\frac{\partial u}{\partial t}=ce^{\frac{u}{2}},&\text{ on }\partial\mathbb R_+^2, 
\end{cases}
\end{equation}  
where $b, c\in\mathbb R$. Up to the transformation $w(y)=u(y)+\log\mu_2(y)$ (see details below), for $b=0,~c=-1$, equation \eqref{I-2hL-eq} is just the Euler-Lagrange equation of the Onofri trace inequality \eqref{Ineq-trace-2D}, for which there have been many classification results of solutions under assumption of finite mass 
\begin{equation}\label{I-2fm-assu}
\int_{\mathbb R_+^2}e^{u}\mathrm ds \mathrm dt+\int_{\partial\mathbb R_+^2}e^{\frac{u}{2}}\mathrm ds<\infty.
\end{equation}

Li and Zhu in \cite{LZ1995} provided the proof of the classification results for \eqref{I-2hL-eq} for the case $b=1$, $c\in\mathbb R$ by using the method of moving spheres. They proved that the $C^2(\mathbb R_+^2)\cap C^{1}(\overline{\mathbb R_+^2})$ solutions are the form of 
\begin{equation}\label{I-LZ-resu}
u(s,t)=\log\frac{8\lambda^2}{\big[\lambda^2+(s-s_0)^2+(t-t_0)^2\big]^2}
\end{equation}
for some $\lambda>0$, $s_0\in\mathbb R$ and $t_0=c\frac{\lambda}{\sqrt{2}}$.

Varying from the method of Li and Zhu, by the classical complex analysis technique, Ou in \cite{O2000} proved the classification results to \eqref{I-2hL-eq} for $b=0$ and $c<0$ only with the finite mass assumption $\int_{\mathbb R_+^2}e^{u}\mathrm ds\mathrm dt<\infty$. His results showed that the $C^2$-solutions $u$ must take the following form 
\begin{eqnarray}\label{I-Ou-resu}
u(s,t)=2\log\frac{2t_0}{(s-s_0)^2+(t+t_0)^2}+2\log\frac{2}{|c|},
\end{eqnarray}
where $s_0\in\mathbb R$ and $t_0>0$.

Later, Zhang in \cite{Z2003} weakened the assumption \eqref{I-2fm-assu} and classified the solutions to \eqref{I-2hL-eq} by the method of moving spheres. 
%
%
%
%
%
Later, G\'{a}lvez and Mira \cite{GM2009} generalized Zhang's results and found all the solutions to \eqref{I-2hL-eq} without any integral finiteness hypothesis using the complex analysis method.

An interesting question: Can we extend the Onofri trace inequality \eqref{Ineq-trace-2D}
to the upper half space $\overline{\mathbb{R}_+^n}~(n\geq 3)$ and classify the nontrivial extremal functions?

In this work, we exploit the above question, that is, we establish the sharp Onofri trace inequality on $\overline{\mathbb{R}_+^n}$ through the limiting argument of Sobolev trace inequality, and classify the quasi-linear Liouville equation with Neumann boundary  by integral identity method, which is related to Euler-Lagrange equation of the Onofri trace inequality. This leads to a full classification of the nontrivial extremal functions.

\subsection{Our main results}
In order to state our main results, we first introduce some notation.

Let $x=(x',t)\in\mathbb R_+^n=\mathbb R^{n-1}\times\{t>0\}$ and
$$\mathrm \mu_n(x)=\mathrm \mu_n(x',t):=\frac{2}{\sigma_{n-1}\big((1+t)^2+|x'|^2\big)^{\frac{n}{2}}},$$
where $\sigma_{n-1}:=|\mathbb S^{n-1}|=\frac{2\pi^{\frac n 2}}{\Gamma(\frac n 2)}$. 
Define
\begin{align*}
W_{\mu_n}(\mathbb R_+^n):=&\{w\in L^1(\partial\mathbb R_+^n,\mathrm d\mu_n(x'))
:\ |\nabla w|\in L^n(\mathbb R_+^n,\mathrm dx),\nonumber\\
&\quad |\nabla w|^2|\nabla \log\mu_n|^{n-2}\in L^1(\mathbb R_+^n,\mathrm dx)\},
\end{align*}
equipped with the norm
\begin{equation*}
\|w\|_{\mu_n}:=\int_{\partial\mathbb R_+^n}|w|\mathrm d\mu_n(x')+\|\nabla w\|_{L^n(\mathbb R_+^n)}+\big(\int_{\mathbb R_+^n}|\nabla w|^2|\nabla\log\mu_n|^{n-2}\mathrm dx\big)^{\frac 1 2},
\end{equation*}
where $\mathrm d\mu_n(x'):=\mu_n(x',0)\mathrm dx'=
\frac{2}{\sigma_{n-1}(1+|x'|^2)^{\frac{n}{2}}}\mathrm dx'$. 

For any integer $0\leq k\leq\infty$, denote
$$C_{c}^{k}(\overline{\mathbb R_+^n}):=\{u\in C^{k}(\overline{\mathbb R_+^n}),  \text{supp}\, u\subset\subset\mathbb R^n\}.$$

Write
\[
K_n(x,y):=R_n\big(-\frac{n(x', 1+t)}{(1+t)^2 + |x'|^2},y\big),\quad(x,y)\in \mathbb R_+^n\times\mathbb R^n,
\] 
where $R_n(X,Y)$ is defined in the previous section.
 
We firstly establish the Onofri trace inequality on the $n$-dimensional upper half space for $n\geq2$ as follows.
\begin{theorem}\label{Ot-inq-thm}
Suppose  $w\in W_{\mu_n}(\mathbb R_+^n)$, the following Onofri trace inequality holds
\begin{equation}\label{I-Ot-inq}
\log\big(\int_{\partial\mathbb R_+^n}e^{w}\mathrm d\mu_n(x')\big)- \int_{\partial\mathbb R_+^n}w\mathrm d\mu_n(x') \leq\alpha_n\int_{\mathbb R_+^n}K_n(x,\nabla w)\mathrm dx,
\end{equation}
where $\alpha_n:=\frac{2}{n^n\sigma_{n-1}}$ is the best constant.
Moreover, the equality holds in \eqref{I-Ot-inq} if and only if
\begin{equation}\label{Extre-f}
w(x)=\log\frac{\big(|x'|^2+(1+t)^2\big)^{\frac{n}{2}}\lambda}{\big(|x'-x_0'|^2+(t+\lambda)^2\big)^{\frac{n}{2}}}+\tilde C,
\end{equation}
for $x=(x',t)\in\overline{\mathbb R_+^n}$, $x_0'\in\mathbb R^{n-1}$, $\lambda>0$ and $\tilde C\in\mathbb R$.
\end{theorem} 

The Euler-Lagrange equation of the Onofri trace inequality \eqref{I-Ot-inq}, up to a constant multiple, can be written as
\begin{equation}\label{EL eq-1}
\begin{cases}
-\operatorname{div}\big[|X+\nabla w|^{n-2}(X+\nabla w)\big]=0, &\quad\text{in }\mathbb R_+^n,\\
|X+\nabla w|^{n-2}(X+\nabla w)\cdot\mathbf{n}=\mathrm{e}^{w}\mu_n(x',0),&\quad\text{on }\partial\mathbb R_+^n,
\end{cases}
\end{equation} 
where $X=-\frac{n(x', 1+t)}{(1+t)^2 + |x'|^2}=\nabla \big(\log \mu_n(x)\big)$.

Using the transformation $u(x)=\log \mu_n(x)+w(x)$, the above equation can be simplified to the following quasi-linear Liouville equation with Neumann boundary condition
\begin{equation}\label{I-hL-eq}
\begin{cases}
-\Delta_n u=0, &\quad \text{in}~\mathbb R_+^n,\\
|\nabla u|^{n-2}\frac{\partial u}{\partial t}=-e^{u}, &\quad \text{on}~\partial\mathbb R_+^n,\\
\end{cases}
\end{equation}
where $u\in W^{1,n}_{\operatorname{loc}}(\overline{\mathbb R^n_+}), n\ge2$.

\begin{definition}\label{Liouv-weak-solu-def}
We say that $u$ is a weak solution of equation \eqref{I-hL-eq} if ~$u\in W^{1,n}_{\operatorname{loc}}(\overline{\mathbb R_+^n})$ and
\begin{equation*}
\int_{\mathbb R_+^n}|\nabla u|^{n-2}\langle\nabla u,\nabla \varphi\rangle \mathrm dx=\int_{\partial\mathbb R_+^n}e^{u}\varphi \mathrm dx',
\end{equation*}
for all $\varphi\in W^{1,n}_0(\Omega\cap\overline{\mathbb R_+^n})\cap L^{\infty}(\Omega\cap \partial \mathbb R_+^n)$ with $\Omega\subset\mathbb R^n$ bounded.
\end{definition}

Under the assumption of finite mass
\begin{equation}\label{I-fm-assu}
\int_{\mathbb R_+^n}e^{\frac{n}{n-1}u}\mathrm dx+\int_{\partial\mathbb R_+^n}e^u\mathrm dx'<\infty,
\end{equation} 
we have the following classification theorem.
\begin{theorem}\label{hL-Classf-thm}
Let $u\in W^{1,n}_{\operatorname{loc}}(\overline{\mathbb R^n_+})$ be a weak solution of equation \eqref{I-hL-eq} satisfying the assumption \eqref{I-fm-assu}. Then $u$ must be the form
\begin{equation}\label{I-hL-Solution}
u(x)=\log \frac{n^{n-1}\lambda}{\big(|x'-x_0'|^2+(t+\lambda)^2\big)^{\frac{n}{2}}},
\end{equation}
for any $x=(x',t)\in\mathbb R_+^n$ and some $\lambda>0,~x_0'\in\mathbb R^{n-1}$.
\end{theorem}

\begin{remark}
(i) Combining Theorem \ref{hL-Classf-thm} with the transformation $u(x)=\log \mu_n(x)$ $+w(x)$, we can classify the solutions to \eqref{EL eq-1} and get all extremal functions \eqref{Extre-f} in $W_{\mu_n}(\mathbb R_+^n)$. Similar procedure applied to the case on $\mathbb{R}^n$ can give an affirmative answer to the open question proposed by del Pino and Dolbeault in \cite{PD2013}.
   
(ii) Constant functions $C\in W_{\mu_n}(\mathbb R_+^n)$ are extremal functions of \eqref{I-Ot-inq}, but $C\not\in W^{1,n}(\mathbb{R}_+^n)$. This indicates that $W_{\mu_n}(\mathbb R_+^n)$ is a natural Sobolev space to study the inequality \eqref{I-Ot-inq}.
\end{remark}




Our approach to proving Theorem \ref{hL-Classf-thm} employs integral identity method, which does not depend on conformal invariance, and also does not make use of moving planes/spheres. The strategy of the proof can be explained as follows. Following the idea of \cite{CL2024, DGL2024}, we establish the first-order regularity results and logarithmic asymptotic estimate. Based on these results, by the limiting argument, we show the asymptotic estimates of $u$ and $\nabla u$ at infinity. Using a similar method to that of \cite{Z2024}, namely the Caccioppoli-type inequality, we get the second-order regularity result and the asymptotic integral estimate on second-order derivatives. Finally, using the regularization argument as \cite{CFR2020} and \cite{Z2024} for solutions of approximating equations and establishing the Serrin-Zou differential identity (see \cite{SZ2002}) and the Pohozaev differential identity, we complete the classification of the solutions to \eqref{I-hL-eq}. This approach has also been applied to discuss the symmetry results for the critical anisotropic $p$-Laplacian equations in convex cones by Ciraolo, Figalli and Roncoroni \cite{CFR2020}.

The paper is organized as follows. In Section \ref{Section 2}, we establish the Onofri trace inequality on the upper half-space and show the best constant. In Sections \ref{Section 3}-\ref{Section 5}, we  classify the solutions of the quasi-linear Liouville equation \eqref{I-hL-eq}. More specifically, in Section \ref{Section 3}, we prove the $C^{1,\alpha }$ regularity to equation \eqref{I-hL-eq} by the idea of Brezis-Merle. In Section \ref{Section 4}, we give the sharp asymptotic estimates of $u$ and $\nabla u$ and the asymptotic integral estimates of the second-order derivatives of $u$. In Section \ref{Section 5}, we classify the solutions of \eqref{I-hL-eq} by employing the Serrin-Zou type identity and the Pohozaev type identity, and after that we establish a complete classification of the extremal functions of inequality \eqref{I-Ot-inq}.

After obtaining the classification of the notrivial extremal functions of inequality \eqref{I-Ot-inq},  we will discuss an interesting question in the future: the stability of the Onofri trace inequality \eqref{I-Ot-inq}.

Throughout this paper, we define $a(\xi):=|\xi|^{n-2}\xi$ for $\xi\in\mathbb R^n$ and denote by $C$ a generic positive constant that may vary in different terms and by $\mathbf{n}$ the unit outward normal vector on the boundary. For $x\in\mathbb R^n,~r>0$, the Euclidean ball of radius $r$ centered at $x$ is denoted by
$$B_{r}(x)=\{y\in\mathbb R^n\,:\,|x-y|<r\}$$ and $B_{r}^{+}(x):=B_{r}(x)\cap\mathbb R_+^n,~\Sigma_r(x):=B_{r}^{+}(x)\cap \partial\mathbb R_+^n$, $\partial B_r^+(x):=\partial B_r(x)\cap\mathbb R_+^n$. In particular, if $x=0$, we write $B_r=B_r(0),~B_r^+=B_r^+(0),~\Sigma_r=\Sigma_r(0),~\partial B_r^+=\partial B_r^+(0)$ unless otherwise specified. We denote by $\omega_n=|B_1|$ the volume of the unit ball $B_1\subset\mathbb R^n$.

 \section{\textbf{Sharp Onofri trace inequality on the upper half space}}\label{Section 2}
In this section, we establish the $n$-dimensional sharp Onofri trace inequality on the upper half space.

To do this, we first recall a sharp Sobolev trace inequality as 
\begin{equation}\label{O-St-inq}
(\int_{\partial\mathbb R_+^n}|u|^{\frac{p(n-1)}{n-p}}\mathrm dx')^{\frac{n-p}{p(n-1)}}\leq\frac{1}{S(n,p)}(\int_{\mathbb R_+^n}|\nabla u|^p\mathrm dx)^{\frac{1}{p}}
\end{equation}
for $u\in \mathcal D^{1,p}(\mathbb R_+^n):=\{u\in L^{\frac{np}{n-p}}(\mathbb R_+^n):\nabla u\in L^p(\mathbb R_+^n)\}$, where $n\geq 2$, $1<p<n$ and
$$S(n,p)=\big(\frac{n-p}{p-1}\big)^{\frac{p-1}{p}}\big[\frac1 2\sigma_{n-2}B(\frac{n-1}{2},\frac{n-1}{2(p-1)})\big]^\frac{p-1}{p(n-1)}.$$ 
Moreover, the equality holds in \eqref{O-St-inq} if and only if
\begin{equation}\label{O-St-extremal}
u(x',t)=\big(\frac{\lambda^{\frac{2}{p}}}{|x'-x_0'|^2+(t+\lambda)^2}\big)^{\frac{n-p}{2(p-1)}}
\end{equation} for $(x',t)\in\mathbb R_+^n, x_0'\in\mathbb R^{n-1}$ and $\lambda>0$. 
Escobar \cite{E1988} classified the extremal function when $p=2$ and proposed the conjecture: the extremal functions for $1<p<n$ is the form \eqref{O-St-extremal}. Also, he raised an open question in \cite{E1988}: Does every positive critical point to the Sobolev trace inequality have the form \eqref{O-St-extremal}? By means of mass transportation techniques, Nazaret \cite{N2006} and Maggi, Neumayer \cite{MN2017} obtained the best constants in the Sobolev trace inequality for $p\in(1, n)$, as well as the extremal functions, respectively. Recently, Zhou \cite{Z2024}  classified the positive solutions of the Euler-Lagrange equation of \eqref{O-St-inq} by Serrin-Zou's method and then answered Escobar's question.

\subsection{Onofri trace inequality} We exploit a similar approach to that in \cite{PD2013} by limiting procedures of sharp Sobolev trace inequality \eqref{O-St-inq} for any smooth compactly supported function, and then, arguing as \cite{BCM2025}, we extend the result to a suitable weighted Sobolev space $W_{\mu_n}({\mathbb R^n_+})$ by the density argument.

\begin{proposition}\label{Ot-inq-prop}
For any $w\in W_{\mu_n}(\mathbb R_+^n)$, the inequality \eqref{I-Ot-inq} holds.
\end{proposition}
\begin{proof}
\noindent \textbf{Step 1.} Show that the inequality \eqref{I-Ot-inq} holds for any $w\in C_c^\infty(\overline{\mathbb R_+^n})$. 

Set $u_{\ast}(x',t):=(\frac{1}{(1+t)^2+|x'|^2})^{\frac{n-p}{2(p-1)}}$, an extremal function of the Sobolev trace inequality \eqref{O-St-inq},  and $\delta:=\frac{n-p}{p(n-1)},\, h:=u_{\ast}(1+\delta w)$. It yields
\begin{equation}\label{O-ast-eq}
\big(\int_{\partial\mathbb R_+^n}|u_{\ast}|^{\frac{p(n-1)}{n-p}}\mathrm dx'\big)^{\frac{n-p}{p(n-1)}}=\frac{1}{S(n,p)}\big(\int_{\mathbb R_+^n}|\nabla u_{\ast}|^p\mathrm dx\big)^{\frac{1}{p}},
\end{equation}
and
\begin{equation}\label{O-ST-EL}
\begin{cases}
\Delta_p u_{\ast}=0,&\text{in }\mathbb R_+^n,\\
|\nabla u_{\ast}|^{p-2}\frac{\partial u_{\ast}}{\partial \mathbf{n}}=(\frac{n-p}{p-1})^{p-1}u_{\ast}^{p_{\ast}-1},&\text{on }  \partial\mathbb R_+^n,\\ 
\end{cases}    
\end{equation}
where $\mathbf{n}=(0,\cdots,0,-1)\in\mathbb R^n$ and $p_{\ast}=\frac{(n-1)p}{n-p}$. Since $h\in\mathcal D^{1,p}(\mathbb R_+^n)$, by the Sobolev trace inequality \eqref{O-St-inq}, we have
\begin{equation*}\label{O-hSt-inq}
\big(\int_{\partial\mathbb R_+^n}|h|^{\frac{p(n-1)}{n-p}}\mathrm dx'\big)^{\frac{n-p}{p(n-1)}}\leq\frac{1}{S(n,p)}\big(\int_{\mathbb R_+^n}|\nabla h|^p\mathrm dx\big)^{\frac{1}{p}}.
\end{equation*}
Combining the above with \eqref{O-ast-eq}, it holds
\begin{equation}\label{O-hu-inq}
\log\big(\frac{\int_{\partial\mathbb R_+^n}|h|^{\frac{p(n-1)}{n-p}}\mathrm dx'}{\int_{\partial\mathbb R_+^n}|u_{\ast}|^{\frac{p(n-1)}{n-p}}\mathrm dx'}\big)\leq\frac{1}{p}\log\big(\frac{\int_{\mathbb R_+^n}|\nabla h|^p\mathrm dx}{\int_{\mathbb R_+^n}|\nabla u_{\ast}|^p\mathrm dx}\big)^{\frac{1}{\delta}}.
\end{equation}
By a straightforward computation, one has
\begin{equation}\label{O-bound-inte}
\int_{\partial\mathbb R_+^n}|u_{\ast}(x',0)|^{\frac{p(n-1)}{n-p}}\mathrm dx'=\frac{1}{2}\sigma_{n-2}B\big(\frac{n-1}{2},\frac{n-1}{2(p-1)}\big)=:C_{0,p}
\end{equation}
and 
\begin{eqnarray}\label{O-grad-inte}
\int_{\mathbb R_+^n}|\nabla u_{\ast}|^p\mathrm dx=\frac{1}{2}\big(\frac{n-p}{p-1}\big)^{p-1}\sigma_{n-2}B\big(\frac{n-1}{2},\frac{n-1}{2(p-1)}\big)=\delta^{p-1}C_{1,p},
\end{eqnarray}
where $C_{1,p}=\frac{1}{2}(\frac{p(n-1)}{p-1})^{p-1}\sigma_{n-2}B(\frac{n-1}{2},\frac{n-1}{2(p-1)})$. 

Denote 
$$C_0:=\lim\limits_{p\rightarrow n}C_{0,p}=\frac{1}{2}\sigma_{n-1}\quad\text{and}\quad C_{1}:=\lim\limits_{p\rightarrow n}C_{1,p}=\frac{1}{2}n^{n-1}\sigma_{n-1}.$$ Also notice that 
\begin{eqnarray}\label{O-grad-ast}
\nabla u_{\ast}&=&-\frac{n-p}{p-1}\big[(1+t)^2+|x'|^2\big]^{-\frac{n-p}{2(p-1)}-1}(x',1+t)\nonumber\\
&=&-\delta\frac{p(n-1)}{p-1}\big[(1+t)^2+|x'|^2\big]^{-\frac{n-p}{2(p-1)}-1}(x',1+t).
\end{eqnarray}
Hence,
\begin{eqnarray}\label{O-grad-p}
|\nabla u_{\ast}|^p&=&\big(\frac{n-p}{p-1}\big)^p\big[(1+t)^2+|x'|^2\big]^{-\frac{p(n-1)}{2(p-1)}}\nonumber\\
&=&\delta^p\big(\frac{p(n-1)}{p-1}\big)^p[(1+t)^2+|x'|^2]^{-\frac{p(n-1)}{2(p-1)}}.
\end{eqnarray}
Recalling the definition of $h$, we deduce that
\begin{eqnarray*}
\int_{\partial\mathbb R_+^n}h^{\frac{p(n-1)}{n-p}}\mathrm dx' = \int_{\partial\mathbb R_+^n}\frac{1}{(1+|x'|^2)^{\frac{p(n-1)}{2(p-1)}}}(1+\delta w)^{\frac{1}{\delta}}\mathrm dx'. 
\end{eqnarray*}
Letting $p\rightarrow n$, we have $\delta\rightarrow 0$ and 
\begin{equation*}
\int_{\partial\mathbb R_+^n}h^{\frac{p(n-1)}{n-p}}\mathrm dx'\rightarrow C_0\int_{\partial\mathbb R_+^n}e^{w}\mathrm d\mu_n(x')
\end{equation*}
with
$\mathrm d\mu_n(x')
=\frac{2\mathrm dx}{\sigma_{n-1}(1+|x'|^2)^{\frac{n}{2}}}.$
We arrive at
\begin{equation}\label{O-boundary-limit}
\log\big(\frac{\int_{\partial\mathbb R_+^n}|h|^{\frac{p(n-1)}{n-p}}\mathrm dx'}{\int_{\partial\mathbb R_+^n}|u_{\ast}|^{\frac{p(n-1)}{n-p}}\mathrm dx'}\big)\rightarrow\log\big(\int_{\partial\mathbb R_+^n}e^{w}\mathrm d\mu_n(x')\big),\qquad\text{as }p\rightarrow n.
\end{equation}
Define 
$$X_{\delta}:=\nabla u_{\ast}(1+\delta w)\quad \text{and}\quad Y_{\delta}:=\delta u_{\ast}\nabla w.$$
We have
$$\nabla h=\nabla u_{\ast}(1+\delta w)+\delta u_{\ast}\nabla w=X_{\delta}+Y_{\delta}.$$ 
By the definition of $R_{p}$, one has
\begin{eqnarray}\label{O-grad-h}
|\nabla h|^p &=& |X_{\delta}+Y_{\delta}|^p = |X_\delta|^p + p |X_\delta|^{p-2}X_\delta\cdot Y_{\delta}+R_p(X_\delta, Y_{\delta})\nonumber\\
&=& |\nabla u_{\ast}|^p (1+\delta w)^p + p |\nabla u_{\ast}|^{p-2}(1+\delta w)^{p-1}\nabla u_{\ast}\cdot \delta u_{\ast}\nabla w+ R_p(X_\delta, Y_{\delta})\nonumber\\
&=&|\nabla u_{\ast}|^p (1+\delta w)^p + u_{\ast} |\nabla u_{\ast}|^{p-2}\nabla u_{\ast} \cdot \nabla (1+\delta w)^p+ R_p(X_\delta,Y_\delta).
\end{eqnarray}
Using integration by parts together with equations \eqref{O-ST-EL}, \eqref{O-bound-inte} and \eqref{O-grad-inte}, we obtain
\begin{eqnarray}\label{O-2nd-int}
&&\int_{\mathbb R_+^n} u_{\ast} |\nabla u_{\ast}|^{p-2} \nabla u_{\ast}\cdot\nabla(1+\delta w)^p\mathrm dx \nonumber\\
&=& -\int_{\mathbb R_+^n} |\nabla u_{\ast}|^p (1+\delta w)^p\mathrm dx-\int_{\mathbb R_+^n} u_{\ast}\Delta_p u_{\ast} (1+\delta w)^p\mathrm dx\nonumber\\
&&+\int_{\partial\mathbb R_+^n} u_{\ast} |\nabla u_{\ast}|^{p-2} \frac{\partial u_{\ast}}{\partial \mathbf{n}} (1+\delta w)^p\mathrm dx'\nonumber\\
&=&-\int_{\mathbb R_+^n} |\nabla u_{\ast}|^p (1+\delta w)^p\mathrm dx+\int_{\partial\mathbb R_+^n} u_{\ast} |\nabla u_{\ast}|^{p-2} \frac{\partial u_{\ast}}{\partial \mathbf{n}} (1+\delta w)^p\mathrm dx'\nonumber\\
&=&-\int_{\mathbb R_+^n} |\nabla u_{\ast}|^p (1+\delta w)^p\mathrm dx+\int_{\mathbb R_+^n}|\nabla u_{\ast}|^p\mathrm dx\nonumber\\
&&+\big(\frac{p(n-1)}{p-1}\big)^{p-1}\delta^{p-1}\int_{\partial\mathbb R_+^n} u_{\ast}^{p_{\ast}}\big(p\delta w+o(\delta)\big)\mathrm dx',\qquad\text{as }p\rightarrow n.
\end{eqnarray}

Combining \eqref{O-grad-h} and \eqref{O-2nd-int} yields
\begin{eqnarray}\label{O-h-int}
\int_{\mathbb R_+^n} |\nabla h|^p \mathrm dx
&=&\int_{\mathbb R_+^n}|\nabla u_{\ast}|^p\mathrm dx+\big(\frac{p(n-1)}{p-1}\big)^{p-1}\delta^{p-1}\int_{\partial\mathbb R_+^n} u_{\ast}^{p_{\ast}}\big(p\delta w+o(\delta)\big)\mathrm dx'\nonumber\\
&&+\int_{\mathbb R_+^n}R_p(X_\delta,Y_\delta)\mathrm dx,\qquad\text{as }p\rightarrow n.
\end{eqnarray}

Note that 
\begin{eqnarray}\label{O-boundary-int}
&&\big(\frac{p(n-1)}{p-1}\big)^{p-1}\delta^{p-1}\int_{\partial\mathbb R_+^n} u_{\ast}^{p_{\ast}}\big(p\delta w+o(\delta)\big)\mathrm dx'\nonumber\\
&=& \big(\frac{p(n-1)}{p-1}\big)^{p-1}\delta^pp\int_{\partial\mathbb R_+^n}\frac{w\mathrm dx'}{(1+|x'|^2)^{\frac{p(n-1)}{2(p-1)}}}+o(\delta^p),\quad\text{as }p\rightarrow n
\end{eqnarray}
and, by homogeneity,
\begin{eqnarray}\label{O-Rp-limits}
\frac{1}{\delta^p} R_p(X_\delta, Y_{\delta}) = R_p\big(\frac{1}{\delta}X_\delta, \frac{1}{\delta}Y_{\delta}\big)\rightarrow R_n\big(-\frac{n(x', 1+t)}{(1+t)^2 + |x'|^2} , \nabla w\big),\text{ as }p\rightarrow n.
\end{eqnarray}

Now, recalling \eqref{O-grad-inte} and substituting \eqref{O-boundary-int} and \eqref{O-Rp-limits} into \eqref{O-h-int}, it holds
\begin{equation}\label{O-inter-limit}
\frac{1}{p}\log\big(\frac{\int_{\mathbb R_+^n}|\nabla h|^p\mathrm dx}{\int_{\mathbb R_+^n}|\nabla u_{\ast}|^p\mathrm dx}\big)^{\frac{1}{\delta}}=\frac{1}{p}\log\big(1+\delta A_{\delta} +o(\delta)\big)^{\frac{1}{\delta}}\rightarrow\frac{1}{n}A,\qquad\text{as }p\rightarrow n,
\end{equation}
where 
\begin{eqnarray*}
A_{\delta}&=&\big(\frac{p(n-1)}{p-1}\big)^{p-1}p\frac{1}{C_{1,p}}\int_{\partial\mathbb R_+^n}\frac{1}{(1+|x'|^2)^{\frac{p(n-1)}{2(p-1)}}}w\mathrm dx'\nonumber\\
&&+\frac{1}{C_{1,p}}\int_{\mathbb R_+^n}R_p\big(\frac{1}{\delta}X_\delta,\frac{1}{\delta}Y_{\delta}\big)\mathrm dx\nonumber\\ 
&\rightarrow&\frac{\sigma_{n-1} n^n}{2}\frac{1}{C_{1}}\int_{\partial\mathbb R_+^n}w\mathrm d\mu_n(x')+\frac{1}{C_{1}}\int_{\mathbb R_+^n}R_n\big(-\frac{n(x', 1+t)}{(1+t)^2 + |x'|^2}, \nabla w\big)\mathrm dx\nonumber\\
&=&n\int_{\partial\mathbb R_+^n}w\mathrm d\mu_n(x')+\frac{2}{n^{n-1}\sigma_{n-1}}\int_{\mathbb R_+^n}K_n(x,\nabla w)\mathrm dx=:A
\end{eqnarray*}
as $p\rightarrow n$. It follows from \eqref{O-hu-inq}, \eqref{O-boundary-limit} and \eqref{O-inter-limit} that
\begin{eqnarray*}
\log\big(\int_{\partial\mathbb R_+^n}e^{w}\mathrm d\mu_n(x')\big)-\int_{\partial\mathbb R_+^n}w\mathrm d\mu_n(x') \leq\frac{2}{n^n\sigma_{n-1}}\int_{\mathbb R_+^n}K_n(x,\nabla w)\mathrm dx.
\end{eqnarray*}

\noindent\textbf{Step 2.}
Extend the inequality \eqref{I-Ot-inq} to 
$W_{\mu_n}(\mathbb R_+^n)$ by the density argument. 

For any $w\in W_{\mu_n}(\mathbb R_+^n)$, we know by Proposition \ref{density-prop} that there exists a sequence $\{w_k\}_{k\in \mathbb N}\subset C_c^\infty(\overline{\mathbb R_+^n})$ such that $\|w-w_k\|_{\mu_n}\rightarrow 0$. 
Then, there exists a subsequence (still denoted by $\{w_k\}_{k\in\mathbb N}$) and a nonnegative function $f\in L^n(\mathbb R_+^n,\mathrm dx)$ with $f^2|\nabla \log\mu_n|^{n-2}\in L^1(\mathbb R_+^n,\mathrm dx)$ such that

$(i)$ $w_k(x)\rightarrow w(x)$ a.e. in $\partial\mathbb R_+^n$ as $k\rightarrow+\infty$;

$(ii)$ $\nabla w_k(x)\rightarrow \nabla w(x)$ a.e. in $\mathbb R_+^n$ as $k\rightarrow+\infty$;

$(iii)$ for any $k$, $|\nabla w_k(x)|\leq f(x)$ a.e. in $\mathbb R_+^n$;

$(iv)$ $K_n\left(x,\nabla w_k\right)\rightarrow K_n\left(x,\nabla w\right)$ a.e. in $\mathbb R_+^n$ as $k\rightarrow+\infty$.\newline 
Recalling the fact that $0\leq R_n(X,Y)\leq C(n)(|Y|^n+|Y|^2|X|^{n-2})$ for any $X,Y\in\mathbb R^n$ (see Lemma 2.1 in \cite{BCM2025}) and the definition of $K_n$, for any $k\in\mathbb N$, one has
\begin{eqnarray*}
0\leq K_n\left(x,\nabla w_k(x)\right)&\leq& C\left(|\nabla w_k(x)|^n+|\nabla w_k(x)|^2|\nabla \log\mu_n(x)|^{n-2}\right)\nonumber\\
&\leq&C\left(f^n(x)+f^2(x)|\nabla\log\mu_n(x)|^{n-2}\right)\quad a.e.\text{ in }\mathbb R_+^n.
\end{eqnarray*}
By the dominated convergence theorem, it yields 
\begin{equation*}
\int_{\mathbb R_+^n}K_n\left(x,\nabla w_k\right)\mathrm dx\rightarrow\int_{\mathbb R_+^n}K_n\left(x,\nabla w\right)\mathrm dx,\quad\text{as }k\rightarrow\infty
\end{equation*}
and from Fatou's Lemma one has
$$\int_{\partial\mathbb R_+^n} e^w\mathrm d\mu_n(x') \leq\liminf\limits_{k\rightarrow+\infty}\big(\int_{\partial\mathbb R_+^n}e^{w_k}\mathrm d\mu_n(x')\big).$$
Hence,
\begin{eqnarray*}
\log\int_{\partial\mathbb R_+^n} e^w\mathrm d\mu_n(x') &\leq&\liminf\limits_{k\rightarrow+\infty}\big(\log\int_{\partial\mathbb R_+^n}e^{w_k}\mathrm d\mu_n(x')\big)\\
&\leq &\liminf\limits_{k\rightarrow+\infty} \big(\alpha_n\int_{\mathbb R_+^n}K_n(x, \nabla w_k)\mathrm dx+\int_{\partial\mathbb R_+^n}w_k\mathrm d\mu_n(x')\big)\\
&=& \alpha_n\int_{\mathbb R_+^n}K_n(x, \nabla w)\mathrm dx+\int_{\partial\mathbb R_+^n}w\mathrm d\mu_n(x').
\end{eqnarray*}
The proposition is proved.
\end{proof}

\subsection{Optimality of Onofri trace inequality}
In the rest of this section, we shall show that $\alpha_n$ is the best constant and we also present the sufficiency of extremal functions of the inequality \eqref{I-Ot-inq}.

Define
$$Q_n[w]:=\frac{\int_{\mathbb R_+^n}K_n(x,\nabla w)\mathrm dx}{\log\big(\int_{\partial\mathbb R_+^n}e^{w}\mathrm d\mu_n(x')\big)- \int_{\partial\mathbb R_+^n}w\mathrm d\mu_n(x')}.$$
To show the optimality of $\alpha_n$, 
it is sufficient to prove the existence of some $w_0\in W_{\mu_n}(\mathbb R_+^n)$ 
such that $\frac1{\alpha_n}=Q_n[w_0]$.
  
\begin{proposition}\label{Best-const}The equality holds in inequality \eqref{I-Ot-inq} for function $w$ satisfying the form \eqref{Extre-f}. Moreover, the constant $\alpha_n=\frac{2}{n^n\sigma_{n-1}}$ in \eqref{I-Ot-inq} is optimal.
\end{proposition}
\begin{proof}
Let $w$ be the function of form \eqref{Extre-f}. Firstly, we check that $w$ lies in the space $W_{\mu_n}(\mathbb R_+^n)$. Indeed, 
by the asymptotic behavior of $w$ as $|x|\rightarrow +\infty$, there exist constants $C>0$ depending only on $\lambda,x_0',n,\tilde C$ and $M>0$ such that
$$\int_{\partial\mathbb R_+^n}|w|\mathrm d\mu_n(x')\leq C\int_{\mathbb R^{n-1}}\left[\log(|x'|^2+1)+1\right]\mathrm d\mu_n(x')<+\infty,$$
and 
\begin{eqnarray*}
\int_{\mathbb R_+^n}|\nabla w|^n\mathrm dx&\leq& C\int_{\mathbb R_+^n\cap\{|x|\leq M\}}\frac{\mathrm dx}{\left(|x'|^2+(1+t)^2\right)^{\frac n 2}}+C\int_{\mathbb R_+^n\cap\{|x|> M\}}\frac{\mathrm dx}{|x|^{2n}}\nonumber\\
 &<&+\infty,   
\end{eqnarray*}
\begin{eqnarray*}
\int_{\mathbb R_+^n}|\nabla w|^2|\nabla\log\mu_n|^{n-2}\mathrm dx&\leq& C\int_{\mathbb R_+^n\cap\{|x|\leq M\}}\frac{\mathrm dx}{\left(|x'|^2+(1+t)^2\right)^{\frac n 2}}\nonumber\\
&&+C\int_{\mathbb R_+^n\cap\{|x|> M\}}\frac{\mathrm dx}{|x|^{n+2}}<+\infty,    
\end{eqnarray*}
by noting that
$$|\nabla w|=O(|x|^{-2})\text{ and }|\nabla\log\mu_n|=O(|x|^{-1}),\quad\text{as }|x|\rightarrow+\infty.$$
This proves $w\in W_{\mu_n}(\mathbb R_+^n)$. 

It is easy to observe that if $w$ is an extremal function, then $w+C$ is also one ($C$ maybe any constant). So, without loss of generality, we take $\tilde C=\log\frac{n^{n-1}\sigma_{n-1}}{2}$ in \eqref{Extre-f}. 
That is,
\begin{equation}\label{T-EL-solution2}
w(x)=\log\frac{n^{n-1}\sigma_{n-1}\big(|x'|^2+(1+t)^2\big)^{\frac{n}{2}}\lambda}{2\big(|x'-x_0'|^2+(t+\lambda)^2\big)^{\frac{n}{2}}}=u(x)-\log\mu_n(x),
\end{equation}
where $u(x)=\log\frac{n^{n-1}\lambda}{\big(|x'-x_0'|^2+(t+\lambda)^2\big)^{\frac{n}{2}}}$ satisfies \eqref{I-hL-eq} and $\mu_n(x)=\frac{2}{\sigma_{n-1}(|x'|^2+(1+t)^2)^{\frac{n}{2}}}$ satisfies
\begin{equation}\label{T-mu-eq}
\begin{cases}
-\Delta_n(\log\mu_n)=0,&\qquad\text{in }\mathbb R_+^n,\\
a(\nabla\log\mu_n)\cdot\mathbf{n}=\frac{n^{n-1}\sigma_{n-1}}{2}\mu_n,&\qquad\text{on }\partial\mathbb R_+^n.
\end{cases}    
\end{equation}

Next, we show $\frac1{\alpha_n}=Q_n[w]$.

Inserting $\eqref{T-EL-solution2}$ into \eqref{I-Ot-inq} yields that all terms are as follows 
\begin{equation}\label{T-1-term}
\log\int_{\partial\mathbb R_+^n}e^w\mathrm d\mu_n(x')=\log\int_{\partial\mathbb R_+^n}e^u\mathrm dx'=\log\frac{n^{n-1}\sigma_{n-1}}{2},
\end{equation}
\begin{equation}\label{T-2-term}
\int_{\partial\mathbb R_+^n}w\mathrm d\mu_n(x')=\int_{\partial\mathbb R_+^n}u\mu_n\mathrm dx'-\int_{\partial\mathbb R_+^n}\mu_n\log\mu_n\mathrm dx',
\end{equation}
and
\begin{eqnarray}\label{T-3-term}
&&\int_{\mathbb R_+^n}K_n(x,\nabla w)\mathrm dx\nonumber\\
&=&\int_{\mathbb R_+^n}\big[|\nabla u|^n-|\nabla\log\mu_n|^n-n|\nabla\log\mu_n|^{n-2}\nabla\log\mu_n\cdot(\nabla u-\nabla\log\mu_n)\big]\mathrm dx.\nonumber\\
&&
\end{eqnarray}
By using integration by parts, \eqref{I-hL-eq} and \eqref{T-mu-eq}, we obtain
\begin{eqnarray}\label{T-3-R}
&&\int_{B_R^+}\big[|\nabla u|^n-|\nabla\log\mu_n|^n-n|\nabla\log\mu_n|^{n-2}\nabla\log\mu_n\cdot(\nabla u-\nabla\log\mu_n)\big]\mathrm dx\nonumber\\
&=&\int_{\partial B_R^+}\big[ua(\nabla u)-\log\mu_n a(\nabla\log\mu_n)-n(u-\log\mu_n)a(\nabla\log\mu_n) \big]\cdot\mathbf{n}\mathrm d\mathcal{H}^{n-1}\nonumber\\
&&+\int_{\Sigma_R}e^u u\mathrm dx'+\frac{n^{n-1}(n-1)\sigma_{n-1}}{2}\int_{\Sigma_R}\mu_n\log\mu_n\mathrm dx'\nonumber\\
&&-\frac{n^n\sigma_{n-1}}{2}\int_{\Sigma_R}\mu_n u\mathrm dx',\qquad \forall R>0.
\end{eqnarray}
From the explicit form of $u$ and $\mu_n$, one has
\begin{eqnarray}\label{T-3-R1}
&&\int_{\partial B_R^+}\big[ua(\nabla u)-\log\mu_n a(\nabla\log\mu_n)-n(u-\log\mu_n)a(\nabla\log\mu_n) \big]\cdot\mathbf{n}\mathrm d\mathcal{H}^{n-1}\nonumber\\
&&\rightarrow\frac{n^{n-1}(n-1)\sigma_{n-1}}{2}\log\frac{n^{n-1}\lambda\sigma_{n-1}}{2},\qquad\text{as }R\rightarrow\infty,   
\end{eqnarray}
by 
$$a(\nabla u)\cdot\mathbf{n}=-n^{n-1}R^{1-n}(1+O(R^{-1})),$$
$$a(\nabla \log\mu_n)\cdot\mathbf{n}=-n^{n-1}R^{1-n}(1+O(R^{-1})),$$
$$u-\log\mu_n=\log\frac{n^{n-1}\lambda\sigma_{n-1}}{2}+O(R^{-1})$$
on $\partial B_R^+$ as $R\rightarrow\infty$. Combining \eqref{T-3-term}, \eqref{T-3-R} and \eqref{T-3-R1} and letting $R\rightarrow\infty$, we have
\begin{eqnarray}\label{T-3-term1}
 \frac{2}{n^n\sigma_{n-1}}\int_{\mathbb R_+^n}K_n(x,\nabla w)\mathrm dx 
&=&\frac{n-1}{n}\log\frac{n^{n-1}\lambda\sigma_{n-1}}{2}+\frac{2}{n^n\sigma_{n-1}}\int_{\partial\mathbb R_+^n}e^u u\mathrm dx'\nonumber\\
&&+\frac{n-1}{n}\int_{\partial\mathbb R_+^n}\mu_n\log\mu_n\mathrm dx'-\int_{\partial\mathbb R_+^n}\mu_n u\mathrm dx'.
\end{eqnarray}
It follows from \eqref{T-1-term}, \eqref{T-2-term} and \eqref{T-3-term1} that the equality holds in \eqref{I-Ot-inq} if 
\begin{eqnarray*}
&&\log\frac{n^{n-1}\sigma_{n-1}}{2}\nonumber\\
&=&\frac{n-1}{n}\log\frac{n^{n-1}\lambda\sigma_{n-1}}{2}+\frac{2}{n^n\sigma_{n-1}}\int_{\partial\mathbb R_+^n}e^u u\mathrm dx'-\frac{1}{n}\int_{\partial\mathbb R_+^n}\mu_n\log\mu_n\mathrm dx',
\end{eqnarray*}
that is,
\begin{eqnarray*}
&&\frac{\sigma_{n-1}}{2}\log\frac{n^{n-1}\lambda^{1-n}\sigma_{n-1}}{2}\nonumber\\
&=&\int_{\partial\mathbb R_+^n}\frac{\lambda}{(|x'-x_0'|^2+\lambda^2)^\frac{n}{2}}\log\frac{n^{n-1}\lambda}{(|x'-x_0'|^2+\lambda^2)^\frac{n}{2}}\mathrm dx'\nonumber\\
&&-\int_{\partial\mathbb R_+^n}\frac{1}{(|x'|^2+1)^\frac{n}{2}}\log\frac{2}{\sigma_{n-1}(|x'|^2+1)^\frac{n}{2}}\mathrm dx'\nonumber\\
&=&\int_{\mathbb R^{n-1}}\frac{1}{(|y|^2+1)^\frac{n}{2}}\big[\log\frac{n^{n-1}\lambda^{1-n}}{(|y|^2+1)^\frac{n}{2}}-\log\frac{2}{\sigma_{n-1}(|y|^2+1)^\frac{n}{2}}\big]\mathrm dy\nonumber\\
&=&\log\frac{n^{n-1}\lambda^{1-n}\sigma_{n-1}}{2}\int_{\mathbb R^{n-1}}\frac{1}{(|y|^2+1)^\frac{n}{2}}\mathrm dy\nonumber\\
&=&\frac{\sigma_{n-1}}{2}\log\frac{n^{n-1}\lambda^{1-n}\sigma_{n-1}}{2},
\end{eqnarray*}
which implies that for any $x_0'\in\mathbb R^{n-1}$, $\lambda>0$ and $\tilde C>0$, the form \eqref{Extre-f} of $w$ can make the equality hold in inequality \eqref{I-Ot-inq}. Also note that for any $\lambda\neq1$ and $x_0'=0$, we have $$\int_{\mathbb R_+^n}K_n(x,\nabla w)\mathrm dx\neq 0\quad
\text{and}\quad\log\int_{\partial\mathbb R_+^n}e^w\mathrm d\mu_n(x')-\int_{\partial\mathbb R_+^n}w\mathrm d\mu_n(x')\neq 0.$$  
It implies that $Q_n[w]=\frac{1}{\alpha_n}$, and hence $\inf\limits_{w\in W_{\mu_n}(\mathbb R_+^n)}Q_n[w]=\frac1{\alpha_n}$.

The above process indicates that $\alpha_n$ is the best constant 
and $w\in W_{\mu_n}(\mathbb R_+^n)$ is the extremal function of inequality \eqref{I-Ot-inq} when $w$ is of the form \eqref{Extre-f}.
\end{proof}

\section{\textbf{Regularity}}\label{Section 3}

In this section, we show the regularity of the solution to \eqref{I-hL-eq}.

Consider the following equation
\begin{equation}\label{R-lb-eq}
\begin{cases}
-\Delta_n u=f, &\quad\text{in}~B^+_R(y),\\
|\nabla u|^{n-2}\frac{\partial u}{\partial t}=-g, &\quad\text{on}~\Sigma_R(y),
\end{cases} 
\end{equation}
where $f\in L^p(B^+_R(y)),~g\in L^p(\Sigma_R(y))$, $p>1$ and $y\in\overline{\mathbb R_+^n}$. 

\begin{definition}\label{R-lb-weak}
We say $u\in W^{1,n}( B^+_R(y))$ is a weak solution to \eqref{R-lb-eq} if 
\begin{equation}\label{R-lb-weak-eq}
\int_{B^+_R(y)}|\nabla u|^{n-2}\langle\nabla u,\nabla\varphi\rangle \mathrm dx= \int_{\Sigma_R(y)}g\varphi \mathrm dx'+
\int_{B^+_R(y)}f\varphi \mathrm dx
\end{equation}
for any $\varphi\in C^{\infty}_c(B_R(y)\cap\overline{\mathbb R_+^n})$.
\end{definition}

The following two lemmas on the local boundedness of the solution and the Brezis-Merle-type exponential inequality, respectively, are essential for the proof of the global boundedness from above of the solution to \eqref{I-hL-eq}.

\begin{lemma}\label{loc-bdd-lem}
Suppose $u\in W^{1,n}( B^+_R(y))$ is a weak solution to \eqref{R-lb-eq}. Then we have 
\begin{equation*}\label{R-lb-inq}
\|u^{+}\|_{L^{\infty}(B^+_{\frac{R}{2}}(y)})\leq C(\|u^{+}\|_{L^n(B^+_R(y))}+\|f\|_{L^p(B^+_R(y))}+\|g\|_{L^p(\Sigma_R(y))})
\end{equation*}
where $C$ depends only on $R,n,p$.
\end{lemma}

\begin{proof}
We employ the standard Moser iteration technique to establish the desired result. For some $k,~m>0$, set $\overline{u}=u^{+}+k$ and
\begin{equation*}
\overline{u}_{m}=
\begin{cases}
\overline{u}, &\quad \text{if }u<m,\\
k+m, &\quad \text{if }u\geq m.
\end{cases} 
\end{equation*}
For any $y\in\overline{\mathbb R_+^n}, R>0$, we have
\begin{equation*}
\overline{u},\overline{u}_{m}\in W^{1,n}(B^+_R(y)),~k\leq\overline{u}_{m}\leq k+m, \overline{u}_{m}\leq\overline{u},~\nabla \overline{u}_{m}=
\begin{cases}
\nabla u, & \text{if }0\leq u\leq m,\\
0,& \text{otherwise}.
\end{cases}
\end{equation*}
Let $\eta\in C_c^{\infty}(B_R)$, $\eta\geq0$, $\beta>0$ and  $\varphi=\eta^n(\overline{u}_{m}^{\beta}\overline{u}-k^{\beta+1})$, a direct computation yields
\begin{eqnarray*}\label{R-phi-eq}
\nabla \varphi &=& \beta\eta^n\overline{u}_{m}^{\beta-1}\nabla \overline{u}_{m}\overline{u}+ \eta^n\overline{u}_{m}^{\beta}\nabla \overline{u}+n\eta^{n-1}\nabla\eta(\overline{u}_{m}^{\beta}\overline{u}-k^{\beta+1})\nonumber\\
&=& \eta^n\overline{u}_{m}^{\beta}(\beta\nabla \overline{u}_{m}+\nabla \overline{u})+n\eta^{n-1}\nabla\eta(\overline{u}_{m}^{\beta}\overline{u}-k^{\beta+1}),
\end{eqnarray*}
and $\varphi=0$ and $\nabla \varphi=0$ in $\{u\leq0\}$.
 
By \eqref{R-lb-weak-eq}, we have
\begin{eqnarray*}
&&\beta\int_{B^+_R(y)}\eta^n|\nabla \overline{u}_{m}|^n\overline{u}_{m}^{\beta}\mathrm dx+\int_{B^+_R(y)}\eta^n|\nabla \overline{u}|^n\overline{u}_{m}^{\beta}\mathrm dx \nonumber\\
&=&-n\int_{B^+_R(y)}|\nabla \overline{u} |^{n-2}\nabla\overline{u}\cdot\nabla\eta\eta^{n-1}(\overline{u}_{m}^{\beta}\overline{u}-k^{\beta+1})\mathrm dx\nonumber\\&&+\int_{\Sigma_R(y)}g\eta^n(\overline{u}_{m}^{\beta}\overline{u}-k^{\beta+1})\mathrm dx'+\int_{B^+_R(y)}f\eta^n(\overline{u}_{m}^{\beta}\overline{u}-k^{\beta+1})\mathrm dx.
\end{eqnarray*}
By the definition of $\bar{u}$, it holds
\begin{eqnarray}\label{R-phi-inq}
& &\beta\int_{B^+_R(y)}\eta^n|\nabla \overline{u}_{m}|^n\overline{u}_{m}^{\beta}\mathrm dx+\int_{B^+_R(y)}\eta^n|\nabla \overline{u}|^n\overline{u}_{m}^{\beta}\mathrm dx \nonumber\\
&\leq&C(n)\big[\int_{B^+_R(y)}|\nabla \overline{u}|^{n-1}|\nabla\eta|\eta^{n-1}(\overline{u}_{m}^{\beta}\overline{u}-k^{\beta+1})\mathrm dx\nonumber\\&&+\int_{\Sigma_R(y)}|g|\eta^n(\overline{u}_{m}^{\beta}\overline{u}-k^{\beta+1})\mathrm dx'+\int_{B^+_R(y)}|f|\eta^n(\overline{u}_{m}^{\beta}\overline{u}-k^{\beta+1})\mathrm dx\big].
\end{eqnarray}
By Young's inequality, we have 
\begin{equation}\label{R-Young-inq}
|\nabla \overline{u}|^{n-1}|\nabla \eta|\eta^{n-1}\overline{u}_{m}^{\beta}\overline{u}\leq\varepsilon|\nabla \overline{u}|^n\eta^n\overline{u}_{m}^{\beta}+C(\varepsilon)|\nabla \eta|^n\overline{u}_{m}^{\beta}\overline{u}^n.
\end{equation}
For sufficiently small $\varepsilon>0$, by \eqref{R-phi-inq}, \eqref{R-Young-inq} and $\overline{u},~\overline{u}_{m}\geq k$, we obtain
\begin{eqnarray}\label{R-eta-inq}
&&\beta\int_{B^+_R(y)}\eta^n|\nabla \overline{u}_{m}|^n\overline{u}_{m}^{\beta}\mathrm dx+\int_{B^+_R(y)}\eta^n|\nabla \overline{u}|^n\overline{u}_{m}^{\beta}\mathrm dx \nonumber\\
&\leq& C(n)\big(\int_{B^+_R(y)}|\nabla \eta|^n\overline{u}_{m}^{\beta}\overline{u}^n\mathrm dx+\int_{\Sigma_R(y)}|g|\eta^n\overline{u}_{m}^{\beta}\overline{u}\mathrm dx'\nonumber\\
&&+\int_{B^+_R(y)}|f|\eta^n\overline{u}_{m}^{\beta}\overline{u}\mathrm dx\big).
\end{eqnarray}
Set $w=\overline{u}_{m}^{\frac{\beta}{n}}\overline{u}$. It is easy to see
\begin{eqnarray}\label{R-w-inq}
|\nabla (\eta w)|^n&\leq&C(n)(|\nabla \eta|^nw^n+\eta^n|\nabla w|^n)\nonumber\\
&\leq&C(n)\big[|\nabla \eta|^n\overline{u}_{m}^{\beta}\overline{u}^n+(1+\beta)^{n-1}
(\beta\eta^n\overline{u}_{m}^{\beta}|\nabla\overline{u}_{m}|^n+\eta^n\overline{u}_{m}^{\beta}|\nabla\overline{u}|^n)\big].\nonumber\\
&&
\end{eqnarray}
From \eqref{R-eta-inq}, \eqref{R-w-inq}, for $\overline{u}\geq k$, it yields 
\begin{eqnarray}\label{R-w-int}
&&\int_{B^+_R(y)}|\nabla (\eta w)|^n\mathrm dx \nonumber\\
&\leq& C(n)(1+\beta)^{n-1}\big(\int_{B^+_R(y)}|\nabla \eta|^n\overline{u}_{m}^{\beta}\overline{u}^n\mathrm dx+\int_{\Sigma_R(y)}|g|\eta^n\overline{u}_{m}^{\beta}\overline{u}\mathrm dx' \nonumber\\
&&+\int_{B^+_R(y)}|f|\eta^n\overline{u}_{m}^{\beta}\overline{u}\mathrm dx\big)\nonumber\\
&\leq&C(n)(1+\beta)^{n-1}\big(\int_{B^+_R(y)}|\nabla \eta|^n\overline{u}_{m}^{\beta}\overline{u}^n\mathrm dx+\int_{\Sigma_R(y)}\frac{|g|}{k^{n-1}}\eta^n\overline{u}_{m}^{\beta}{\overline{u}}^n\mathrm dx' \nonumber\\
&&+\int_{B^+_R(y)}\frac{|f|}{k^{n-1}}\eta^n\overline{u}_{m}^{\beta}{\overline{u}}^n\mathrm dx\big).
\end{eqnarray}
Choose $k^{n-1}=\|f\|_{L^p(B^+_R(y))}+\|g\|_{L^p(\Sigma_R(y))}$ if $f$ or $g$ is not identically $0$. Otherwise, choose arbitrarily $k>0$ and eventually let $k\rightarrow0^{+}$. By assumption we have 
\begin{equation*}\label{f-g-norm-inq}
\|\frac{f}{k^{n-1}}\|_{L^p(B^+_R(y))}\leq1~\text{and}~\|\frac{g}{k^{n-1}}\|_{L^p(\Sigma_R(y))}\leq1.
\end{equation*}
Hence, the H\"{o}lder inequality implies 
\begin{eqnarray}\label{R-f-int}
\int_{B^+_R(y)}\frac{|f|}{k^{n-1}}\eta^n\overline{u}_{m}^{\beta}{\overline{u}}^n\mathrm dx&=&\int_{B^+_R(y)}\frac{|f|}{k^{n-1}}\eta^nw^n\mathrm dx \nonumber\\
&& \leq(\int_{B^+_R(y)}(\frac{|f|}{k^{n-1}})^p\mathrm dx)^{\frac{1}{p}}( \int_{B^+_R(y)}(\eta w)^{\frac{np}{p-1}}\mathrm dx)^{1-\frac{1}{p}}\nonumber\\
&&\leq( \int_{B^+_R(y)}(\eta w)^{\frac{np}{p-1}}\mathrm dx)^{1-\frac{1}{p}}
\end{eqnarray}
and 
\begin{eqnarray}\label{R-g-int}
\int_{\Sigma_R(y)}\frac{|g|}{k^{n-1}}\eta^n\overline{u}_{m}^{\beta}{\overline{u}}^n\mathrm dx'  
&=& \int_{\Sigma_R(y)}\frac{|g|}{k^{n-1}}\eta^nw^n\mathrm dx' \nonumber\\
&&\leq (\int_{\Sigma_R(y)}(\frac{|g|}{k^{n-1}})^p\mathrm dx')^{\frac{1}{p}}( \int_{\Sigma_R(y)}(\eta w)^{\frac{np}{p-1}}\mathrm dx')^{1-\frac{1}{p}}\nonumber\\
&&\leq(\int_{\Sigma_R(y)}(\eta w)^{\frac{np}{p-1}}\mathrm dx')^{1-\frac{1}{p}}.
\end{eqnarray}
For fixed $q>\frac{np}{p-1}$, employing interpolation inequality and Sobolev inequality, it yields 
\begin{eqnarray}\label{R-interior-inq}
\|\eta w\|_{L^{\frac{np}{p-1}}(B^+_R(y))}&\leq&\varepsilon\|\eta w\|_{L^{q}(B^+_R(y))}+C(n,p,q)\varepsilon^{-\tau}\|\eta w\|_{L^n(B^+_R(y))}\nonumber\\
&\leq&C\varepsilon\|\nabla (\eta w)\|_{L^n(B^+_R(y))}+C(n,p,q)\varepsilon^{-\tau}\|\eta w\|_{L^n(B^+_R(y))},\nonumber\\
&&
\end{eqnarray}
where $\tau=\frac{q}{q(p-1)-np}>0$ and $C=C(n,p,q, R)$. Similarly, it follows from interpolation inequality and Sobolev trace inequality that 
\begin{eqnarray}\label{R-boundary-inq}
\|\eta w\|_{L^{\frac{np}{p-1}}(\Sigma_R(y))}&\leq&\varepsilon\|\eta w\|_{L^{q}(\Sigma_R(y))}+C(n,p,q)\varepsilon^{-\tau}\|\eta w\|_{L^n(\Sigma_R(y))}\nonumber\\
&\leq&C\varepsilon\|\nabla (\eta w)\|_{L^n(B^+_R(y))}+C(n,p,q)\varepsilon^{-\tau}\|\eta w\|_{L^n(\Sigma_R(y))},\nonumber\\
&&
\end{eqnarray} 
where $C=C(n,p,q,R)$. Combining \eqref{R-w-int} and \eqref{R-f-int}-\eqref{R-boundary-inq} yields
\begin{eqnarray}\label{R-eta-int}
&&\int_{B^+_R(y)}|\nabla (\eta w)|^n\mathrm dx \nonumber\\
&\leq&C\big[(1+\beta)^{n-1}\int_{B^+_R(y)}|\nabla \eta|^nw^n\mathrm dx+(1+\beta)^{(n-1)(\tau+1)}\int_{B^+_R(y)}\eta^n w^n\mathrm dx\nonumber\\
&+&(1+\beta)^{(n-1)(\tau+1)}\int_{\Sigma_R(y)}\eta^n w^n\mathrm dx'\big],
\end{eqnarray}
where $C$ is independent of $\beta$. Using Sobolev inequality and Sobolev trace inequality and fixing $s>n$, by \eqref{R-eta-int}, we obtain 
\begin{eqnarray}\label{R-w-inq1}
&&\|\eta w\|_{L^{s}(B^+_R(y))}+\|\eta w\|_{L^{s}(\Sigma_R(y))}\nonumber\\
&\leq& C(1+\beta)^{\tau_1}\big[\big(\int_{B^+_R(y)}(|\nabla \eta|^n+\eta^n)w^n\mathrm dx\big)^{\frac 1 n}+\big(\int_{\Sigma_R(y)}\eta^nw^n\mathrm dx'\big)^{\frac 1 n}\big],~~~~
\end{eqnarray}
where $\tau_1=\frac{(n-1)(\tau+1)}{n}$ and $C$ is a positive constant independent of $\beta$. For $0<r<R$, choose the nonnegative cut-off function $\eta$ as 
\begin{equation*}
\eta =
\begin{cases}
1,\quad&\text{in }B_r,\\
0,\quad&\text{in }\mathbb R^n\setminus B_R,\\ 
\end{cases}
\end{equation*}
which also satisfies $0\leq\eta\leq 1$ and $|\nabla\eta|\leq C(R-r)^{-1}$. From \eqref{R-w-inq1}, we infer that
\begin{equation}\label{R-w-inq2}
\|w\|_{L^{s}(B^+_r(y))}+\|w\|_{L^{s}(\Sigma_r(y))}\leq C(1+\beta)^{\tau_1}\big(1+\frac{1}{R-r}\big)\big(\|w\|_{L^n(B^+_R(y))}+\|w\|_{L^n(\Sigma_R(y))}\big),
\end{equation}
where $C$ is independent of $\beta$. Following the idea of \cite{HL2011}, by \eqref{R-w-inq2}, a classical Moser iteration argument yields the desired result.
\end{proof}

\begin{lemma}\label{BM-inq-lem}
Let $u\in W^{1,n}_{\operatorname{loc}}(\overline{\mathbb R_+^n})$ be a weak solution of \eqref{I-hL-eq} satisfying \eqref{I-fm-assu} and $h\in W^{1,n}(B^+_R(y))$ be a weak solution of 
\begin{equation}\label{R-uh-eq}
\begin{cases}
-\Delta_n h=0, &\qquad\text{in}~B^+_R(y),\\
h=u, &\qquad\text{on}~\partial B_R^+(y),\\
\langle a(\nabla h),\mathbf{n}\rangle=0, &\qquad\text{on}~\Sigma_R(y),
\end{cases}
\end{equation}
where $y\in\overline{\mathbb R_+^n}$ and $R>0$. Then for any $\lambda\in(0,\Lambda_1), y\in\overline{\mathbb R_+^n}$, it holds
\begin{equation}\label{R-interior-BM}
\int_{B^+_R(y)}e^{\lambda|u-h|}\mathrm dx\leq\frac{|B^+_R(y)|}{1-\lambda\Lambda_1^{-1}}.
\end{equation}
In addition, for any $\mu\in(0,\frac{n-1}{2n}\Lambda_1), y\in\partial\mathbb R_+^n$, it holds
\begin{equation}\label{R-boundary-BM}
\int_{\Sigma_R(y)}e^{\mu|u-h|}\mathrm dx'\leq |\Sigma_R(y)|\big(1+\frac{C\mu}{\frac{n-1}{2n}\Lambda_1-\mu}\big),
\end{equation}
for some constant $C=C(n)>0$, where 
$$\Lambda_1=2^{-\frac{n-2}{n-1}}S_1^{\frac{n}{n-1}}(\int_{\Sigma_R(y)}e^{u}\mathrm dx')^{-\frac{1}{n-1}}$$
and $S_1$ is the best constant of the embedding $\mathcal D^{1,1}(\mathbb R^n_+)\hookrightarrow L^{\frac{n}{n-1}}(\mathbb R^n_+)$, where
$$\mathcal D^{1,1}(\mathbb R^n_+)=\{u\in L^{\frac{n}{n-1}}(\mathbb R^n_+):\nabla u\in L^1(\mathbb R^n_+)\}.$$
\end{lemma}
\begin{proof}
For fixed positive constant $k$, define
\begin{equation*}\label{R-Tk-def}
T_{k}(u):= 
\begin{cases}
u, &\quad\text{if}~|u|\leq k,\\
k\frac{u}{|u|}, &\quad\text{if}~|u|>k.
\end{cases}
\end{equation*}
Since $u\in W^{1,n}_{\operatorname{loc}}(\overline{\mathbb R_+^n})$, we have $T_{k}(u-h)\in W^{1,n}_0(B_R(y)\cap \overline{\mathbb R_+^n})\cap L^{\infty}(B_R^+(y))$, where $T_{k}(u-h)\in W^{1,n}_0(B_R(y)\cap \overline{\mathbb R_+^n})$ means $T_{k}(u-h)|_{\partial B_R^+(y)}=0$ and $T_{k}(u-h)\in W^{1,n}(B_R^+(y))$.
   
Given $l>0$, testing \eqref{I-hL-eq} and \eqref{R-uh-eq} with $T_{k+l}(u-h)-T_{k}(u-h)$, it yields
\begin{eqnarray}\label{R-hu-eq}
&&\int_{B^+_R(y)}\bigl\langle a(\nabla u)-a(\nabla h),\nabla\big[T_{k+l}(u-h)- T_{k}(u-h)\big]\bigl\rangle\mathrm dx \nonumber\\
&=&\int_{\Sigma_R(y)}e^{u}\big[T_{k+l}(u-h)-T_{k}(u-h)\big]\mathrm dx'.
\end{eqnarray}
Notice that
\begin{equation*}
T_{k+l}(u-h)-T_{k}(u-h)=
\begin{cases}
0,&\qquad\text{if}~|u-h|\leq k,\\
u-h-k\frac{u-h}{|u-h|},&\qquad\text{if}~k<|u-h|\leq k+l,\\
l\frac{u-h}{|u-h|},&\qquad\text{if}~|u-h|>k+l.
\end{cases}
\end{equation*}
Hence, $T_{k+l}(u-h)-T_{k}(u-h)\leq l$ and
\begin{equation*}
\nabla\big[T_{k+l}(u-h)- T_{k}(u-h)\big]=
\begin{cases}
\nabla (u-h),&\qquad\text{if }k<|u-h|\leq k+l,\\
0,&\qquad\text{otherwise.}
\end{cases}
\end{equation*}
Let $p,q\in\mathbb R^n$, recall the inequality 
$$\langle |p|^{n-2}p-|q|^{n-2}q,p-q\rangle\geq C_n|p-q|^n,$$ where $C_n$ is a positive constant depending only on $n$ (in fact, we can choose $C_n=2^{2-n}$, see Chapter 12 in \cite{L2017}). For simplicity, write $H_k:=\{x\in\mathbb R^n_+\,:\,k<|u-h|\leq k+l\}$. Combining the above with \eqref{R-hu-eq} yields
\begin{eqnarray}\label{R-grad-esti}
\frac{1}{l}\int_{B^+_R(y)\cap{H_k} }|\nabla (u-h)|^n\mathrm dx\leq2^{n-2}\int_{\Sigma_R(y)}e^{u}\mathrm dx'=:C_0.
\end{eqnarray}
 
Now we prove the property \eqref{R-interior-BM}. Set 
\begin{equation}\label{R-Phi1-def}
\Phi_1(k):=|\{x\in B_R^+(y): |u-h|>k\}|.
\end{equation}
By Sobolev inequality, we obtain
\begin{eqnarray}\label{R-Phi1-Sobolev}
\Phi_1(k+l)^{\frac{n-1}{n}}&\leq& \frac{1}{l}\big(\int_{B^+_R(y)}|T_{k+l}(u-h)-T_{k}(u-h)|^{\frac{n}{n-1}}\mathrm dx\big)^{\frac{n-1}{n}} \nonumber\\
&\leq& \frac{1}{lS_1}\int_{B^+_R(y)}\big|\nabla\big[ T_{k+l}(u-h)-T_{k}(u-h)\big]\big|\mathrm dx\nonumber\\
&=&\frac{1}{lS_1}\int_{B^+_R(y)\cap{H_k}}|\nabla (u-h)|\mathrm dx,
\end{eqnarray} 
where $S_1$ is the Sobolev constant of the embedding $D^{1,1}(\mathbb R^n_+)\hookrightarrow L^{\frac{n}{n-1}}(\mathbb R^n_+)$. Using H\"{o}lder inequality, \eqref{R-grad-esti} and the definition \eqref{R-Phi1-def} of $\Phi_1$, one has
\begin{eqnarray}\label{R-Holder-1}
& &\int_{B^+_R(y)\cap{H_k}}|\nabla (u-h)|\mathrm dx\nonumber\\
&\leq&\big(\int_{B^+_R(y)\cap{H_k}}|\nabla (u-h)|^n\mathrm dx\big)^{\frac 1 n}\big(\int_{B^+_R(y)\cap{H_k}}1\mathrm dx\big)^{\frac{n-1}{n}}\nonumber\\
&=&\big(\int_{B^+_R(y)\cap{H_k}}|\nabla (u-h)|^n\mathrm dx\big)^{\frac 1 n}\big| B^+_R(y)\cap{H_k}\big|^{\frac{n-1}n}\nonumber\\
&\leq&(lC_0)^{\frac 1 n}\big(\Phi_1(k)-\Phi_1(k+l)\big)^{\frac{n-1}{n}}.
\end{eqnarray}
Combining \eqref{R-Phi1-Sobolev} and \eqref{R-Holder-1} yields
\begin{equation}\label{R-Phi1-esti}
\Phi_1(k+l)\leq\frac{\Phi_1(k)-\Phi_1(k+l)}{l\Lambda_1},
\end{equation} 
where $$\Lambda_1=(\frac{S_1^n}{C_0})^{\frac{1}{n-1}}=2^{-\frac{n-2}{n-1}}S_1^{\frac{n}{n-1}}(\int_{\Sigma_R(y)}e^{u}\mathrm dx')^{-\frac{1}{n-1}}.$$ 
Since $\Phi_1$ is non-increasing, by passing to the limit as $l\rightarrow 0+$, we obtain 
\begin{equation*}\label{R-Phi1-derive}
\Phi_1(k)\leq-\frac{1}{\Lambda_1}\Phi_1'(k)\quad \text{a.e.}\, k\in \mathbb{R}_+.
\end{equation*}
Integrating with respect to $k$ from 0 to $k$, it can be deduced that
\begin{equation*}\label{log-Phi-1-inq}
\ln\frac{\Phi_1(k)}{\Phi_1(0)}\leq-\Lambda_1k,
\end{equation*}
that is,
\begin{equation}\label{R-Phi1-e}
\Phi_1(k)\leq |B_R^+(y)|e^{-\Lambda_1k}.
\end{equation}
It follows from \eqref{R-Phi1-e} that
\begin{eqnarray}\label{R-interior-esti}
\int_{B^+_R(y)}e^{\lambda|u-h|}\mathrm dx &=& |B^+_R(y)|+\lambda\int_{B^+_R(y)}\mathrm dx\int^{|u-h|(x)}_0e^{\lambda k}\mathrm dk \nonumber\\
&=&|B^+_R(y)|+\lambda\int^{\infty}_0e^{\lambda k}\Phi_1(k)\mathrm dk\nonumber\\
&\leq&|B^+_R(y)|+\lambda\int^{\infty}_0e^{\lambda k}|B^+_R(y)|e^{-\Lambda_1k}\mathrm dk\nonumber\\
&=&\frac{|B^+_R(y)|}{1-\lambda\Lambda_1^{-1}}.
\end{eqnarray}
This proves \eqref{R-interior-BM}.

Next, we consider the proof of \eqref{R-boundary-BM}. Set
\begin{equation*}\label{R-Phi2-def}
\Phi_2(k):=|\{x\in\Sigma_R(y):|u-h|>k\}|.
\end{equation*}
By Sobolev trace inequality, we get
\begin{eqnarray*}\label{R-Phi2-esti}
\Phi_2(k+l) &\leq& \frac{1}{l}\int_{\Sigma_R(y)}|T_{k+l}(u-h)-T_{k}(u-h)|\mathrm dx' \nonumber\\
&\leq& \frac{1}{lS_2}\int_{B_R^+(y)}|\nabla T_{k+l}(u-h)-\nabla T_{k}(u-h)|\mathrm dx\nonumber\\
&\leq&\frac{1}{lS_2}\int_{B_R^+(y)\cap{H_k}}|\nabla (u-h)|\mathrm dx,
\end{eqnarray*}
where $S_2$ is the Sobolev trace constant of the embedding $D^{1,1}(\mathbb R^n_+)\hookrightarrow L^{1}(\partial\mathbb R^n_+)$. Again using \eqref{R-Holder-1}, we obtain
\begin{equation}\label{R-Phi2-Holder}
\Phi_2(k+l)\leq\big(\frac{\Phi_1(k)-\Phi_1(k+l)}{l}\big)^{\frac{n-1}{n}}\frac{C_0^{\frac 1 n}}{S_2}.
\end{equation}
Taking $l=k$ in \eqref{R-Phi2-Holder}, one has
\begin{equation}\label{R-2k-esti}
\Phi_2(2k)\leq\big(\frac{\Phi_1(k)-\Phi_1(2k)}{k}\big)^{\frac{n-1}{n}}\frac{C_0^{\frac 1 n}}{S_2}\leq\big(\frac{\Phi_1(k)}{k}\big)^{\frac{n-1}{n}}\frac{C_0^{\frac 1 n}}{S_2}.
\end{equation}
It follows from \eqref{R-Phi1-e} and \eqref{R-2k-esti} that
\begin{equation*}\label{R-Phi2-esti1}
\Phi_2(k) \leq \big(\frac{2\Phi_1(\frac{k}{2})}{k}\big)^{\frac{n-1}{n}}\frac{C_0^{\frac 1 n}}{S_2}\leq C_1C_0^{\frac 1 n}|B_R^+(y)|^{\frac{n-1}{n}}\frac{e^{-\frac{n-1}{2n}\Lambda_1k}}{k^{\frac{n-1}{n}}},
\end{equation*}
where $C_1$ is the positive constant depending only on $n$. Recalling the definition of $\Phi_2(k)$, it yields
\begin{eqnarray*}\label{R-Phi2-esti2}
\Phi_2(k) &\leq& \min\big\{|\Sigma_R(y)|, C_1 C_0^{\frac 1 n}|B_R^+(y)|^{\frac{n-1}{n}}\frac{e^{-\frac{n-1}{2n}\Lambda_1k}}{k^{\frac{n-1}{n}}}\big\} \nonumber\\
&\leq&\min\big\{\mathcal |\Sigma_R(y)|, C_1 C_0^{\frac 1 n}|\Sigma_R(y)|\frac{e^{-\frac{n-1}{2n}\Lambda_1k}}{k^{\frac{n-1}{n}}}\big\}\nonumber\\
&\leq& C\mathcal|\Sigma_R(y)|e^{-\frac{n-1}{2n}\Lambda_1k},
\end{eqnarray*}
where $C=C(n)$ and the second inequality holds as it is sufficient to consider $B_R(y)$, $y\in\partial \mathbb R_+^n$ for the inequality \eqref{R-boundary-BM}. The last inequality holds because $\Phi_2(k) \leq|\Sigma_R(y)|\leq e^{-\frac{n-1}{2n}}|\Sigma_R(y)|e^{-\frac{n-1}{2n}\Lambda_1k}$ if $k\leq\Lambda_1^{-1}$,  otherwise, 
\begin{eqnarray*}
\Phi_2(k)&\leq& C_1C_0^{\frac 1 n}|\Sigma_R(y)|\frac{e^{-\frac{n-1}{2n}\Lambda_1k}}{k^{\frac{n-1}{n}}}\nonumber\\
&\leq& C\Lambda_1^{-\frac{n-1}{n}}|\Sigma_R(y)|e^{-\frac{n-1}{2n}\Lambda_1k}\Lambda_1^{\frac{n-1}{n}}\nonumber\\
&\leq& C|\Sigma_R(y)|e^{-\frac{n-1}{2n}\Lambda_1k}. 
\end{eqnarray*}

By the similar argument in \eqref{R-interior-esti}, we obtain
\begin{eqnarray*}
\int_{\Sigma_R(y)}e^{\mu|u-h|}\mathrm dx'
&=& |\Sigma_R(y)|+\mu\int_{\Sigma_R(y)}\mathrm dx'\int^{|u-h|(x)}_0e^{\mu k}\mathrm dk \nonumber\\
&=& |\Sigma_R(y)|+\mu\int^{\infty}_0e^{\mu k}\Phi_2(k)\mathrm dk \nonumber\\
&\leq&|\Sigma_R(y)|+\mu\int^{\infty}_0C|\Sigma_R(y)|e^{(\mu-\frac{n-1}{2n}\Lambda_1)k}\mathrm dk \nonumber\\
&=&|\Sigma_R(y)|\big(1+\frac{C\mu}{\frac{n-1}{2n}\Lambda_1-\mu}\big).
\end{eqnarray*}
This finishes the proof of Lemma \ref{BM-inq-lem}.
\end{proof}

Based on the above lemmas and Serrin's local $L^{\infty}$ estimate, we will show any solution of \eqref{I-hL-eq} is globally bounded from above, of class $C^{1,\alpha}_{\operatorname{loc}}$ and has the following upper bound.
\begin{proposition}\label{regul-prop}
Let $u\in W_{\operatorname{loc}}^{1,n}(\overline{\mathbb R_+^n})$ be a weak solution of \eqref{I-hL-eq} satisfying \eqref{I-fm-assu}. Then $u^{+}\in L^{\infty}(\mathbb R^n_+)$ and $u\in C^{1,\alpha}_{\operatorname{loc}}(\overline{\mathbb R^n_+})$ for some $\alpha\in(0,1)$. Moreover, there exists a constant $C>0$ such that 
\begin{equation}\label{R-log-esti}
u(x)\leq C-(n-1)\log |x|,\qquad x\in\overline{\mathbb R^n_+}\setminus\{0\}.
\end{equation}
\end{proposition}
\begin{proof}
Firstly, for any $\overline{x}\in\partial\mathbb R^n_+$ and some $r>0$ to be chosen later, let $h\in W^{1,n}(B_r^+(\overline{x}))$ be the weak solution of \eqref{R-uh-eq}. Set
\begin{equation*}
\tilde{h}(x',t):=
\begin{cases}
h(x',t), &\qquad t\geq0,\\
h(x',-t), &\qquad t<0,
\end{cases} 
\end{equation*}
where $x'\in\mathbb R^{n-1}$. In view of \eqref{R-uh-eq}, we obtain that $\tilde{h}(x',t)\in W^{1,n}(B_r(\overline{x}))$ is the weak solution of 
\begin{equation*}\label{R-ht-eq}
\begin{cases}
-\Delta_n \tilde{h}=0, &\qquad \text{in}~B_r(\overline{x}),\\
\tilde{h}=\tilde{u}, &\qquad \text{on}~\partial B_r(\overline{x}),
\end{cases}
\end{equation*}
where $$\tilde{u}(x',t)= 
\begin{cases}
u(x',t), &\qquad t\geq0,\\
u(x',-t), &\qquad t<0.
\end{cases}$$
The comparison principle implies $\tilde{h}\leq \tilde{u}$ in $B_r(\overline{x})$, hence $h\leq u$ in $B_r^+(\overline{x})$. Using Serrin's local $L^{\infty}$-estimate (see \cite{S1964}), we have
\begin{eqnarray}\label{R-ht-esti}
&&\|\tilde{h}^{+}\|_{L^{\infty}(B_{\frac{r}{2}}(\overline{x}))}\leq Cr^{-1}\|\tilde{h}^{+}\|_{L^n(B_r(\overline{x}))}\leq Cr^{-1}\|\tilde{u}^{+}\|_{L^n(B_r(\overline{x}))}\nonumber\\
&&\leq Cr^{-1}\|u^+\|_{L^n(B_r^+(\overline{x}))}\leq Cr^{-1}(\int_{\mathbb R_+^n}e^{\frac{n}{n-1}u}\mathrm dx)^{\frac 1 n}=C(r),
\end{eqnarray}
where $C$ is independent of $\bar x$. Since $\int_{\partial\mathbb R_+^n}e^{u}\mathrm dx'< \infty$, there exists $0<r<1$ small enough such that for any $\overline{x}\in\partial\mathbb R_+^n$, $$\int_{\Sigma_R(\overline{x})}\operatorname{e}^u\mathrm dx'\ll1,$$
which implies $\Lambda_1$ in Lemma \ref{BM-inq-lem} is large enough. Applying Lemma \ref{BM-inq-lem} with $\lambda=2$ and $\mu=2$, there exists a constant $C>0$ depending only on $n,~r$ such that
\begin{equation}\label{R-uh-esti}
\int_{B_r^+(\overline{x})}e^{2|u-h|}\mathrm dx\leq C~\text{and}~\int_{\Sigma_r(\overline{x})}e^{2|u-h|}\mathrm dx'\leq C.
\end{equation}
Combining \eqref{R-ht-esti} and \eqref{R-uh-esti} yields
\begin{equation}\label{R-eu-esti}
\int_{\Sigma_{\frac{r}{2}}(\overline{x})}e^{2u}\mathrm dx'\leq\int_{\Sigma_{\frac{r}{2}}(\overline{x})}e^{2|u-h|}\cdot e^{2h}\mathrm dx'\leq C(r).
\end{equation}
Using Lemma \ref{loc-bdd-lem} and \eqref{R-ht-esti}, we deduce that for any $\overline{x}\in\partial{\mathbb R_+^n}$
\begin{equation}\label{R-infty-esti}
\|u^{+}\|_{L^{\infty}(B_{\frac{r}{4}}^+(\overline{x}))}\leq C(r)(\|u^{+}\|_{L^n(B_{\frac{r}{2}}^+(\overline{x}))}+\|e^{u}\|_{L^2(\Sigma_{\frac{r}{2}}(\overline{x}))})\leq C(r).
\end{equation}
As for $\overline{x}\in \mathbb R_+^n\cap\{x:\operatorname{dist}(x,\partial \mathbb R_+^n)>\frac{r}{4}\}$, using Lemma \ref{loc-bdd-lem} or Serrin's local $L^{\infty}$-estimate, we get
\begin{equation}\label{R-infty-esti2}
\|u^{+}\|_{L^{\infty}(B_{\frac{r}{8}}(\overline{x}))}\leq C(r)\|u^{+}\|_{L^n(B_{\frac{r}{4}}(\overline{x}))}\leq C(r).
\end{equation}
It follows from \eqref{R-infty-esti} and \eqref{R-infty-esti2} that
\begin{equation*}
\|u^{+}\|_{L^{\infty}(\mathbb R^n_+)}\leq C. 
\end{equation*}
Since $u\in L^n_{\operatorname{loc}}(\overline{\mathbb R^n_+})$, the standard Moser iteration implies that $u\in L_{\operatorname{loc}}^{\infty}(\overline{\mathbb R^n_+})$. By elliptic regularity theory (see, e.g. \cite{L1988, L1991, T1984}), we obtain that $u\in C^{1,\alpha}_{\operatorname{loc}}(\overline{\mathbb R^n_+})$ for some $\alpha\in(0,1)$.

Now we prove the property \eqref{R-log-esti}. For $R>0$ and $x\in\mathbb R^n_+$, let $u_R(x)=u(Rx)+(n-1)\log R$. The calculation yields $u_R\in W_{\operatorname{loc}}^{1,n}(\overline{\mathbb R_+^n})$ is the solution of the following equation
\begin{equation*}
\begin{cases}
-\Delta_n u_R=0, &\qquad\text{in}~\mathbb R_+^n,\\
|\nabla u_R|^{n-2}\frac{\partial u_R}{\partial t}=-e^{u_R}, &\qquad\text{on}~\partial\mathbb R_+^n,\\
\int_{\mathbb R_+^n}e^{\frac{n}{n-1}u_R}\mathrm dx+\int_{\partial\mathbb R_+^n}e^{u_R}\mathrm dx'<\infty.
\end{cases}
\end{equation*}
For fixed $x_0\in \mathbb{S}^{n-1}\cap\mathbb R^n_+$, let $h_R\in W^{1,n}(B_{r_0}^+(x_0))$ be the weak solution of 
\begin{equation*}\label{R-hR-eq}
\begin{cases}
-\Delta_n h_R=0, &\qquad\text{in}~B_{r_0}^+(x_0),\\
h_R=u_R, &\qquad\text{on}~\partial B_{r_0}(x_0)\cap\mathbb R_+^n,\\
|\nabla h_R|^{n-2}\frac{\partial h_R}{\partial t}=0, &\qquad\text{on}~ \Sigma_{r_0}(x_0),
\end{cases}
\end{equation*}
where $r_0>0$ is small enough. By the similar argument in \eqref{R-ht-esti}, we can deduce 
\begin{eqnarray*}\label{R-hR-inq}
&&\|\tilde{h}_R^{+}\|_{L^{\infty}(B_{\frac{r_0}{2}}(x_0))}\leq Cr_0^{-1}\|\tilde{h}_R^{+}\|_{L^n(B_{r_0}(x_0))}\leq Cr_0^{-1}\|\tilde{u}_R^{+}\|_{L^n(B_{r_0}(x_0))}\nonumber\\
&&\leq Cr_0^{-1}(\int_{\mathbb R_+^n}e^{\frac{n}{n-1}u_R}\mathrm dx)^{\frac 1 n}=Cr_0^{-1}(\int_{\mathbb R_+^n}e^{\frac{n}{n-1}u}\mathrm dx)^{\frac 1 n}=C(r_0),
\end{eqnarray*}
where notation $f^{+}=\max\{f,0\}$  and $C$ is independent of $r_0$, $R$ and $x_0$. Since $r_0$ is small enough, by the argument in \eqref{R-uh-esti}-\eqref{R-infty-esti}, there exists a positive constant $C$ independent of $R$ and $x_0$ such that
\begin{equation}\label{R-uR-esti}
\|u_R^{+}\|_{L^{\infty}(B_{\frac{r_0}{2}}^+(x_0))}\leq C.
\end{equation}
Notice that $\mathbb{S}^{n-1}\cap\mathbb R_+^n$ can be covered by a finite number of sets $B_{\frac{r_0}{2}}^+(x_0)$, $x_0\in \mathbb{S}^{n-1}$, combining with \eqref{R-uR-esti}, we obtain
\begin{equation}\label{R-uR-esti1}
\|u_R^{+}\|_{L^{\infty}(\mathbb{S}^{n-1}\cap\mathbb R_+^n)}\leq C,
\end{equation}
where $C$ is independent of $R$. By \eqref{R-uR-esti1}, we deduce 
$$u(x)+(n-1)\log |x|\leq C.$$ This shows \eqref{R-log-esti} and Proposition \ref{regul-prop}.
\end{proof}

\section{\textbf{Asymptotic estimates}}\label{Section 4}

In this section, we shall prove the sharp asymptotic estimates on both $u$ and $\nabla u$ at infinity, and the asymptotic integral estimate on second-order derivatives of $u$. 

\subsection{Logarithmic asymptotics behavior at infinity}  We first show the following logarithmic asymptotic behavior at infinity. It is similar to the method of \cite{CL2024}.

\begin{proposition}\label{log-esti-lem}
Let $u\in W_{\operatorname{loc}}^{1,n}(\overline{\mathbb R_+^n}\setminus\overline{B_{1}})$ be a nonnegative solution of the following equation
\begin{equation}\label{A-log-eq}
\begin{cases}
-\Delta_n u=0, &\qquad\text{in}~\mathbb R_+^n\setminus \overline{B_{1}},\\
a(\nabla u)\cdot\mathbf{n}=f(x), &\qquad\text{on}~\partial\mathbb R_+^n\setminus \overline{B_{1}},
\end{cases} 
\end{equation}
where $f\in L^{1}(\partial\mathbb R_+^n\setminus \overline{B_{1}})$ is a scalar function satisfying
\begin{equation}\label{A-f-condition}
-\frac{C}{|x|^{n-1}}\leq f(x)\leq0\qquad a.e.x\in\partial\mathbb R_+^n\setminus \overline{B_{1}},
\end{equation}
for some constant $C>0$. If $u(x)\rightarrow+\infty$ as $|x|\rightarrow+\infty$ with $x\in\mathbb R_+^n$, then there exists a constant $L>1$ such that 
\begin{equation}\label{A-log-esti}
\frac{1}{L}\log |x|\leq u(x)\leq L\log |x|
\end{equation}
for any $x\in\mathbb R^n_+$ with $|x|$ large enough.   
\end{proposition}

We prove this through the following lemmas. Let $J$ denote the class of Lipschitz continuous functions defined on $\overline{\mathbb R_+^n}\setminus B_{1}$ which are identically one in some neighborhood of infinity but vanish in some neighborhood of $\partial B_{1}\cap\overline{\mathbb R_+^n}$.
 
\begin{lemma}\label{const-int-lem}
Let $u\in W_{\operatorname{loc}}^{1,n}(\overline{\mathbb R_+^n}\setminus\overline{B_{1}})$ be a nonnegative solution of the following equation
\begin{equation}\label{A-hat-eq}
\begin{cases}
\operatorname{div}(\widehat{a}(\nabla u))=0, &\qquad\text{in}~\mathbb R_+^n\setminus \overline{B_{1}},\\
\widehat{a}(\nabla u)\cdot\mathbf{n}=f(x), &\qquad\text{on}~\partial\mathbb R_+^n\setminus \overline{B_{1}},
\end{cases} 
\end{equation}
where $\widehat{a}:\mathbb R^n\rightarrow\mathbb R^n$ satisfies
\begin{equation}\label{A-a-condition}
|\widehat{a}(\xi)|\leq c_0|\xi|^{n-1}~\text{and}~\langle \widehat{a}(\xi),\xi\rangle\geq c_{1}|\xi|^n,\qquad\forall \xi\in\mathbb R^n
\end{equation}
for some constants $c_0,~c_{1}>0$ and $f\in L^{1}(\partial\mathbb R_+^n\setminus \overline{B_{1}})$ is a scalar function satisfying \eqref{A-f-condition}. Then there exists a constant $K>0$ such that
\begin{equation*}
\int_{\mathbb R_+^n\setminus B_{1}}\langle \widehat{a}(\nabla u),\nabla \vartheta\rangle \mathrm dx-\int_{\partial\mathbb R_+^n\setminus B_{1}}f\vartheta \mathrm dx'=K
\end{equation*}
for any $\vartheta\in J$.
\end{lemma}
\begin{proof}
Given any $\vartheta_{1}, \vartheta_{2}\in J$, substituting their difference $\vartheta_{1}-\vartheta_{2}$ into the weak formulation of equation \eqref{A-hat-eq}, we have
\begin{eqnarray*}
&&\int_{\mathbb R_+^n\setminus B_{1}}\langle \widehat{a}(\nabla u),\nabla \vartheta_{1}\rangle \mathrm dx-\int_{\partial\mathbb R_+^n\setminus B_{1}}f\vartheta_{1} \mathrm dx' \\
&=&\int_{\mathbb R_+^n\setminus B_{1}}\langle \widehat{a}(\nabla u),\nabla \vartheta_{2}\rangle \mathrm dx-\int_{\partial\mathbb R_+^n\setminus B_{1}}f\vartheta_{2} \mathrm dx',
\end{eqnarray*}
which implies the existence of $K$. Next, we demonstrate that $K>0$. Without loss of generality, we may assume $u<0$ on $\{x\in\overline{\mathbb R_+^n}: |x|=2\}$, since otherwise we can consider the shifted function $u-k$ instead of $u$ for some appropriate constant $k$. Define
\begin{equation*}
\hat{u}(x)=
\begin{cases}
0, &\quad\text{if}~u(x)\leq0,\\
u, &\quad\text{if}~0<u(x)<1,\\
1, &\quad\text{if}~u(x)\geq1,
\end{cases}
\end{equation*}
for $|x|>2$ and $\hat{u}\equiv0$ for $|x|\leq2$. The fact $u(x)\rightarrow+\infty$ as $|x|\rightarrow+\infty$ implies that $\hat{u}\in J$. Using the conditions \eqref{A-a-condition} and \eqref{A-f-condition}, we deduce
\begin{eqnarray*}
K&=&\int_{\mathbb R_+^n\setminus B_{1}}\langle \widehat{a}(\nabla u),\nabla \hat{u}\rangle \mathrm dx-\int_{\partial\mathbb R_+^n\setminus B_{1}}f\hat{u} \mathrm dx'\\
&\geq&\int_{\mathbb R_+^n\setminus B_{1}}\langle \widehat{a}(\nabla \hat{u}),\nabla \hat{u}\rangle \mathrm dx\\
&\geq&c_{1}\int_{\mathbb R_+^n\setminus B_{1}}|\nabla \hat{u}|^n\mathrm dx>0.
\end{eqnarray*} 
The proof is complete.
\end{proof}

Next, we establish the following Harnack inequality.
\begin{lemma}\label{Harnack-lem}
Let $u$ be as in Proposition \ref{log-esti-lem} replacing $B_{1}$ with $B_{r_0}$, $r_0>0$. Then for any $R>r_0$ it holds
\begin{equation}\label{A-Harnack-inq}
\sup\limits_{B_{2R}^{+}\setminus B_{R}^{+}} u\leq C\inf\limits_{B_{2R}^{+}\setminus B_{R}^{+}} u+C,
\end{equation}
where $C$ is a positive constant independent of $R$.
\end{lemma}
\begin{proof}
For $k>0$ and $R>\theta r_0$ where $\theta>1$, define $v:=u+k$ and $v_{R}(x):=v(Rx)=u(Rx)+k$. We have the following two conclusions hold.

1. \textbf{Local boundedness estimate.}
There exists a constant $\mu_0>0$ independent of $R$ such that for any $k\geq\mu_0$ and $t>0$, the following inequality holds:
\begin{equation}\label{A-lb-inq}
\|v_{R}\|_{L^{\infty}(B_{2}^{+}\setminus B_{1}^{+})}\leq C\|v_{R}\|_{L^{t}(B_{3}^{+}\setminus B_{\frac{1}{\theta_0}}^{+})},
\end{equation}
for $1<\theta_0<\theta$ and some constant $C$ independent of $R$.

2. \textbf{Weak Harnack inequality.}
There exist $\mu_{1}>0$ and $p_{1}>0$ independent of $R$ such that for any $k\geq\mu_{1}$,
\begin{equation}\label{A-wH-inq}
\inf_{B_{2}^{+}\setminus B_{1}^{+}}v_{R}\geq C\|v_{R}\|_{L^{p_{1}}(B_{3}^{+}\setminus B_{\frac{1}{\theta_0}}^{+})},
\end{equation} 
for $1<\theta_0<\theta$ and some constant $C$ independent of $R$.

Now, we are ready to prove it.
 
\noindent\textbf{Step 1. Proof of \eqref{A-lb-inq}.} We  show \eqref{A-lb-inq} by the classical Moser iteration technique. Combining  equation \eqref{A-log-eq} with the definition of $v_{R}$ yields $v_{R}\in W_{\operatorname{loc}}^{1,n}(\overline{\mathbb R_+^n}\setminus \overline{B_{\frac{1}{\theta}}})$ satisfies the following equation
\begin{equation}\label{A-vR-eq}
\begin{cases}
-\Delta_n v_{R}=0, &\qquad\text{in}~\mathbb R_+^n\setminus \overline{B_{\frac{1}{\theta}}},\\
|\nabla v_{R}|^{n-2}\frac{\partial v_{R}}{\partial t}=-f_{R}, &\qquad\text{on}~\partial\mathbb R_+^n\setminus \overline{B_{\frac{1}{\theta}}},
\end{cases}
\end{equation}
where $f_{R}(x) := R^{n-1}f(Rx)=:\widetilde{f}(x)v_{R}^{n-1}$. This implies that
\begin{equation*}
\int_{\mathbb R_+^n\setminus \overline{B_{\frac{1}{\theta}}}}a(\nabla v_{R})\cdot \nabla \psi\mathrm dx=\int_{\partial\mathbb R_+^n\setminus \overline{B_{\frac{1}{\theta}}}}\widetilde{f}(x)v_{R}^{n-1}\psi\mathrm dx',\qquad\forall\psi\in C_{c}^{\infty}(\overline{\mathbb R_+^n}\setminus \overline{B_{\frac{1}{\theta}}}).
\end{equation*}
The condition \eqref{A-f-condition} and the definition of $v_{R}$ imply that 
\begin{equation*}
|\widetilde{f}(x)|=\big|\frac{R^{n-1}f(Rx)}{v_{R}^{n-1}}\big|\leq\frac{C}{k^{n-1}},\quad|x|>\frac{1}{\theta},
\end{equation*}
where $C$ is a positive constant which depends only on $n$ and $\theta$. Since $v_R\ge k$, by the standard Moser iteration, we deduce that there exists a constant $\mu_0>0$ independent of $R$ such that for any $k\geq\mu_0$,
\begin{eqnarray*}
\|v_{R}\|_{L^{\infty}(B_{2}^{+}\setminus B_{1}^{+})}&\leq &C(\|v_{R}\|_{L^{t}(B_{3}^{+}\setminus B_{\frac{1}{\theta_0}}^{+})}+\|\widetilde{f}\|_{L^{\infty}(B_{3}^{+}\setminus B_{\frac{1}{\theta_0}}^{+})})\nonumber\\
&\leq& C\|v_{R}\|_{L^{t}(B_{3}^{+}\setminus B_{\frac{1}{\theta_0}}^{+})}+\frac{1}{2}\|v_{R}\|_{L^{\infty}(B_{2}^{+}\setminus B_{1}^{+})},
\end{eqnarray*}
which implies \eqref{A-lb-inq}.

\noindent\textbf{Step 2. Proof of \eqref{A-wH-inq}.} Denote $w:=v_{R}^{-1}$. A direct computation yields that $w$ satisfies the following equation
\begin{equation}\label{A-w-eq}
\begin{cases}
-\Delta_n w\leq 0,&\qquad\text{in}~\mathbb R_+^n\setminus \overline{B_{\frac{1}{\theta}}},\\
|\nabla w|^{n-2}\frac{\partial w}{\partial t}=\widetilde{f}w^{n-1}, &\qquad\text{on}~\partial\mathbb R_+^n\setminus \overline{B_{\frac{1}{\theta}}}.\\
\end{cases}
\end{equation}
It follows from equation \eqref{A-w-eq} that
\begin{equation*}
\int_{\mathbb R_+^n\setminus \overline{B_{\frac{1}{\theta}}}}a(\nabla w)\cdot\nabla \psi\mathrm dx\leq-\int_{\partial\mathbb R_+^n\setminus \overline{B_{\frac{1}{\theta}}}}\widetilde{f}(x)w^{n-1}\psi\mathrm dx',\quad\forall\psi\in C_{c}^{\infty}(\overline{\mathbb R_+^n}\setminus \overline{B_{\frac{1}{\theta}}}).   
\end{equation*}
Using the standard Moser iteration yields that there exists $\mu_{1}>0$ independent of $R$ such that for any $k\geq\mu_{1}$, 
\begin{equation*}
\|w\|_{L^{\infty}(B_{2}^{+}\setminus B_{1}^{+})}\leq C\|w\|_{L^{t}(B_{3}^{+}\setminus B_{\frac{1}{\theta_0}}^{+})},
\end{equation*}
for some constant $C>0$ independent of $R$ and any $t>0$. Therefore, recalling the definition of $w$, we obtain that
\begin{equation}\label{A-vR-inq}
\inf\limits_{B_{2}^{+}\setminus B_{1}^{+}}v_{R}\geq C\|v_{R}^{-1}\|_{L^{t}(B_{3}^{+}\setminus B_{\frac{1}{\theta_0}}^{+})}^{-1}.
\end{equation}
Define 
$$\varphi:=\log v_{R}-\frac{1}{|B_{3}^{+}\setminus B_{\frac{1}{\theta_0}}^{+}|}\int_{B_{3}^{+}\setminus B_{\frac{1}{\theta_0}}^{+}}\log v_{R}\mathrm dx$$ 
and 
$$m:=\frac{1}{|B_{3}^{+}\setminus B_{\frac{1}{\theta_0}}^{+}|}\int_{B_{3}^{+}\setminus B_{\frac{1}{\theta_0}}^{+}}\log v_{R}\mathrm dx.$$
We claim that $\varphi\in \operatorname{BMO}$. Indeed, from the definition of $\varphi$ and equation \eqref{A-vR-eq}, it holds
\begin{eqnarray}\label{A-varphi-int}
&&(n-1)\int_{\mathbb R_+^n\setminus \overline{B_{\frac{1}{\theta}}}}|\nabla \varphi|^n\phi^n\mathrm dx=(n-1)\int_{\mathbb R_+^n\setminus \overline{B_{\frac{1}{\theta}}}}|\nabla v_{R}|^nv_{R}^{-n}\phi^n\mathrm dx\nonumber\\
&&=-\int_{\mathbb R_+^n\setminus \overline{B_{\frac{1}{\theta}}}}a(\nabla v_{R})\cdot\nabla(v_{R}^{1-n}\phi^n)\mathrm dx+\int_{\mathbb R_+^n\setminus \overline{B_{\frac{1}{\theta}}}}a(\nabla v_{R})\cdot\nabla\phi^nv_{R}^{1-n}\mathrm dx\nonumber\\
&&=-\int_{\partial\mathbb R_+^n\setminus \overline{B_{\frac{1}{\theta}}}}\widetilde{f}\phi^n\mathrm dx'+n\int_{\mathbb R_+^n\setminus \overline{B_{\frac{1}{\theta}}}}\phi^{n-1}|\nabla \varphi|^{n-2}\nabla \varphi\cdot\nabla\phi\mathrm dx 
\end{eqnarray}
for $\phi\in C_{c}^{\infty}(\mathbb R^n)$. By \eqref{A-varphi-int} and Young's inequality, one has 
\begin{equation}\label{A-varphi-inq}
\int_{\mathbb R_+^n\setminus \overline{B_{\frac{1}{\theta}}}}|\nabla \varphi|^n\phi^n\mathrm dx\leq C(\int_{\partial\mathbb R_+^n\setminus \overline{B_{\frac{1}{\theta}}}}|\widetilde{f}|\phi^n\mathrm dx'+\int_{\mathbb R_+^n\setminus \overline{B_{\frac{1}{\theta}}}}|\nabla \phi|^n\mathrm dx),
\end{equation}
where $C$ is a positive constant depending only on $n$. For any $B_{r}(y)\subset B_{3}\setminus B_{\frac{1}{\theta_0}}$, choose $\phi$ as
\begin{equation*}
\phi=
\begin{cases}
1,\quad&   \text{ in }B_{r}(y),\\
0,\quad&\text{ in }B_{2r}^{c}(y),\\
\end{cases}
\end{equation*}
with  $0\leq\phi\leq 1,~|\nabla\phi|\leq\frac{2}{r}$. Using \eqref{A-varphi-inq} and Poincar\'{e} inequality, we have
\begin{equation*}
\int_{B_{r}^{+}(y)}|\varphi-(\varphi)_{r,y}|^n\mathrm dx\leq Cr^n\int_{B_{r}^{+}(y)}|\nabla \varphi|^n\mathrm dx\leq Cr^n(r^{n-1}+1)\leq Cr^n,
\end{equation*}
for some positive constant $C$ depending only on $n$ and $\theta$,  which implies $\varphi\in\operatorname{BMO}(B_{3}^{+}\setminus B_{\frac{1}{\theta_0}}^{+})$. 

Applying the John-Nirenberg lemma yields that there exist positive constants $p_{1}$ and $C$ depending only on $n$ and $\theta$ such that 
\begin{equation*}
\int_{B_{3}^{+}\setminus B_{\frac{1}{\theta_0}}^{+}}e^{p_{1}|\varphi|}\mathrm dx\leq C.
\end{equation*}
Noticing that 
$$e^{p_{1}\varphi}=e^{-p_{1}m}v_{R}^{p_{1}}\text{ and }e^{-p_{1}\varphi}=e^{p_{1}m}v_{R}^{-p_{1}},$$
 it yields
\begin{eqnarray}\label{A-p1-inq}
\|v_{R}^{-1}\|_{L^{p_{1}}(B_{3}^{+}\setminus B_{\frac{1}{\theta_0}}^{+})}\cdot\|v_{R}\|_{L^{p_{1}}(B_{3}^{+}\setminus B_{\frac{1}{\theta}}^{+})}
&=& (\int_{B_{3}^{+}\setminus B_{\frac{1}{\theta_0}}^{+}}e^{-p_{1}\varphi}\mathrm dx)^{\frac{1}{p_{1}}}(\int_{B_{3}^{+}\setminus B_{\frac{1}{\theta_0}}^{+}}e^{p_{1}\varphi}\mathrm dx)^{\frac{1}{p_{1}}}\nonumber\\
&\leq&(\int_{B_{3}^{+}\setminus B_{\frac{1}{\theta_0}}^{+}}e^{p_{1}|\varphi|}\mathrm dx)^{\frac{2}{p_{1}}}\leq C.
\end{eqnarray}
Combining \eqref{A-vR-inq} with \eqref{A-p1-inq}, we obtain \eqref{A-wH-inq}. 

As a direct consequence of \eqref{A-lb-inq} and \eqref{A-wH-inq}, we deduce that inequality \eqref{A-Harnack-inq} holds.
\end{proof}

From the Harnack inequality in Lemma \ref{Harnack-lem}, we derive the following Harnack inequality on the sphere, which will be applied directly in the proof of the logarithmic asymptotic estimate \eqref{A-log-esti}.
\begin{lemma}\label{sphere-Harnack-lem}
Let $u$ satisfy the assumptions of Proposition \ref{log-esti-lem}. Then there exist constants $C>0$ and $\kappa_{1}>0$ such that
\begin{equation}\label{A-H-inq2}
\max\limits_{|x|=\kappa,~x\in\overline{\mathbb R_+^n}}u\leq C \min\limits_{|x|=\kappa,~x\in\overline{\mathbb R_+^n}}u
\end{equation}
for any $\kappa\geq\kappa_{1}$.
\end{lemma}
\begin{proof}
Define the rescaled functions $u_{\kappa}$ and $f_{\kappa}$ as 
$$u_{\kappa}(x) := u(\kappa x), \quad f_{\kappa}(x) := \kappa^{n-1}f(\kappa x) \quad \text{for } x \in \mathbb R^n_+ \setminus \overline{B_{\frac{1}{2}}},$$
where $\kappa>2$. The equation \eqref{A-log-eq} can become
\begin{equation}\label{A-kappa-eq}
\begin{cases}
-\Delta_n u_{\kappa} = 0, \quad & \text{in } \mathbb R_+^n \setminus \overline{B_{\frac{1}{2}}}, \\
a(\nabla u_{\kappa}) \cdot \mathbf{n} = f_{\kappa}(x), \quad &  \text{on } \partial\mathbb R_+^n \setminus \overline{B_{\frac{1}{2}}},
\end{cases}
\end{equation}
 where $f_{\kappa}$ still satisfies the condition \eqref{A-f-condition}.

Applying the Harnack inequality in Lemma \ref{Harnack-lem} to \eqref{A-kappa-eq} with $r_0=\frac{1}{2}$ and $\frac{1}{2}<R<1$ which implies $\mathbb S^{n-1}\subset B_{2R}\setminus B_{R}$, we obtain the following estimate  
\begin{equation}\label{A-kappa-inq}
\max_{\mathbb S^{n-1} \cap \overline{\mathbb R_+^n}} u_{\kappa} \leq \hat C \min_{\mathbb S^{n-1} \cap \overline{\mathbb R_+^n}} u_{\kappa} + \hat C,
\end{equation}
where $\hat C$ is a positive constants independent of $\kappa$.

Furthermore, since $u(x) \to +\infty$ as $ |x| \to +\infty $, there exists $\kappa_{1} > 0$ such that 
\begin{equation*}
\min_{|x| = \kappa,~ x\in\overline{\mathbb R_+^n}} u \geq \hat C \quad \text{for all } \kappa \geq \kappa_{1}.
\end{equation*}
Combining this with \eqref{A-kappa-inq}, we obtain the desired inequality \eqref{A-H-inq2}.
\end{proof}

\noindent\textbf{Proof of Proposition \ref{log-esti-lem}.}
Without loss of generality, we may assume that $u(x)<0$ on $|x|=2$ by subtracting a suitable constant. For $\kappa>2$, define
$$m(\kappa):=\min\limits_{|x|=\kappa, x\in\overline{\mathbb R_+^n}}u(x)~\text{and}~M(\kappa):=\max\limits_{|x|=\kappa, x\in\overline{\mathbb R_+^n}}u(x).$$ 
The fact $u(x)\rightarrow +\infty$ as $|x|\rightarrow +\infty$ implies that there exists $\kappa_0>2$ such that $m(\kappa)>0$ for any $\kappa\geq\kappa_0$. By Lemma \ref{const-int-lem} and the same approach as the proof of (6.8) in \cite{CL2024}, we have 
\begin{equation}\label{A-m-esti}
m(\kappa)\leq C\log \frac{\kappa}{2},
\end{equation}
where $C>0$ independent of $\kappa$. The details are as follows.

Define
\begin{equation*}
u_{1}(x)=   
\begin{cases}
0,&\quad\text{in }\overline{B_2^+},\\
\min\big\{\max\{0, u(x)\},m(\kappa)\big\},&\quad\text{in }\overline{B_\kappa^+\backslash B_2^+},\\
m(\kappa),&\quad\text{in }\overline{\mathbb R_+^n}\backslash B_\kappa,\\
\end{cases}
\end{equation*}
where $\kappa\geq\kappa_0$. One can readily see that $\frac{u_{1}}{m(\kappa)}\in J$. Using Lemma \ref{const-int-lem}, we obtain
\begin{eqnarray}\label{A-mK-inq}
m(\kappa)K=\int_{\mathbb R_+^n\setminus B_{1}}\langle a(\nabla u),\nabla u_{1}\rangle \mathrm dx-\int_{\partial\mathbb R_+^n\setminus B_{1}}fu_{1}\mathrm dx'\geq\int_{\mathbb R_+^n}|\nabla u_{1}|^n\mathrm dx.
\end{eqnarray}

We claim 
\begin{equation}\label{A-u1-inq}
\int_{\mathbb R_+^n}|\nabla u_{1}|^n\mathrm dx\geq \frac{1}{2}|\mathbb S^{n-1}|(m(\kappa))^n(\log\frac{\kappa}{2})^{1-n}.   
\end{equation}
Indeed, consider minimizing the functional $\int_{\mathbb R_+^n}|\nabla \phi|^n\mathrm dx$ over all Lipschitz continuous functions $\phi(x)$ satisfying
\begin{equation*}
\phi(x)=
\begin{cases}
0,&\quad\text{in }\overline{B_2^+},\\
m(\kappa),&\quad\text{in }\overline{\mathbb R_+^n}\backslash B_\kappa.\\
\end{cases}
\end{equation*}
This variational problem admits a unique minimizer $\varphi_{1}$ which satisfies the following equation
\begin{equation*}
\begin{cases}
\Delta_n\varphi_{1}=0,&\quad\text{in }B_{\kappa}^+\setminus \overline{B_{2}^+},\\
\varphi_{1}=0,&\quad\text{on }\partial B_2^+,\\
\varphi_{1}=m(\kappa),&\quad\text{on }\partial B_{\kappa}^+,\\
\end{cases}
\end{equation*}
which implies
\begin{equation*}
\varphi_{1}(x)=
\begin{cases}
0,&\quad\text{in }\overline{B_2^+},\\
m(\kappa)\frac{\log|x|-\log2}{\log\kappa-\log2},&\quad\text{in }\overline{B_\kappa^+\backslash B_2^+},\\
m(\kappa),&\quad\text{in }\overline{\mathbb R_+^n}\backslash B_\kappa.\\
\end{cases}   
\end{equation*}
Therefore, we conclude
\begin{equation*}
\int_{\mathbb R_+^n}|\nabla u_{1}|^n\mathrm dx
\geq \int_{\mathbb R_+^n}|\nabla \varphi_{1}|^n\mathrm dx 
=\frac{1}{2}|\mathbb S^{n-1}|(m(\kappa))^n(\log\frac{\kappa}{2})^{1-n},
\end{equation*}
which finishes the proof of the claim. It follows from \eqref{A-mK-inq} and \eqref{A-u1-inq} that \eqref{A-m-esti} holds.

Let $\kappa_0$ be sufficiently large such that $-\int_{\partial\mathbb R_+^n\setminus B_{\kappa_0}}f\mathrm dx'\leq\frac{K}{2}$, where $K>0$ is given by Lemma \ref{const-int-lem}. Also, from Lemma \ref{const-int-lem} and the approach of the proof of (6.11) in \cite{CL2024}, there exists $C>0$ independent of $\kappa$ such that
\begin{equation}\label{A-max-esti}
M(\kappa)\geq C\log \frac{\kappa}{\kappa_0},\qquad\forall~\kappa\geq\kappa_0.
\end{equation}
The details are as follows. 

Let
\begin{equation*}
\varphi_{2}(x)=
\begin{cases}
0,&\quad\text{in }\overline{B_{\kappa_0}^+},\\
\frac{\log|x|-\log\kappa_0}{\log\kappa-\log\kappa_0},&\quad\text{in }\overline{B_\kappa^+\backslash B_{\kappa_0}^+},\\
1,&\quad\text{in }\overline{\mathbb R_+^n}\backslash B_\kappa,\\
\end{cases}   
\end{equation*}
where $\kappa>\kappa_0$. Using Lemma \ref{const-int-lem} and H\"{o}lder inequality yields  
\begin{eqnarray*}
K&=&\int_{\mathbb R_+^n\setminus B_{1}}\langle a(\nabla u),\nabla \varphi_{2}\rangle \mathrm dx-\int_{\partial\mathbb R_+^n\setminus B_{1}}f\varphi_{2}\mathrm dx'\nonumber\\
&\leq&\big(\int_{B_{\kappa}^+\setminus B_{\kappa_0}^+}|\nabla \varphi_{2}|^n\mathrm dx\big)^{\frac{1}{n}}\big(\int_{B_{\kappa}^+\setminus B_{\kappa_0}^+}|\nabla u|^n\mathrm dx\big)^{\frac{n-1}{n}}+\frac{K}{2},
\end{eqnarray*}
which implies 
\begin{equation}\label{A-Ku-inq}
C(n)K^{\frac{n}{n-1}}\log\frac{\kappa}{\kappa_0}\leq\int_{B_{\kappa}^+\setminus B_{\kappa_0}^+}|\nabla u|^n\mathrm dx.
\end{equation}
Define 
\begin{equation*}
u_{2}(x)=
\begin{cases}
0,&\quad\text{in }\overline{B_2^+},\\
\max\{0,u(x)\},&\quad\text{in }\overline{B_\kappa^+\backslash B_2^+},\\
\min\{M(\kappa),u(x)\},&\quad\text{in }\overline{\mathbb R_+^n}\backslash B_\kappa.\\
\end{cases}   
\end{equation*}
It is clear that $\frac{u_{2}}{M(\kappa)}\in J$. Again, by Lemma \ref{const-int-lem} and $u=u_2$ in $B_{\kappa}^+\backslash B_{\kappa_0}^+$ for $\kappa_0$ large enough, we deduce
\begin{eqnarray}\label{A-KM-inq}
KM(\kappa)&=&\int_{\mathbb R_+^n\setminus B_{1}}\langle a(\nabla u),\nabla u_{2}\rangle \mathrm dx-\int_{\partial\mathbb R_+^n\setminus B_{1}}fu_{2}\mathrm dx'\nonumber\\
&\geq&\int_{\{u=u_{2}\}}|\nabla u|^n\mathrm dx\geq\int_{B_{\kappa}^+\setminus B_{\kappa_0}^+}|\nabla u|^n\mathrm dx.
\end{eqnarray}
By \eqref{A-Ku-inq} and \eqref{A-KM-inq}, we obtain the validity of \eqref{A-max-esti}.

Combining \eqref{A-m-esti}, \eqref{A-max-esti} and Lemma \ref{sphere-Harnack-lem}, it holds
\begin{equation*}
C\log\frac{\kappa}{\kappa_0}\leq m(\kappa)\leq M(\kappa)\leq \bar{C}\log\frac{\kappa}{2},
\end{equation*}
for any $\kappa\geq\kappa_0$, where $C$ and $\bar{C}$ are some positive constants independent of $\kappa$. Hence, this proves \eqref{A-log-esti}, which completes the proof.
\qed

\subsection{Asymptotic estimates on $u$ and $\nabla u$ at infinity}
Based on Proposition \ref{regul-prop} and Proposition \ref{log-esti-lem}, arguing as \cite{CL2024, DGL2024}, we can show the following sharp asymptotic estimates for $u$ and $\nabla u$ at infinity.
\begin{proposition}\label{1order-esti-prop}
Let $u\in W_{\operatorname{loc}}^{1,n}(\overline{\mathbb R_+^n})$ be a weak solution of \eqref{I-hL-eq} satisfying the assumption \eqref{I-fm-assu}. Then there exist $R_0>0$ and $\beta>n-1$ such that
\begin{equation}\label{A-u-esti}
u(x)+\beta \log |x|\in L^{\infty}(\overline{\mathbb R^n_+}\setminus B_{R_0}),
\end{equation}
and
\begin{equation}\label{A-gradu-esti}
\lim_{x\in\overline{\mathbb R_+^n},|x|\rightarrow \infty}|x||\nabla (u(x)+\beta \log |x|)|=0.
\end{equation}
\end{proposition}

To establish the above asymptotic estimates, we also need the following comparison principle.
\begin{lemma}\label{compa-prin-lem}
Let $\Omega\subset\mathbb R^n$ be a bounded domain such that $|\Gamma_0|>0$ where $\Gamma_0:=\partial\Omega\cap\mathbb R^n_+$, and assume that $\mathbb R^n_+\cap\Omega$ is connected. Assume that $u, v\in W^{1,n}(\mathbb R^n_+\cap\Omega)\cap C^{0}((\mathbb R^n_+\cap\Omega)\cup\Gamma_0)$ satisfy
\begin{equation}\label{A-com-eq}
\begin{cases}
-\Delta_n u\leq -\Delta_n v, &\qquad\text{in}~\mathbb R_+^n\cap\Omega,\\
u\leq v, &\qquad\text{on}~\Gamma_0,\\
\langle a(\nabla u), \mathbf{n}\rangle\leq \langle a(\nabla v), \mathbf{n}\rangle, &\qquad \text{on}~\partial(\mathbb R_+^n\cap\Omega)\setminus\Gamma_0.
\end{cases}
\end{equation}
Then $u\leq v$ in $\mathbb R_+^n\cap\Omega$.
\end{lemma}
\begin{proof}
Following the idea of Lemma 2.4 in \cite{CL2022}, we test \eqref{A-com-eq} with $(u-v)^{+}$, and then by integration by parts, it holds
\begin{eqnarray}\label{A-uv-inq}
&&\int_{\mathbb R_+^n\cap\Omega}\langle a(\nabla u)-a(\nabla v), \nabla(u-v)^{+}\rangle \mathrm dx \nonumber\\
&\leq& \int_{\partial(\mathbb R_+^n\cap\Omega)\setminus\Gamma_0}\langle a(\nabla u)-a(\nabla v),\mathbf{n}\rangle(u-v)^{+} \mathrm dx'\leq0.
\end{eqnarray} 
Noting that 
\begin{equation*}
\langle a(p)-a(q),p-q\rangle\geq C|p-q|^n,\qquad \forall p,q\in\mathbb R^n,
\end{equation*}
for some $C>0$ depending only on $n$, and
by \eqref{A-uv-inq} we arrive at
\begin{equation*}
\int_{\mathbb R_+^n\cap\Omega}|\nabla(u-v)^{+}|^n\mathrm dx\leq0,
\end{equation*}
which implies that $(u-v)^{+}$ is constant in $\mathbb R_+^n\cap\Omega$.
Since $(u-v)^{+}=0$ on $\Gamma_0$, we get $(u-v)^{+}=0$ in $\mathbb R_+^n\cap\Omega$, that is, $u\leq v$ in $\mathbb R_+^n\cap\Omega$.
\end{proof}

\begin{lemma}[\cite{CL2022} Lemma 3.2]\label{const-lem}
Let $\sigma\in\mathbb R$ be a constant. Assume that $G(x)\in W^{1,n}_{\operatorname{loc}}(\overline{\mathbb R^n_+}\setminus\{0\})\cap L^{\infty}(\mathbb R^n_+)$ and the function $\sigma \log |x|+G(x)$ satisfies 
\begin{equation*}
\begin{cases}
\Delta_n (\sigma\log |x|+G(x))=0,&\qquad\text{in}~\mathbb R_+^n,\\
\langle a(\nabla(\sigma \log |x|+G(x))),\mathbf{n}\rangle=0, &\qquad \text{on}~\partial\mathbb R_+^n\setminus\{0\}.\\
\end{cases}
\end{equation*}
Then $G(x)$ is a constant function.
\end{lemma}

\noindent\textbf{Proof of Proposition \ref{1order-esti-prop}.}
Let $u$ be a weak solution of \eqref{I-hL-eq}. 
For fixed $R_0>1$ and any $x\in\overline{\mathbb R^n_+}\setminus B_{\frac{R_0}{R}}, \forall R>R_0$, we define 
$$u_{R}(x):=\frac{u(Rx)-\mu_0}{\log R}~\text{and}~U_0:=\sup\limits_{x\in\mathbb R^n_+} u,$$ 
where $\mu_0:=\inf\limits_{x\in\partial B_{R_0}^+} u(x)$.
By \eqref{R-log-esti} in Proposition \ref{regul-prop}, i.e., $u(x)\leq C-(n-1)\log |x|, ~\forall x\in \overline{\mathbb R^n_+}\setminus\{0\}$, we have
\begin{equation}\label{A-eu-esti}
e^{u(x)}\leq e^{C-(n-1)\log |x|}=\frac{C}{|x|^{n-1}},\qquad \text{in}~\overline{\mathbb R^n_+}\setminus\{0\}.
\end{equation} 
Letting $\hat{u}=U_0-u$, we obtain $\hat{u}(x)\geq(n-1)\log |x|+U_0-C\rightarrow+\infty$ as $|x|\rightarrow+\infty$.
Equation \eqref{I-hL-eq} implies that $\hat{u}$ is a solution of the following equation 
\begin{equation*}
\begin{cases}
-\Delta_n\hat{u}=0, &\qquad\text{in} ~\mathbb R_+^n,\\
|\nabla \hat{u}|^{n-2}\frac{\partial\hat{u}}{\partial t}=e^u, &\qquad \text{on}~\partial\mathbb R_+^n.
\end{cases}
\end{equation*}
For $\hat{u}$, it follows from Proposition \ref{log-esti-lem} that there exists a constant $L\geq n-1$ such that
\begin{equation*}
(n-1)\log |x|+U_0-C\leq\hat{u}(x)\leq L\log |x|,
\end{equation*}
and hence
\begin{equation}\label{A-u-inq}
U_0-L\log |x|\leq u(x)\leq C-(n-1)\log |x|
\end{equation}
for any $x\in\overline{\mathbb R^n_+}$ with $|x|$ large enough. By $u_{R}(x)=\frac{u(Rx)-\mu_0}{\log  R}$, we have 
\begin{equation}\label{A-uR-esti}
|u_{R}(x)|\leq\frac{L\log |Rx|+C}{\log  R}\leq\frac{L}{\log R_0}\log|x|+L+\frac{C}{\log R_0}
\end{equation}
for some constant $C>0$ independent of $R$, which indicates that $u_{R}$ is bounded in $L^{\infty}_{\operatorname{loc}}(\overline{\mathbb R^n_+}\setminus\{0\})$ uniformly in $R$. A direct computation shows that $u_{R}$ satisfies the equation 
\begin{eqnarray}\label{A-uR-eq}
&&\begin{cases}
-\Delta_nu_{R}=0, &\qquad \text{in} ~\mathbb R_+^n\setminus B_{\frac{R_0}{R}},\\
a(\nabla u_{R})\cdot\mathbf{n}=(\frac{R}{\log  R})^{n-1}e^{u(Rx)}, &\qquad \text{on}~\partial\mathbb R_+^n\setminus B_{\frac{R_0}{R}}.
\end{cases}
\end{eqnarray}
It follows from \eqref{A-eu-esti} that for any $x\in\overline{\mathbb R_+^n}\setminus\{0\}$,
\begin{equation}\label{A-euR-inq}
 (\frac{R}{\log R})^{n-1}e^{u(Rx)}\leq\frac{C}{(\log R)^{n-1}|x|^{n-1}}\leq\frac{C}{(\log R_0)^{n-1}|x|^{n-1}},   
\end{equation}
which implies $(\frac{R}{\log R})^{n-1}e^{u(Rx)}$ is bounded in $L^{\infty}_{\operatorname{loc}}(\overline{\mathbb R^n_+}\setminus\{0\})$, uniformly in $R$. For \eqref{A-uR-eq}, using the elliptic regularity theory in \cite{L1988, L1991}, we obtain that $u_{R}$ is uniformly bounded in $C^{1,\alpha}_{\operatorname{loc}}(\overline{\mathbb R_+^n}\setminus\{0\})$. By the Ascoli-Arzel\'{a}'s Theorem and a diagonal argument, we can find a sequence $R_{j}\rightarrow+\infty$ such that $u_{R_{j}}\rightarrow u_{\infty}$ in $C^{1}_{\operatorname{loc}}(\overline{\mathbb R_+^n}\setminus\{0\})$, where $u_{\infty}$ satisfies 
\begin{equation}\label{A-infty-eq}
\begin{cases}
\Delta_nu_{\infty} =0, &\qquad\text{in}~\mathbb R_+^n, \\
a(\nabla u_{\infty})\cdot\mathbf{n}=0, &\qquad\text{on}~\partial\mathbb R_+^n,
\end{cases}
\end{equation}
by combining \eqref{A-uR-eq} and the first inequality in \eqref{A-euR-inq}. It follows from the first inequality of \eqref{A-uR-esti} that $|u_{\infty}|\leq L$ for any $x\in\mathbb R_+^n$. By Lemma \ref{const-lem} with $\sigma=0$, $G(x)=u_{\infty}$ and by \eqref{A-u-inq}, we have $u_{\infty}\equiv-\beta$ for some positive constant $n-1\leq\beta\leq L$. Moreover, by $u_{R}\rightarrow -\beta$ in $C^{1}_{\operatorname{loc}}(\overline{\mathbb R_+^n}\setminus\{0\})$ as $R\rightarrow\infty$, we have 
\begin{equation}\label{A-u-limit}
\lim\limits_{x\in\overline{\mathbb R_+^n}, |x|\rightarrow+\infty}\frac{u(x)}{\log |x|}=-\beta.
\end{equation}

Next, we will prove that
\begin{equation*}
u(x)+\beta \log |x|\in L^{\infty}(\overline{\mathbb R_+^n}\setminus B_{R_0}).
\end{equation*}  
Write
$$I(R):=\inf\limits_{R_0\leq|x|\leq R,~x\in\mathbb R_+^n}\frac{u(x)-\mu_0}{\log|x|}.$$
Obviously, $I(R)$ is non-increasing and $I(R)\leq0$, and $u(x)\rightarrow-\infty$ as $|x|\rightarrow\infty$, hence there exists $\hat{R}>R_0$ such that $I(R)<0$ for all $R> \hat{R}$. Noting that 
\begin{equation*}
\begin{cases}
-\Delta_n(u-\mu_0)=-\Delta_n(I(R)\log|x|)=0,&\qquad\text{in }B_R^+\setminus B_{R_0}^+,\\
u-\mu_0\geq I(R)\log|x|,&\qquad\text{on }\partial B_{R_0}^+\cup\partial B_{R}^+,\\
a(\nabla(u-\mu_0))\cdot\mathbf{n}> a(\nabla(I(R)\log|x|))\cdot\mathbf{n}=0,&\qquad\text{on }\Sigma_R\setminus\Sigma_{R_0},
\end{cases}
\end{equation*}
the comparison principle in Lemma \ref{compa-prin-lem} and Serrin's strong maximum principle \cite{S1970} imply 
\begin{equation*}
I(R)=\inf\limits_{\partial B_R^+}\frac{u(x)-\mu_0}{\log R}
\end{equation*}for any $R>\hat{R}$.
In view of \eqref{A-u-limit}, we deduce $\lim_{R\rightarrow\infty}I(R)=-\beta$, which implies $\inf\limits_{\mathbb R_+^n\setminus B_{R_0}}\frac{u(x)-\mu_0}{\log|x|}=-\beta$, namely, 
\begin{equation*}
u(x)\geq C-\beta\log|x|,\qquad x\in\mathbb R_+^n\setminus B_{R_0},
\end{equation*}
for some constant $C>0$. Note that $\int_{\mathbb R_+^n}e^{\frac{n}{n-1}u}\mathrm dx<\infty$ implies $\beta>n-1$.

Now, it remains to prove $u(x)\leq C-\beta\log|x|$ in $\mathbb R_+^n\setminus B_{R_0}$. First, by $u_{\infty}=-\beta$ and $u_{R}\rightarrow u_{\infty}$ in $C^{1,\alpha}_{\operatorname{loc}}(\overline{\mathbb R_+^n}\setminus\{0\})$ as $R\rightarrow\infty$,  there exists $R(\varepsilon)\rightarrow+\infty$ as $\varepsilon\rightarrow0$ such that for any $R\geq R(\varepsilon)>R_0$,
\begin{equation*}
\big|\frac{u(Rx)-\mu_0}{\log R}+\beta\big|\leq\varepsilon,\qquad \text{in}~\mathbb R_+^n\setminus B_{\frac{R(\varepsilon)}{R}},
\end{equation*}
which implies
\begin{equation}\label{A-asy-esti}
u(x)-\mu_0\leq(-\beta+\varepsilon)\log|x|,\qquad\text{in}~\mathbb R_+^n\setminus B_{R(\varepsilon)}.
\end{equation}

Arguing as in Proposition 3.2 in \cite{DGL2024}(see also Proposition 3.1 in \cite{CL2024}), we construct suitable supersolutions to \eqref{I-hL-eq}. Taking $\varepsilon=\varepsilon_0:=\frac{\beta-(n-1)}{2}$ in \eqref{A-asy-esti}, we have
\begin{equation*}
e^u\leq\frac{C_0}{|x|^{\frac{\beta+n-1}{2}}},\qquad\text{in}~\mathbb R_+^n\setminus B_{R_0},
\end{equation*}
for some constant $C_0>0$.

For $\varepsilon\in(0,\beta-n+1)$ small enough, define
\begin{equation}\label{A-supsolution-defi}
\overline{u}_{\varepsilon}(x):=C_{1}^{\frac{1}{n-1}}\phi_{a_{\varepsilon},b}(|x|+\frac{t}{|x|^{\delta}}),\qquad x=(x',t)\in\overline{\mathbb R_+^n},
\end{equation}
where $C_{1}>0$ and
$$\phi_{a_{\varepsilon},b}(r):=-(\frac{\gamma-n}{2})^{\frac{1}{1-n}}\int_{R_{1}}^{r}\frac{(a_{\varepsilon}-s^{\frac{n-\gamma}{2}})^{\frac{1}{n-1}}}{s}\mathrm ds+b$$
is the solution of the ODE $r^{1-n}(r^{n-1}(-\phi')^{n-1})'=r^{-\frac{n+\gamma}{2}}$($r\geq R_{1}\gg1$).

The parameters are chosen as follows
\begin{equation*}
\begin{cases}
0<\frac{\gamma-n}{2}<\delta<\min\{1,\frac{\beta-(n-1)}{2}\},&\\
R_{1}>R_0,~(n-1)\log \frac{R_{1}}{2}>\mu_0,&\\
a_{\varepsilon}=\frac{(\beta-\varepsilon)^{n-1}(\gamma-n)}{2C_{1}}>\frac{(n-1)^{n-1}(\gamma-n)}{2C_{1}}\geq 2R_{1}^{\frac{n-\gamma}{2}},&\\
2C_{1}(\gamma-n)^{-1}(a_{\varepsilon}-R_{1}^{\frac{n-\gamma}{2}})R_{1}^{\frac{\beta-(n-1)}{2}-\delta}&\\
\quad>\big[(n-1)^{n-1}-2C_{1}(\gamma-n)^{-1}R_{1}^{\frac{n-\gamma}{2}}\big]R_{1}^{\frac{\beta-(n-1)}{2}-\delta}\geq C_0,&\\
C_{1}^{\frac{1}{n-1}}b\geq \max\big\{0,\sup\limits_{\partial B_{R_1}^+}u+\beta\log2\big\}.&\\
\end{cases}
\end{equation*}
Note that for any $r\geq R_{1}$, 
$$\phi_{a_{\varepsilon},b}'(r)=-(\frac{\gamma-n}{2})^{\frac{1}{1-n}}\frac{(a_{\varepsilon}-r^{\frac{n-\gamma}{2}})^{\frac{1}{n-1}}}{r}\leq0,$$
\begin{equation}\label{A-phi-eq}
\big[\big(-\phi_{a_{\varepsilon},b}'\big)^{n-1}\big]'+(n-1)\frac{\big(-\phi_{a_{\varepsilon},b}'\big)^{n-1}}{r}=r^{-\frac{n+\gamma}{2}}   
\end{equation}
and 
\begin{equation}\label{A-phi-inq}
(\frac{\gamma-n}{2})^{\frac{1}{1-n}}\frac{(a_{\varepsilon}-r^{\frac{n-\gamma}{2}})^{\frac{1}{n-1}}}{r}\leq-\phi_{a_{\varepsilon},b}'(r)\leq(\frac{\gamma-n}{2})^{\frac{1}{1-n}}\frac{(a_{\varepsilon})^{\frac{1}{n-1}}}{r}.
\end{equation}
We can verify that (see Appendix \ref{Section AP-B} for the details)
\begin{equation}\label{A-supsolution-inq}
\begin{cases}
-\Delta_n\overline{u}_{\varepsilon}\geq C|x|^{-\frac{n+\gamma}{2}}\geq0,&\qquad\text{in }B_R^+\setminus B_{R_{1}}^+,\\
a(\nabla \overline{u}_{\varepsilon})\cdot\mathbf{n}\geq \frac{C_0}{|x|^{\frac{n-1+\beta}{2}}}\geq e^u,&\qquad\text{on } \Sigma_R\setminus \Sigma_{R_1}.\\
\end{cases}
\end{equation}
Recalling \eqref{A-phi-inq} and the definition of $\overline{u}_{\varepsilon}$ and integrating \eqref{A-phi-inq} from $R_{1}$ to $r$, where $r=|x|+\frac{t}{|x|^{\delta}}$, we obtain for any $x\in\mathbb R_+^n\setminus B_{R_1}$,
\begin{equation}\label{A-supsolution-esti1}
\overline{u}_{\varepsilon}\geq(-\beta+\varepsilon)\log\frac{|x|+\frac{t}{|x|^{\delta}}}{R_{1}}+C_{1}^{\frac{1}{n-1}}b\geq(-\beta+\varepsilon)\log\frac{2|x|}{R_1}+C_{1}^{\frac{1}{n-1}}b
\end{equation}
and
\begin{equation}\label{A-supsolution-esti2}
\overline{u}_{\varepsilon}\leq(-\beta+\varepsilon)\log\frac{|x|+\frac{t}{|x|^{\delta}}}{R_{1}}+C_{1}^{\frac{1}{n-1}}b+C_{2}\leq(-\beta+\varepsilon)\log\frac{|x|}{R_{1}}+C_{1}^{\frac{1}{n-1}}b+C_{2},
\end{equation}
where 
\begin{eqnarray*}
C_{2}&=&\frac{1}{n-1}C_{1}^{\frac{1}{n-1}}(\frac{2}{\gamma-n})^{\frac{n}{n-1}}R_{1}^{\frac{n-\gamma}{2(n-1)}}\nonumber\\
&\geq& C_{1}^{\frac{1}{n-1}}(\frac{\gamma-n}{2})^{\frac{1}{1-n}}\int_{R_{1}}^{r}\frac{a_{\varepsilon}^{\frac{1}{n-1}}-(a_{\varepsilon}-r^{\frac{n-\gamma}{2}})^{\frac{1}{n-1}}}{r}\mathrm dr.  
\end{eqnarray*}

Together, \eqref{A-asy-esti} and \eqref{A-supsolution-esti1} imply
\begin{equation}\label{A-boundary-inq2}
\overline{u}_{\varepsilon}\geq u,\qquad\text{on}~\partial B_R^+\cup \partial B_{R_{1}}^+.
\end{equation}
From \eqref{A-supsolution-inq} and \eqref{A-boundary-inq2}, an application of the comparison principle in Lemma \ref{compa-prin-lem} to $u$ and $\overline{u}_{\varepsilon}$ yields
\begin{equation}\label{A-compa-inq}
u\leq\overline{u}_{\varepsilon},\qquad\text{in }B_R^+\setminus B_{R_1}^+.
\end{equation}
Combining \eqref{A-compa-inq} and \eqref{A-supsolution-esti2} and letting $R\rightarrow\infty$, $\varepsilon\rightarrow 0$, it yields 
$$u(x)\leq C-\beta\log|x|, \qquad\text{in}~\mathbb R_+^n\setminus B_{R_1}.$$ 
By the continuity of $u$, we establish the validity of \eqref{A-u-esti} .

Set $F(x):=u(x)+\beta \log |x|$, it follows that
\begin{equation}\label{A-uF-eq}
u(x)=-\beta \log |x|+F(x),\qquad x\in\overline{\mathbb R_+^n}\setminus \{0\}.
\end{equation}
Using \eqref{A-u-esti}, we get $F(x)\in L^{\infty}(\overline{\mathbb R_+^n}\setminus B_{R_0})$. Consider the sequence of functions
\begin{equation*}
\hat{u}_{R_{j}}(x):=u(R_{j}x)+\beta \log R_{j},\qquad x\in \overline{\mathbb R_+^n}\setminus \{0\},
\end{equation*}
where $R_{j}\rightarrow+\infty$ as $j\rightarrow\infty$. The equation \eqref{A-uF-eq} implies 
\begin{equation}\label{A-ujF-eq}
\hat{u}_{R_{j}}(x)=-\beta \log |x|+F_{R_{j}}(x),\qquad x\in \overline{\mathbb R_+^n}\setminus \{0\},
\end{equation}
where $F_{R_{j}}(x):=F(R_{j}x)$. Also note that $\hat{u}_{R_{j}}(x)$ is the solution of 
\begin{equation}\label{A-uRj-eq}
\begin{cases}
-\Delta_n\hat{u}_{R_{j}}=0, &\qquad \text{in} ~\mathbb R_+^n\setminus B_{\frac{R_0}{R_{j}}},\\
a(\nabla \hat{u}_{R_{j}})\cdot\mathbf{n}=R_{j}^{n-1-\beta}\frac{e^{F_{R_{j}}}}{|x|^{\beta}}, &\qquad \text{on}~\partial\mathbb R_+^n\setminus B_{\frac{R_0}{R_{j}}}.
\end{cases}
\end{equation}
It can be observed that $\hat{u}_{R_{j}}$ and $R_{j}^{n-1-\beta}\frac{e^{F_{R_{j}}}}{|x|^{\beta}}$ are bounded in $L_{\operatorname{loc}}^{\infty}(\overline{\mathbb R_+^n}\setminus\{0\})$, uniformly in $j$. By the elliptic regularity estimates \cite{L1988, L1991}, we know that, up to a subsequence, $\hat{u}_{R_{j}}\rightarrow \hat{u}_{\infty}$ in $C^{1,\alpha}_{\operatorname{loc}}(\overline{\mathbb R_+^n}\setminus\{0\})$ as $j\rightarrow\infty$ for some $\hat{u}_{\infty}$, which satisfies
\begin{equation*}
\begin{cases}
-\Delta_n\hat{u}_{\infty}=0, &\qquad\text{in}~\mathbb R_+^n,\\
a(\nabla \hat{u}_{\infty})\cdot\mathbf{n}=0, &\qquad\text{on}~\partial\mathbb R_+^n,
\end{cases}
\end{equation*}
by \eqref{A-uRj-eq} and $\beta>n-1$. \eqref{A-ujF-eq} entails that
\begin{equation*}
\hat{u}_{\infty}(x)=-\beta \log |x|+F_{\infty}(x),\qquad x\in \overline{\mathbb R_+^n},
\end{equation*}
where $F_{\infty}$ is the limit of $F_{R_{j}}$ in $C^{1}_{\operatorname{loc}}(\overline{\mathbb R_+^n}\setminus\{0\})$. The uniform boundedness of $F_{R_{j}}$ ensures that $F_{\infty}\in L^{\infty}(\overline{\mathbb R_+^n}\setminus\{0\})$. Applying Lemma \ref{const-lem} to $F_{\infty}$, we conclude that $F_{\infty}$ is a constant function.
It follows from \eqref{A-uF-eq} and $F_{R_{j}}\rightarrow F_{\infty}$ in $C^{1}_{\operatorname{loc}}(\overline{\mathbb R_+^n}\setminus\{0\})$ that 
\begin{eqnarray*}
\sup\limits_{|x|=R_{j},x\in \overline{\mathbb R_+^n}}|x||\nabla(u(x)+\beta \log |x|)|&=&\sup\limits_{|y|=1, y\in \overline{\mathbb R_+^n}}|y||\nabla F_{R_{j}}(y)| \\
&\rightarrow& \sup\limits_{|y|=1, y\in \overline{\mathbb R_+^n}}|y||\nabla F_{\infty}(y)|=0,
\end{eqnarray*}
as $j\rightarrow\infty$ which implies \eqref{A-gradu-esti}. The proof is complete.  
\qed

As an immediate consequence of Lemma 4.2 in \cite{CL2024} or Theorem 4.2 in \cite{IT2012}, we have the following Pohozaev identity.
\begin{lemma}\label{Poho-lem}
Let $\Omega\subset\mathbb R^n$ be a bounded open set and $\Omega\cap\mathbb R_+^n$ be with Lipschitz boundary and $u\in C^{1}(\overline{\Omega\cap\mathbb R^n_+})$ is a solution of 
\begin{equation}\label{A-ufg-eq}
\begin{cases}
-\Delta_{p} u= f(x),&\qquad\text{in }\Omega\cap\mathbb R_+^n,\\
|\nabla u|^{p-2}\frac{\partial u}{\partial t}=-g(x), &\qquad\text{on }\Omega\cap\partial\mathbb R_+^n,\\
\end{cases}
\end{equation}
where $f\in L^{1}(\Omega\cap\mathbb R_+^n)\cap L_{\operatorname{loc}}^{\infty}(\Omega\cap\mathbb R_+^n)$, $g\in L^{1}(\Omega\cap\partial\mathbb R_+^n)$ and $1<p<\infty$. Then the following Pohozaev identity holds for any $y\in\mathbb R^n$: 
\begin{eqnarray}\label{A-Pohozaev-iden}
&&\frac{p-n}{p}\int_{\Omega\cap\mathbb R_+^n}|\nabla u|^p\mathrm dx-\int_{\Omega\cap\mathbb R_+^n}f(x)\langle x-y,\nabla u\rangle \mathrm dx\nonumber\\ 
&=&\int_{\partial\Omega\cap\mathbb R_+^n}|\nabla u|^{p-2}\langle\nabla u,\mathbf{n}\rangle\langle x-y,\nabla u\rangle \mathrm d\mathcal H^{n-1}(x)+\int_{\Omega\cap\partial\mathbb R_+^n}g(x)\langle x-y,\nabla u\rangle \mathrm dx' \nonumber\\
&&-\frac{1}{p}\int_{\partial(\Omega\cap\mathbb R_+^n)}|\nabla u|^p\langle x-y,\mathbf{n}\rangle \mathrm d\mathcal H^{n-1}(x).
\end{eqnarray}
\end{lemma}

Inspired by \cite{CL2024}, \cite{DGL2024} and \cite{E2018}, using Proposition \ref{1order-esti-prop} and the Pohozaev identity, we can obtain the precise value of $\beta$ in Proposition \ref{1order-esti-prop}.
\begin{proposition}\label{beta-value-lem}
Let $\beta$ be the constant in the sharp asymptotic estimate on $u$ and $\nabla u$ in Proposition \ref{1order-esti-prop}. Then $\beta=n$.    
\end{proposition}

\begin{proof}
Let $u$ be the solution of \eqref{I-hL-eq}. We shall need the following  two conclusions to complete the proof, that is,
\begin{equation}\label{A-beta-iden1}
\beta=\big(\frac{2\int_{\partial\mathbb R^n_+}e^u\mathrm dx'}{n\omega_n}\big)^{\frac{1}{n-1}}
\end{equation}
and
\begin{equation}\label{A-eu-value}
\int_{\partial\mathbb R^n_+}e^u\mathrm dx'=\frac{1}{2}n^n\omega_n.
\end{equation}
\noindent\textbf{Step 1.} Proof of \eqref{A-beta-iden1}. By \eqref{A-gradu-esti} in Proposition \ref{1order-esti-prop}, we have 
\begin{equation}\label{A-gu-esti}
\nabla u=-\beta |x|^{-2}x+o(|x|^{-1})=-\beta R^{-2}x+o(R^{-1})
\end{equation}
and 
\begin{equation}\label{A-gul-esti}
|\nabla u|=\beta|x|^{-1}+o(|x|^{-1})=\beta R^{-1}+o(R^{-1}),
\end{equation}
uniformly for $x\in \partial B_R^+$, as $R\rightarrow\infty$. Using the divergence theorem, it yields 
\begin{eqnarray}\label{A-div-esti}
\int_{\Sigma_R}e^u\mathrm dx' &=& \int_{\partial B_{R}^+}|\nabla u|^{n-2}\langle\nabla u,-\mathbf{n}\rangle \mathrm d\mathcal H^{n-1}\nonumber\\
&=& \int_{\partial B_R^+}|\nabla u|^{n-2}\langle\nabla u,-\frac{x}{|x|}\rangle \mathrm d\mathcal H^{n-1} \nonumber\\
&=&\frac{R}{\beta}\int_{\partial B_{R}^+}|\nabla u|^{n-2}\langle\nabla u,\nabla u+o(R^{-1})\rangle \mathrm d\mathcal H^{n-1} \nonumber\\
&=&\int_{\partial B_{R}^+}(\beta R^{-1}+o(R^{-1}))^{n-1}(1+o(1))\mathrm d\mathcal H^{n-1} \nonumber\\
&=&\beta^{n-1}\int_{\partial B_{1}^+}(1+o(1))^n \mathrm d\mathcal H^{n-1}.
\end{eqnarray}
Letting $R\rightarrow\infty$, we obtain
\begin{equation*}
\int_{\partial\mathbb R^n_+}e^u\mathrm dx'=\frac{1}{2}\beta^{n-1}n\omega_n,
\end{equation*}
that is, 
\begin{equation*}
\beta=\big(\frac{2\int_{\partial\mathbb R^n_+}e^u\mathrm dx'}{n\omega_n}\big)^{\frac{1}{n-1}}.
\end{equation*}
This completes the proof of \eqref{A-beta-iden1}.
   
\noindent \textbf{Step 2.} Proof of \eqref{A-eu-value}. Using \eqref{A-Pohozaev-iden} with $f=0$, $g=e^u$, $\Omega=B_{R}$, $y=0$ and $p=n$, we deduce
\begin{eqnarray}\label{A-sP-iden}
0 &=& \int_{\partial B_{R}^+}|\nabla u|^{n-2}\langle \nabla u,\mathbf{n}\rangle\langle x,\nabla u\rangle \mathrm d\mathcal H^{n-1}+ \int_{\Sigma_R}e^u\langle x,\nabla u\rangle \mathrm dx' \nonumber\\
&&-\frac{1}{n}\int_{\partial B_{R}^+}|\nabla u|^n\langle x,\mathbf{n}\rangle \mathrm d\mathcal H^{n-1}.
\end{eqnarray} 
Recalling Proposition \ref{1order-esti-prop}, it yields 
\begin{equation}\label{A-xu-esti}
\langle x,\nabla u\rangle=-\beta+o(1),
\end{equation}
uniformly for $x\in\partial B_R^+$, as $R\rightarrow\infty$. Using \eqref{A-gu-esti}, \eqref{A-gul-esti} and \eqref{A-xu-esti}, as done in \eqref{A-div-esti}, it holds 
\begin{eqnarray}\label{A-lim-eq}
&&\int_{\partial B_R^+}|\nabla u|^{n-2}\langle \nabla u,\mathbf{n}\rangle\langle x,\nabla u\rangle \mathrm d\mathcal H^{n-1}-\frac{1}{n}\int_{\partial B_R^+}|\nabla u|^n\langle x,\mathbf{n}\rangle \mathrm d\mathcal H^{n-1}\nonumber\\
&&\rightarrow \frac{1}{2}\beta^n(n-1)\omega_n
\end{eqnarray}
as $R\rightarrow\infty$. The divergence theorem implies 
\begin{eqnarray}\label{dive-thm-equ}
\int_{\partial B_{R}\cap\partial\mathbb R_+^n}e^u\langle x,\mathbf{n}\rangle\mathrm d\mathcal{H}^{n-2} 
&=& \int_{ \Sigma_R}\operatorname{div}(e^ux) \mathrm dx' \nonumber\\
&=&\int_{ \Sigma_R}e^u\langle\nabla u,x\rangle \mathrm dx'+(n-1)\int_{ \Sigma_R}e^u \mathrm dx'. ~~~~~~~~
\end{eqnarray} 
It follows from $\int_{\mathbb R_+^n}e^{\frac{n}{n-1}u}\mathrm dx$ $<\infty$ and \eqref{A-u-esti} that
\begin{equation}\label{A-eu-inq}
\int_{\partial B_{R}\cap\partial\mathbb R_+^n}e^u\langle x,\mathbf{n}\rangle\mathrm d\mathcal{H}^{n-2}\leq CR^{n-1-\beta}\rightarrow0, \qquad \text{as}~R\rightarrow\infty.
\end{equation}
From \eqref{A-sP-iden} and \eqref{A-lim-eq}-\eqref{A-eu-inq}, letting $R\rightarrow\infty$, one has
\begin{equation}\label{A-beta-iden2}
\int_{\partial\mathbb R_+^n}e^u \mathrm dx'=\frac{1}{2}\beta^n\omega_n.
\end{equation}
On the other hand, from \eqref{A-beta-iden1} and \eqref{A-beta-iden2}, it holds
\begin{equation*}
\int_{\partial\mathbb R^n_+}e^u\mathrm dx'=\frac{1}{2}n^n\omega_n.
\end{equation*}
This proves \eqref{A-eu-value}.

Finally, combining \eqref{A-beta-iden1} with \eqref{A-eu-value} yields $\beta=n$.         
\end{proof}

\subsection{Asymptotic integral estimate on second order derivatives}
Following the approach in \cite{CFR2020, DGL2024, Z2024}, that is, by using the approximate method and the Caccioppoli type inequality, we can prove the following second-order regularity and asymptotic integral estimate on the second-order derivatives of weak solutions to \eqref{I-hL-eq}.
\begin{proposition}\label{2order-esti-prop}
Let $u$ be a solution to \eqref{I-hL-eq} satisfying \eqref{I-fm-assu}. Then $a(\nabla u)\in W_{\operatorname{loc}}^{1,2}(\overline{\mathbb R^n_+} )$, and for any $\gamma\in \mathbb R$ the following asymptotic estimate holds:
\begin{equation}\label{A-2order-esti}
\int_{B_{2R}^+\setminus B_R^+}|\nabla[a(\nabla u)]|^2e^{\gamma u}\mathrm dx\leq CR^{-(\gamma+1)n},\qquad\forall R>1,
\end{equation}
where $C$ is a positive constant independent of $R$. 
\end{proposition}

We also need a linear algebra lemma as follows.

\begin{lemma}[\cite{AKM2018} Lemma 4.5]\label{matrix-lem}
Let the $n$-th order matrix $A$ be symmetric with positive eigenvalues and let $\lambda_{\min}$ and $\lambda_{\max}$ be its smallest and largest eigenvalues, respectively. Let $B$ be a $n$-th order symmetric matrix, then
$$\operatorname{tr}[(AB)^2]\leq\operatorname{tr}[AB(AB)^T]\leq n\left(\frac{\lambda_{\max}}{\lambda_{\min}}\right)^2\operatorname{tr}[(AB)^2].$$
\end{lemma}

\noindent\textbf{Proof of Proposition \ref{2order-esti-prop}.}
Let $\{\phi_{j}\}$ be a family of radially symmetric smooth mollifiers and
\begin{equation*}
a(x):=|x|^{n-2}x,\quad a^{j}(x):=(a\ast \phi_j)(x),
\end{equation*}
\begin{equation*}
b(x):=|x|^n,\quad b^j(x):=(b\ast \phi_j)(x)
\end{equation*}
for any $x\in\mathbb R^n$.
By Lemma 2.4 in \cite{FF1997} and its proof therein, we deduce that $a^{j}$ satisfies 
\begin{equation}\label{A-ak-inq}
\langle\nabla a^{j}(x)\xi,\xi\rangle\geq\frac{1}{\lambda}(|x|^2+s_{j}^2)^{\frac{n-2}{2}}|\xi|^2~\text{and}~|\nabla a^{j}(x)|\leq\lambda(|x|^2+s_{j}^2)^{\frac{n-2}{2}},
\end{equation}
for $\forall x,\xi\in\mathbb R^n$, with $s_{j}\neq0,~ s_{j}\to0$ as $j\to\infty$ and for some $\lambda>0$. 

Assume that $u_j^R\in W^{1,n}(B_R^+)$ satisfies
\begin{equation}\label{A-j-eq}
\begin{cases}
\operatorname{div}(a^j(\nabla u_j^R))=0, &\qquad\text{in}~B_R^+,\\
a^j(\nabla u_j^R)\cdot\mathbf{n}=e^u, &\qquad\text{on}~\Sigma_R,\\
u_j^R=u,&\qquad\text{on }\partial B_R^+.
\end{cases}
\end{equation}
It corresponds to the following variational problem
$$\min\limits_{v\in \mathcal X_{n.R}}\{\int_{B_R^+}\frac1 nb^j(\nabla v)\mathrm dx-\int_{\Sigma_R}e^uv\mathrm dx'\},$$
where 
$$\mathcal X_{n.R}:=\{v\in W^{1,n}(B_R^+):v|_{\partial B_R^+}=u\},$$
which implies that the solution $u_j^R$ exists and is unique. Using the Moser iteration, we have $\|u_j^R\|_{L^\infty(B_R^+)}\leq C(R)$. By the elliptic regularity theory \cite{L1988,L1991}, $\|u_j^R\|_{C^{1,\alpha}(B_R^+)}\leq C(R)$. Hence, $u_j^R\rightarrow u_0$ in $C^1(B_R^+)$ as $j\rightarrow+\infty$, where $u_0$ satisfies
\begin{equation*}
\begin{cases}
 \operatorname{div}(a(\nabla u_0))=0, &\qquad\text{in}~B_R^+,\\
a(\nabla u_0)\cdot\mathbf{n}=e^u, &\qquad\text{on}~\Sigma_R,\\
u_0=u,&\qquad\text{on }\partial B_R^+.   
\end{cases}
\end{equation*}
It follows from the uniqueness that $u=u_0$, hence $u_j^R\rightarrow u=u_0$ in $C^1(B_R^+)$ as $j\rightarrow+\infty$. By a diagonal process, there exist diagonal subsequences $\{u_{j_k}^{R_k}\}_{k\in\mathbb N}=:\{u_k\}_{k\in\mathbb N}$ and $\{a^{j_k}\}_{k\in\mathbb N}$ that we still write as $\{a^k\}_{k\in\mathbb N}$, such that $u_k\rightarrow  u$ in $C^1(B_{R_k}^+)$ as $k\rightarrow+\infty$ and
\begin{equation}\label{A-k-eq}
\begin{cases}
\operatorname{div}(a^k(\nabla u_k))=0, &\qquad\text{in}~B_{R_k}^+,\\
a^k(\nabla u_k)\cdot\mathbf{n}=e^u, &\qquad\text{on}~\Sigma_{R_k},\\
u_k=u, &\qquad\text{on }\partial B_{R_k}^+,
\end{cases}
\end{equation}
where $j_k\rightarrow+\infty$ and $R_k\rightarrow+\infty$ as $k\rightarrow+\infty$. By the uniform elliptic theory, $u_k\in C^{2,\alpha}(\overline{\mathbb R_+^n\cap B_{\frac{R_k}{2}}})\cap C^{3,\alpha}_{\operatorname{loc}}(\mathbb R_+^n\cap B_{R_k})$.

Using an analogous argument to that of Proposition 2.7 in \cite{Z2024} and recalling Proposition \ref{1order-esti-prop}, we can verify that $a(\nabla u)=|\nabla u|^{n-2}\nabla u\in W_{\operatorname{loc}}^{1,2}(\overline{\mathbb R^n_+} )$ and the validity of \eqref{A-2order-esti}. The details of the proof are as follows.

Choose a $C^2$ domain $U$ satisfying $B_{4R}\cap\mathbb R_+^n\subset U\subset B_{5R}\cap\mathbb R_+^n$, where $R>1$. For $\delta>0$ small, define the set $U_{\delta}:=\{x\in U:\operatorname{dist}(x,\partial U)>\delta\}$ and $d(x):=\operatorname{dist}(x,\partial U)$. For $k$ large enough, $R_k>5R$.

Firstly, testing \eqref{A-k-eq} with $\psi\in C_{c}^{0,1}(U)$ yields 
\begin{equation}\label{A-kinte-eq}
\int_{U}a^{k}(\nabla u_{k})\cdot\nabla \psi\mathrm dx=0.
\end{equation}
Let $\psi=\partial_{m}\varphi\xi_{\delta}$ in \eqref{A-kinte-eq} where $\varphi\in C^2(U)\cap C_{c}^{1}(B_{4R}\cap \overline{\mathbb R_+^n})$ and $\xi_{\delta}(x)=\gamma(d(x))$,
\begin{equation*}
\gamma(t)=
\begin{cases}
0,\quad&t\in[0,\delta),\\
\frac{t-\delta}{\delta},\quad&t\in[\delta,2\delta),\\
1,\quad
&t\in[2\delta,\infty).
\end{cases}
\end{equation*}
It can be observed that $\xi_\delta\in C_c^{0,1}(U)$ and
$$\xi_\delta=1\text{ in }U_{2\delta},~\xi_\delta=0\text{ in }U\setminus U_{\delta},~\nabla\xi_\delta=-\frac{1}{\delta}\mathbf{n}(y(x))\text{ in }U_\delta\setminus U_{2\delta},$$
where $y=y(x)\in \partial U$ is the unique projection of $x\in U\setminus U_{2\delta}$ on $\partial U$.

Using integration by parts in \eqref{A-kinte-eq}, one has
\begin{eqnarray*}\label{A-varphi-inq1}
0&=&\sum\limits_{i=1}^n\int_{U}a^{k}_{i}(\nabla u_{k})\partial_{i} (\partial_{m}\varphi\xi_{\delta})\mathrm dx\nonumber\\
&=&-\sum\limits_{i=1}^n\int_{U}\partial_{i}\varphi\partial_{m}[a^{k}_{i}(\nabla u_{k})\xi_{\delta}]\mathrm dx+\sum\limits_{i=1}^n\int_{U}a^{k}_{i}(\nabla u_{k})\partial_{m}\varphi\partial_{i}\xi_{\delta}\mathrm dx,
\end{eqnarray*}
that is,
\begin{eqnarray}\label{A-xi-eq}
&&\sum\limits_{i=1}^n\int_{U}\partial_{i}\varphi\partial_{m}[a^{k}_{i}(\nabla u_{k})]\xi_{\delta}\mathrm dx\nonumber\\
&=&-\sum\limits_{i=1}^n\int_{U}\partial_{i}\varphi a^{k}_{i}(\nabla u_{k})\partial_{m}\xi_{\delta}\mathrm dx+\sum\limits_{i=1}^n\int_{U}a^{k}_{i}(\nabla u_{k})\partial_{m}\varphi\partial_{i}\xi_{\delta}\mathrm dx.
\end{eqnarray}
Letting $\delta\rightarrow0$ into \eqref{A-xi-eq}, by the co-area formula, it yields
\begin{eqnarray}\label{A-delta0-eq}
&&\sum\limits_{i=1}^n\int_{U}\partial_{i}\varphi\partial_{m}[a^{k}_{i}(\nabla u_{k})]\mathrm dx\nonumber\\ &=&\sum\limits_{i=1}^n\big(\int_{\partial\mathbb R_+^n}\partial_{i}\varphi a^{k}_{i}(\nabla u_{k})\mathbf{n}^{m}\mathrm dx'-\int_{\partial\mathbb R_+^n}a^{k}_{i}(\nabla u_{k})\partial_{m}\varphi\mathbf{n}^{i}\mathrm dx'\big),
\end{eqnarray}
where $\mathbf{n}=(\mathbf{n}^1,\mathbf{n}^2,\cdots,\mathbf{n}^n)=(0,\cdots,0,-1)$. Set $\varphi=a^{k}_{m}(\nabla u_{k})\rho^2$ in \eqref{A-delta0-eq} where $\rho\in C_{c}^{2}(B_{4R})$. It follows from \eqref{A-delta0-eq} that if $m\neq n$
\begin{eqnarray}\label{A-mn-eq1}
&&\sum\limits_{i=1}^n\int_{U}\partial_{m}[a^{k}_{i}(\nabla u_{k})]\partial_{i}[a^{k}_{m}(\nabla u_{k})\rho^2]\mathrm dx=-\sum\limits_{i=1}^n\int_{\partial\mathbb R_+^n}a^{k}_{i}(\nabla u_{k})\partial_{m}\varphi\mathbf{n}^{i}\mathrm dx'\nonumber\\
&=&-\int_{\partial\mathbb R_+^n}e^u\partial_{m}\varphi\mathrm dx'=\int_{\partial\mathbb R_+^n}\partial_{m}(e^u)a^{k}_{m}(\nabla u_{k})\rho^2\mathrm dx'.
\end{eqnarray}
For the case $m=n$, we have
\begin{eqnarray}\label{A-mn-eq2}
&&\sum\limits_{i=1}^n\int_{U}\partial_n[a^{k}_{i}(\nabla u_{k})]\partial_{i}[a^{k}_n(\nabla u_{k})\rho^2]\mathrm dx\nonumber\\
&=&\sum\limits_{i=1}^n\int_{\partial\mathbb R_+^n}\partial_{i}\varphi a^{k}_{i}(\nabla u_{k})\mathbf{n}^n\mathrm dx'-\sum\limits_{i=1}^n\int_{\partial\mathbb R_+^n}a^{k}_{i}(\nabla u_{k})\partial_n\varphi\mathbf{n}^{i}\mathrm dx'\nonumber\\
&=&-\sum_{i=1}^{n-1}\int_{\partial\mathbb R_+^n}\partial_{i}\varphi a^{k}_{i}(\nabla u_{k})\mathrm dx'=-\sum_{i=1}^{n-1}\int_{\partial\mathbb R_+^n}\partial_{i}[a^{k}_n(\nabla u_{k})\rho^2] a^{k}_{i}(\nabla u_{k})\mathrm dx'\nonumber\\
&=&\sum_{i=1}^{n-1}\int_{\partial\mathbb R_+^n}\partial_{i}(e^u\rho^2) a^{k}_{i}(\nabla u_{k})\mathrm dx',
\end{eqnarray}
where the last equality holds by the fact $a^{k}_n(\nabla u_{k})=-e^u$ on $\partial\mathbb R_+^n$. Combining \eqref{A-mn-eq1} and \eqref{A-mn-eq2} yields
\begin{eqnarray*}
&&\sum\limits_{i=1}^n\int_{U}\partial_{m}[a^{k}_{i}(\nabla u_{k})]\partial_{i}[a^{k}_{m}(\nabla u_{k})\rho^2]\mathrm dx\nonumber\\
&\leq&C\int_{\partial\mathbb R_+^n}e^u|\nabla u||a^{k}(\nabla u_{k})|\rho^2\mathrm dx'+C\int_{\partial\mathbb R_+^n}e^u|a^{k}(\nabla u_{k})||\nabla \rho|\rho\mathrm dx'.
\end{eqnarray*}
We arrive at
\begin{eqnarray}\label{A-rho-inq}
&&\sum\limits_{i,m=1}^n\int_{U}\partial_{m}[a^{k}_{i}(\nabla u_{k})]\partial_{i}[a^{k}_{m}(\nabla u_{k})]\rho^2\mathrm dx\nonumber\\
&\leq&C\int_{U}|\nabla(a^{k}(\nabla u_{k}))||a(\nabla u_{k})||\nabla \rho|\rho\mathrm dx+C\int_{\partial\mathbb R_+^n}e^u|\nabla u||a^{k}(\nabla u_{k})|\rho^2\mathrm dx'\nonumber\\
&&+C\int_{\partial\mathbb R_+^n}e^u|a^{k}(\nabla u_{k})||\nabla \rho|\rho\mathrm dx'.
\end{eqnarray}
By \eqref{A-ak-inq} and Lemma \ref{matrix-lem}, it holds
\begin{equation*}
\big|\nabla [a^{k}(\nabla u_{k})]\big|^2\leq C\sum\limits_{i,m=1}^n\partial_{m}[a^{k}_{i}(\nabla u_{k})]\partial_{i}[a^{k}_{m}(\nabla u_{k})].
\end{equation*}
Using Young's inequality for \eqref{A-rho-inq}, we obtain the following Caccioppoli-type inequality
\begin{eqnarray}\label{A-kCacci-inq}
&&\int_{U}\big|\nabla [a^{k}(\nabla u_{k})]\big|^2\rho^2\mathrm dx\nonumber\\
&\leq&C\big(\int_{U}|a^{k}(\nabla u_{k})|^2|\nabla\rho|^2\mathrm dx+\int_{\partial\mathbb R_+^n}e^u|\nabla u||a^{k}(\nabla u_{k})|\rho^2\mathrm dx'\nonumber\\
&&+\int_{\partial\mathbb R_+^n}e^u|a^{k}(\nabla u_{k})||\nabla \rho|\rho\mathrm dx'\big),
\end{eqnarray}
which implies that $a^{k}(\nabla u_{k})\in W_{\operatorname{loc}}^{1,2}(B_{R_k}\cap\overline{\mathbb R_+^n})$ and $\{a^{k}(\nabla u_{k})\}_{k\in \mathbb N}$ is uniformly bounded in $W_{\operatorname{loc}}^{1,2}(B_{R_k}\cap\overline{\mathbb R_+^n})$. Up to a subsequence, we have $a^{k}(\nabla u_{k})\rightharpoonup a(\nabla u)$ in $W_{\operatorname{loc}}^{1,2}(\overline{\mathbb R_+^n})$ as $k\rightarrow+\infty$ and $a(\nabla u)\in W_{\operatorname{loc}}^{1,2}(\overline{\mathbb R_+^n})$.

Choose $\rho=e^{\frac{\gamma}{2}u_k}\eta$ where $\gamma\in\mathbb R$, $\eta\in C_c^2(B_{4R})$ and
\begin{equation*}
\eta=
\begin{cases}
1,\quad&\text{in }B_{2R}\setminus B_{R},\\
0,\quad&\text{in }(B_{3R}\setminus B_{\frac{R}{2}})^c,\\
\end{cases}
\end{equation*} with $0\leq\eta\leq1,~|\nabla \eta|\leq\frac{C_n}{R}$ in $\mathbb R^n$. Letting $k\rightarrow\infty$ into \eqref{A-kCacci-inq}, one has
\begin{eqnarray}\label{A-eta-inq}
&&\int_{B_{2R}^+\setminus B_{R}^+}|\nabla [a(\nabla u)]|^2e^{\gamma u}\mathrm dx\nonumber\\
&\leq&C\int_{\Sigma_{3R}\setminus\Sigma_{\frac R2}}\big[e^u|\nabla u||a(\nabla u)|(e^{\frac{\gamma}{2}u}\eta)^2+e^u|a(\nabla u)|e^{\frac{\gamma}{2}u}\eta|\nabla (e^{\frac{\gamma}{2}u}\eta)|\big]\mathrm dx'\nonumber\\
&&+ C\int_{B_{3R}^+\setminus B_{\frac R2}^+}|a(\nabla u)|^2|\nabla(e^{\frac{\gamma}{2}u}\eta)|^2\mathrm dx\nonumber\\
&\leq& C\int_{B_{3R}^+\setminus B_{\frac R2}^+}\big(|\nabla u|^{2n}e^{\gamma u}\eta^2+|\nabla u|^{2(n-1)}e^{\gamma u}|\nabla \eta|^2\big)\mathrm dx\nonumber\\
&&+C\int_{\Sigma_{3R}\setminus\Sigma_{\frac R2}}\big(|\nabla u|^ne^{(\gamma+1)u}\eta^2+e^{(\gamma+1)u}|\nabla u|^{n-1}|\nabla \eta|\eta\big)\mathrm dx'.
\end{eqnarray}
According to Proposition \ref{1order-esti-prop} and \eqref{A-eta-inq}, we get
\begin{equation*}
\int_{B_{2R}^+\setminus B_R^+}|\nabla [a(\nabla u)]|^2e^{\gamma u}\mathrm dx\leq CR^{-(\gamma+1)n},\qquad\forall R>1,
\end{equation*}
where $C$ is independent of $R$. The proof is finished.
\qed

\section{\textbf{Classification}}\label{Section 5}
In this section, we classify the solution to equation \eqref{I-hL-eq}. That is, we finish the proof of Theorem \ref{hL-Classf-thm} by the Serrin-Zou type identity and the Pohozaev identity.

We first define the tensor as
\begin{equation*}
E_{ij}=\partial_{j}(|\nabla u|^{n-2}u_{i})-(|\nabla u|^{n-2}u_{i}u_{j}-\frac{1}{n}|\nabla u|^n\delta_{ij}).
\end{equation*}

\subsection{A differential identity} 

According to the limitation of the regularity of the solution, in the proof, we need to consider the regularization of the equation. The following differential identity holds which implies Theorem \ref{hL-Classf-thm}.
\begin{proposition}\label{diff-iden-prop}
Let $u$ be a solution to \eqref{I-hL-eq}. Then we have 
\begin{equation*}
E_{ij}=0,\qquad \qquad \text{in}~L^2_{\operatorname{loc}}(\overline{\mathbb R_+^n}).
\end{equation*}
In particular, we have $a(\nabla u) \in C^{1,\alpha}_{\operatorname{loc}}(\overline{\mathbb R_+^n})$.
\end{proposition}
\begin{proof}
Arguing as \cite{Z2024}, we make the following regularization argument.  

Set 
$$a^{\varepsilon}(x):=(|x|^2+\varepsilon^2)^{\frac{n-2}{2}}x,\qquad\text{for}~x\in\mathbb R^n,$$
$$H(x):=|x|^n~\text{and}~H^{\varepsilon}(x):=(|x|^2+\varepsilon^2)^{\frac{n}{2}},\qquad\text{for}~x\in\mathbb R^n.$$
A direct computation reveals that $a^{\varepsilon}\rightarrow a$ uniformly on compact subset of $\mathbb R^n$ as $\varepsilon\rightarrow 0$ and $a^{\varepsilon}$ satisfies \eqref{A-ak-inq} with $s_j$ replaced by $s_{\varepsilon}$. Arguing as the proof of Proposition \ref{2order-esti-prop}, let $u_{\varepsilon}\in W^{1,n}(B_{R_\varepsilon}^+)$ be a solution of 
\begin{equation}\label{C-regu-eq}
\begin{cases}
-\operatorname{div}(a^{\varepsilon}(\nabla u_{\varepsilon}))= 0,&\qquad\text{in }B_{R_\varepsilon}^+,\\
a^{\varepsilon}(\nabla u_{\varepsilon})\cdot\mathbf{n}=e^{u}, &\qquad\text{on }\Sigma_{R_\varepsilon},\\
u_\varepsilon=u,&\qquad\text{on }\partial B_{R_\varepsilon}^+,
\end{cases}
\end{equation}
where $R_\varepsilon\rightarrow+\infty$ as $\varepsilon\rightarrow0$. Also, the standard elliptic regularity theory \cite{GT1983, L1988, L1991} implies $u_\varepsilon\in C^{2,\alpha}(\overline{\mathbb R_+^n\cap B_{\frac{R_\varepsilon}{2}}})\cap C^{3,\alpha}_{\operatorname{loc}}(\mathbb R_+^n\cap B_{R_\varepsilon})$ and $u_{\varepsilon}\in C^{1,\alpha}_{\operatorname{loc}}(\overline{\mathbb R_+^n}\cap B_{R_\varepsilon})$, uniformly in $\varepsilon$. Hence, up to a subsequence, $u_{\varepsilon}\rightarrow u$ in $C_{\operatorname{loc}}^{1}(\overline{\mathbb R_+^n})$ as $\varepsilon\rightarrow 0$. The proof of Proposition \ref{2order-esti-prop} implies that $\{a^{\varepsilon}(\nabla u_{\varepsilon})\}$ is bounded in $W_{\operatorname{loc}}^{1,2}(\overline{\mathbb R_+^n}\cap B_{R_\varepsilon})$ uniformly in $\varepsilon$, and hence $a^{\varepsilon}(\nabla u_{\varepsilon})\rightharpoonup a(\nabla u)$ in $W_{\operatorname{loc}}^{1,2}(\overline{\mathbb R_+^n})$ up to a subsequence as $\varepsilon\rightarrow 0$. For $\varepsilon$ small enough, $R_\varepsilon>2R$.

Denote
$$X_{i}:=a_{i}(\nabla u),~X_{ij}:=\partial_{j}(X_{i}),~L_{ij}:=X_{i}u_{j}-\frac{1}{n}H(\nabla u)\delta_{ij},$$
$$X_{i}^{\varepsilon}:= a_{i}^{\varepsilon}(\nabla u_{\varepsilon}),~X_{ij}^{\varepsilon}:=\partial_{j}(X_{i}^{\varepsilon}),~
L_{ij}^{\varepsilon}:=X_{i}^{\varepsilon}\partial_{j}u_{\varepsilon}-\frac{1}{n}H^{\varepsilon}(\nabla u_{\varepsilon})\delta_{ij}.$$
Similarly,
\begin{equation*}
E_{ij}=X_{ij}-L_{ij},~\text{and }~E_{ij}^{\varepsilon}:=X_{ij}^{\varepsilon}-L_{ij}^{\varepsilon}.
\end{equation*}
Based on the above notations, the equation \eqref{C-regu-eq} can be rewritten as 
\begin{equation}\label{C-rregu-eq}
\begin{cases}
\sum\limits_{i} X_{ii}^{\varepsilon}=0,&\qquad\text{in}~\mathbb R^n_+,\\
X_n^{\varepsilon}=-e^{u},&\qquad\text{on}~\partial\mathbb R^n_+.
\end{cases}
\end{equation}
Note that we can write the matrix $E_{ij}^{\varepsilon}=X_{ij}^{\varepsilon}-L_{ij}^{\varepsilon}$ as $A^{\varepsilon}B^{\varepsilon}$, where
\begin{eqnarray*}
A^{\varepsilon}&=&\big(A_{ij}^{\varepsilon}\big)_{n\times n}=\big[(\partial_{j}a_{i}^{\varepsilon})(\nabla u_{\varepsilon})\big]_{n\times n}\nonumber\\
&=&(|\nabla u_{\varepsilon}|^2+\varepsilon^2)^{\frac{n-2}{2}}\big[I_n+(n-2)\frac{\nabla u_{\varepsilon}\otimes\nabla u_{\varepsilon}}{|\nabla u_{\varepsilon}|^2+\varepsilon^2}\big],\nonumber\\
B^{\varepsilon}&=&D^2u_{\varepsilon}-(A^{\varepsilon})^{-1}\big[(|\nabla u_{\varepsilon}|^2+\varepsilon^2)^{\frac{n-2}{2}}\nabla u_{\varepsilon}\otimes\nabla u_{\varepsilon}-\frac{1}{n}H^{\varepsilon}(\nabla u_{\varepsilon})I_n\big],
\end{eqnarray*}
where $I_n$ is the $n\times n$ identity matrix. Lemma \ref{matrix-lem} implies that
\begin{equation}\label{Eij-inq}
\sum\limits_{i,j}E_{ij}^{\varepsilon}E_{ji}^{\varepsilon}\geq C\sum\limits_{i,j}|E_{ij}^{\varepsilon}|^2.
\end{equation}
By \eqref{C-regu-eq}, we have the following Serrin-Zou type differential identity
\begin{eqnarray}\label{C-SZ-iden}
&&\sum\limits_{i}\partial_{i}\big(e^{-u}\sum\limits_{j}E_{ij}^{\varepsilon}X_{j}^{\varepsilon}\big)\nonumber\\
&=&\sum\limits_{i}\partial_{i}\big[e^{-u}\big(\sum\limits_{j}X_{ij}^{\varepsilon}X_{j}^{\varepsilon}-X_{i}^{\varepsilon}X_{j}^{\varepsilon}\partial_{j}u_{\varepsilon}+\frac{1}{n}H^{\varepsilon}(\nabla u_{\varepsilon})X_{i}^{\varepsilon}\big)\big]\nonumber\\
&=&e^{-u}\big[\sum\limits_{i,j}\big(X_{ij}^{\varepsilon}X_{ji}^{\varepsilon}-X_{ij}^{\varepsilon}X_{j}^{\varepsilon}u_{i}-X_{i}^{\varepsilon}X_{ji}^{\varepsilon}\partial_{j}u_{\varepsilon}-X_{i}^{\varepsilon}X_{j}^{\varepsilon}\partial_{ij}u_{\varepsilon}\nonumber\\
&&+X_{i}^{\varepsilon}X_{j}^{\varepsilon}\partial_{j}u_{\varepsilon}u_{i}\big)+\frac{1}{n}\sum\limits_{i}\partial_{i}(H^{\varepsilon}(\nabla u_{\varepsilon}))X_{i}^{\varepsilon}-\frac{1}{n}\sum\limits_{i}H^{\varepsilon}(\nabla u_{\varepsilon})X_{i}^{\varepsilon}u_{i}\big]\nonumber\\
&=&e^{-u}\big[\sum\limits_{i,j}\big(X_{ij}^{\varepsilon}X_{ji}^{\varepsilon}-X_{ij}^{\varepsilon}X_{j}^{\varepsilon}u_{i}-X_{ij}^{\varepsilon}X_{j}^{\varepsilon}\partial_{i}u_{\varepsilon}
+X_{i}^{\varepsilon}X_{j}^{\varepsilon}\partial_{j}u_{\varepsilon}u_{i}\big)\nonumber\\
&&-\frac{1}{n}\sum\limits_{i}H^{\varepsilon}(\nabla u_{\varepsilon})X_{i}^{\varepsilon}u_{i}\big]\nonumber\\
&=&e^{-u}\big[\sum\limits_{i,j}\big(X_{ij}^{\varepsilon}X_{ji}^{\varepsilon}-2X_{ij}^{\varepsilon}X_{j}^{\varepsilon}\partial_{i}u_{\varepsilon}\big)+\sum\limits_{i,j}X_{ij}^{\varepsilon}X_{j}^{\varepsilon}(\partial_{i}u_{\varepsilon}-u_{i})\nonumber\\
&&+\sum\limits_{i,j}X_{i}^{\varepsilon}X_{j}^{\varepsilon}\partial_{j}u_{\varepsilon}u_{i}-\frac{1}{n}\sum\limits_{i}H^{\varepsilon}(\nabla u_{\varepsilon})X_{i}^{\varepsilon}u_{i}\big]\nonumber\\
&=&e^{-u}\sum\limits_{i,j}(X_{ij}^{\varepsilon}-L_{ij}^{\varepsilon})(X_{ji}^{\varepsilon}-L_{ji}^{\varepsilon})+P^{\varepsilon}\nonumber\\
&=&e^{-u}\sum\limits_{i,j}E_{ij}^{\varepsilon}E_{ji}^{\varepsilon}+P^{\varepsilon},
\end{eqnarray}
where
\begin{eqnarray*}
P^{\varepsilon}&=&e^{-u}\big[\sum\limits_{i,j}X_{ij}^{\varepsilon}X_{j}^{\varepsilon}(\partial_{i}u_{\varepsilon}-u_{i})+\sum\limits_{i,j}X_{i}^{\varepsilon}X_{j}^{\varepsilon}\partial_{j}u_{\varepsilon}u_{i}-\frac{1}{n}\sum\limits_{i}H^{\varepsilon}(\nabla u_{\varepsilon})X_{i}^{\varepsilon}u_{i}\nonumber\\
&&-\sum\limits_{i,j}L_{ij}^{\varepsilon}L_{ji}^{\varepsilon}\big]\nonumber\\
&=&e^{-u}\big[\sum\limits_{i,j}X_{ij}^{\varepsilon}X_{j}^{\varepsilon}(\partial_{i}u_{\varepsilon}-u_{i})-\sum\limits_{i,j}X_{i}^{\varepsilon}X_{j}^{\varepsilon}\partial_{j}u_{\varepsilon}(\partial_{i}u_{\varepsilon}-u_{i})\nonumber\\
&&+\frac{1}{n}\sum\limits_{i}H^{\varepsilon}(\nabla u_{\varepsilon})X_{i}^{\varepsilon}(\partial_{i}u_{\varepsilon}-u_{i})+\frac1 n\sum\limits_{i}H^{\varepsilon}(\nabla u_{\varepsilon})\big(X_{i}^{\varepsilon}\partial_{i}u_{\varepsilon}-H^{\varepsilon}(\nabla u_{\varepsilon})\big)\big].
\end{eqnarray*}
It follows that
\begin{eqnarray*}
|P^{\varepsilon}|
&\leq& Ce^{-u}\big[|\nabla X^{\varepsilon}||X^{\varepsilon}||\nabla u_{\varepsilon}-\nabla u|+|X^{\varepsilon}|^2|\nabla u_{\varepsilon}||\nabla u_{\varepsilon}-\nabla u|\\
&&+(|\nabla u_{\varepsilon}|^2+\varepsilon^2)^{\frac{n}{2}}|X^{\varepsilon}||\nabla u_{\varepsilon}-\nabla u|+(|\nabla u_{\varepsilon}|^2+\varepsilon^2)^{n-1}\varepsilon^2\big]\\
&\leq&Ce^{-u}(|\nabla u_{\varepsilon}|^2+\varepsilon^2)^{\frac{n-2}{2}}\big[|\nabla X^{\varepsilon}||\nabla u_{\varepsilon}||\nabla u_{\varepsilon}-\nabla u|\\
&&+|X^{\varepsilon}||\nabla u_{\varepsilon}|^2|\nabla u_{\varepsilon}-\nabla u|+(|\nabla u_{\varepsilon}|^2+\varepsilon^2)|X^{\varepsilon}||\nabla u_{\varepsilon}-\nabla u|\\
&&+(|\nabla u_{\varepsilon}|^2+\varepsilon^2)^{\frac{n}{2}}\varepsilon^2\big].
\end{eqnarray*}
Since $\{u_{\varepsilon}\}$ is bounded in $C^{1}_{\operatorname{loc}}(\overline{\mathbb R^n_+})$, $u_{\varepsilon}\rightarrow u$ in $C^{1}_{\operatorname{loc}}(\overline{\mathbb R^n_+})$ and $X^{\varepsilon}\in W^{1,2}_{\operatorname{loc}}(\overline{\mathbb R^n_+})$ uniformly in $\varepsilon$, we have $P^{\varepsilon}\rightarrow 0$ in $L^2_{\operatorname{loc}}(\overline{\mathbb R_+^n})$ as $\varepsilon\rightarrow 0$.

Next, again by \eqref{C-regu-eq}, we get
\begin{equation}\label{C-Pohozaev-iden}
\sum\limits_{i}\partial_n(X_{i}^{\varepsilon}\partial_{i}u_{\varepsilon})-n\sum\limits_{i}\partial_{i}(X_{i}^{\varepsilon}\partial_nu_{\varepsilon})=\sum\limits_{i}X_{in}^{\varepsilon}\partial_{i}u_{\varepsilon}-(n-1)\sum\limits_{i}X_{i}^{\varepsilon}\partial_{ni}u_{\varepsilon}.
\end{equation}
Using the definition of $X_{ij}^{\varepsilon}$, i.e., 
$$X_{ij}^{\varepsilon}=(|\nabla u_{\varepsilon}|^2+\varepsilon^2)^{\frac{n-2}{2}}\big(\partial_{ij}u_{\varepsilon}+(n-2)\sum\limits_{k}\frac{\partial_{i}u_{\varepsilon}\partial_{k}u_{\varepsilon}\partial_{kj}u_{\varepsilon}}{|\nabla u_{\varepsilon}|^2+\varepsilon^2}\big),$$ 
there holds
\begin{eqnarray}\label{C-Xni-eq}
\sum\limits_{i}X_{in}^{\varepsilon}\partial_{i}u_{\varepsilon}&=&(|\nabla u_{\varepsilon}|^2+\varepsilon^2)^{\frac{n-2}{2}}(1+\frac{(n-2)|\nabla u_{\varepsilon}|^2}{|\nabla u_{\varepsilon}|^2+\varepsilon^2})\sum\limits_{i}\partial_{i}u_{\varepsilon}\partial_{ni}u_{\varepsilon}\nonumber\\
&=&(1+\frac{(n-2)|\nabla u_{\varepsilon}|^2}{|\nabla u_{\varepsilon}|^2+\varepsilon^2})\sum\limits_{i}X_{i}^{\varepsilon}\partial_{i}\partial_nu_{\varepsilon}.
\end{eqnarray}
Combining \eqref{C-Pohozaev-iden} and \eqref{C-Xni-eq}, we obtain the following Pohozaev-type differential identity
\begin{eqnarray}\label{C-PXni-eq}
&&\sum\limits_{i}\partial_n(X_{i}^{\varepsilon}\partial_{i}u_{\varepsilon})-n\sum\limits_{i}\partial_{i}(X_{i}^{\varepsilon}\partial_nu_{\varepsilon})=-\frac{(n-2)\varepsilon^2}{|\nabla u_{\varepsilon}|^2+\varepsilon^2}\sum\limits_{i}X_{i}^{\varepsilon}\partial_{ni}u_{\varepsilon}\nonumber\\
&&=-(1+\frac{(n-2)|\nabla u_{\varepsilon}|^2}{|\nabla u_{\varepsilon}|^2+\varepsilon^2})^{-1}\frac{(n-2)\varepsilon^2}{|\nabla u_{\varepsilon}|^2+\varepsilon^2}\sum\limits_{i}X_{in}^{\varepsilon}\partial_{i}u_{\varepsilon}=:Q^{\varepsilon}.
\end{eqnarray}
By the Cauchy-Schwarz inequality $ab\leq\frac{a^2+b^2}{2}$, it holds
\begin{eqnarray*}
|Q^{\varepsilon}|&\leq& (n-2)\frac{\varepsilon^2|\nabla u_{\varepsilon}|}{|\nabla u_{\varepsilon}|^2+\varepsilon^2}|\nabla X^{\varepsilon}|\nonumber\\
&\leq& (n-2)\frac{\varepsilon}{|\nabla u_{\varepsilon}|^2+\varepsilon^2}|\nabla X^{\varepsilon}|\frac{\varepsilon^2+|\nabla u_\varepsilon|^2}{2}=\frac{n-2}{2}\varepsilon|\nabla X^{\varepsilon}|.
\end{eqnarray*}
Recalling the fact $X^{\varepsilon}\in W^{1,2}_{\operatorname{loc}}(\overline{\mathbb R^n_+})$, uniformly in $\varepsilon$, it implies $Q^{\varepsilon}\rightarrow 0$ in $L^2_{\operatorname{loc}}(\overline{\mathbb R^n_+})$.

Define the nonnegative cut-off function $\eta\in C_c^{\infty}(B_{2R})$ as
\begin{equation*}
\eta=
\begin{cases}
 1,\quad  &\text{in}~B_{R},\\
 0, \quad&\text{in}~\mathbb R^n\setminus B_{2R},   
\end{cases}
\end{equation*} with $0\leq\eta\leq1$ and $|\nabla \eta|\leq\frac{2}{R}$ in $\mathbb R^n$.

Multiply \eqref{C-SZ-iden} by $\eta^2$ and integrate over $\mathbb R^n_+$. Using \eqref{Eij-inq} and integration by parts, one has
\begin{eqnarray}\label{C-Eij-int}
& &\sum\limits_{i,j}\int_{\mathbb R^n_+}e^{-u}|E_{ij}^{\varepsilon}|^2\eta^2\mathrm dx\nonumber\\
&\leq &C\sum\limits_{i,j}\int_{\mathbb R^n_+}e^{-u}E_{ij}^{\varepsilon}E_{ji}^{\varepsilon}\eta^2\mathrm dx=C\int_{\mathbb R^n_+}\big[\sum\limits_{i,j}\partial_{i}\big(e^{-u}E_{ij}^{\varepsilon}X_{j}^{\varepsilon}\big)-P^{\varepsilon}\big]\eta^2\mathrm dx\nonumber\\
&\leq& C\big(-2\sum\limits_{i,j}\int_{\mathbb R^n_+}e^{-u}E_{ij}^{\varepsilon}X_{j}^{\varepsilon}\eta\eta_{i}\mathrm dx-\sum\limits_{j}\int_{\partial\mathbb R^n_+}e^{-u}E_{nj}^{\varepsilon}X_{j}^{\varepsilon}\eta^2\mathrm dx'-\int_{\mathbb R^n_+}P^{\varepsilon}\eta^2\mathrm dx\big),\nonumber\\
&&
\end{eqnarray}
where $C$ is a positive constant that depends only on $n$. For the first term on the RHS of \eqref{C-Eij-int}, the H\"{o}lder inequality entails that
\begin{eqnarray}\label{C-first-term}
&&-2\sum\limits_{i,j}\int_{\mathbb R^n_+}e^{-u}E_{ij}^{\varepsilon}X_{j}^{\varepsilon}\eta\eta_{i}\mathrm dx\nonumber\\
&\leq&2\sum\limits_{i,j}(\int_{B_{2R}^+\setminus B_R^+}e^{-u}|E_{ij}^{\varepsilon}|^2\eta^2\mathrm dx)^{\frac{1}{2}}(\int_{B_{2R}^+\setminus B_R^+}e^{-u}|X_{j}^{\varepsilon}|^2|\eta_{i}|^2\mathrm dx)^{\frac{1}{2}}.
\end{eqnarray}
Focusing on the second term on the RHS in \eqref{C-Eij-int}, by \eqref{C-rregu-eq} and integration by parts, we obtain
\begin{eqnarray}\label{C-boundary-iden}
&&\int_{\partial\mathbb R^n_+}e^{-u}\sum\limits_{j}E_{nj}^{\varepsilon}X_{j}^{\varepsilon}\eta^2\mathrm dx'\nonumber\\
&=&\int_{\partial\mathbb R^n_+}e^{-u}\big(\sum\limits_{j}X_{nj}^{\varepsilon}X_{j}^{\varepsilon}-\sum\limits_{j}X_n^{\varepsilon}X_{j}^{\varepsilon}\partial_{j}u_{\varepsilon}+\frac{1}{n}H^{\varepsilon}(\nabla u_{\varepsilon})X_n^{\varepsilon}\big)\eta^2\mathrm dx'\nonumber\\
&=&\int_{\partial\mathbb R^n_+}\big(-\sum\limits_{j=1}^{n-1}u_{j}X_{j}^{\varepsilon}+\sum\limits_{j=1}^{n-1}X_{jj}^{\varepsilon}+\sum\limits_{j}X_{j}^{\varepsilon}\partial_{j}u_{\varepsilon}-\frac{1}{n}H^{\varepsilon}(\nabla u_{\varepsilon})\big)\eta^2\mathrm dx'\nonumber\\
&=&\int_{\partial\mathbb R^n_+}\big[\big(-\sum\limits_{j=1}^{n-1}u_{j}X_{j}^{\varepsilon}++\sum\limits_{j}X_{j}^{\varepsilon}\partial_{j}u_{\varepsilon}-\frac{1}{n}H^{\varepsilon}(\nabla u_{\varepsilon})\big)\eta^2-2\sum\limits_{j=1}^{n-1}X_{j}^{\varepsilon}\eta\eta_{j}\big]\mathrm dx'\nonumber\\
&=&\int_{\partial\mathbb R^n_+}\big(-\sum\limits_{j=1}^{n-1}X_{j}^{\varepsilon}\partial_{j}u_{\varepsilon}+\frac{n-1}{n}\sum\limits_{j=1}^nX_{j}^{\varepsilon}\partial_{j}u_{\varepsilon}+F^{\varepsilon}\big)\eta^2\mathrm dx'\nonumber\\
&&-2\int_{\partial\mathbb R^n_+}\sum\limits_{j=1}^{n-1}X_{j}^{\varepsilon}\eta\eta_{j}\mathrm dx'\nonumber\\
&=&-\frac{1}{n}\int_{\partial\mathbb R^n_+}\big(\sum\limits_{j}X_{j}^{\varepsilon}\partial_{j}u_{\varepsilon}-nX_n^{\varepsilon}\partial_nu_{\varepsilon}\big)\eta^2\mathrm dx'+\int_{\partial\mathbb R^n_+}F^{\varepsilon}\eta^2\mathrm dx'\nonumber\\
&&-2\int_{\partial\mathbb R^n_+}\sum\limits_{j=1}^{n-1}X_{j}^{\varepsilon}\eta\eta_{j}\mathrm dx',
\end{eqnarray}
where
\begin{eqnarray*}
F^{\varepsilon}&=&\sum\limits_{j=1}^{n-1}X_{j}^{\varepsilon}(\partial_{j}u_{\varepsilon}-u_{j})+\sum\limits_{j}X_{j}^{\varepsilon}\partial_{j}u_{\varepsilon}-\frac{1}{n}H^{\varepsilon}(\nabla u_{\varepsilon})-\frac{n-1}{n}\sum\limits_{j}X_{j}^{\varepsilon}\partial_{j}u_{\varepsilon}\nonumber\\
&=&\sum\limits_{j=1}^{n-1}X_{j}^{\varepsilon}(\partial_{j}u_{\varepsilon}-u_{j})-\frac{1}{n}\big[H^{\varepsilon}(\nabla u_{\varepsilon})-\sum\limits_{j}X_{j}^{\varepsilon}\partial_{j}u_{\varepsilon}\big]\nonumber\\
&=&\sum\limits_{j=1}^{n-1}X_{j}^{\varepsilon}(\partial_{j}u_{\varepsilon}-u_{j})-\frac{1}{n}\big(|\nabla u_{\varepsilon}|^2+\varepsilon^2\big)^{\frac{n-2}{2}}\varepsilon^2.
\end{eqnarray*}

The fact that $\{u_{\varepsilon}\}$ is bounded in $C^{1}_{\operatorname{loc}}(\overline{\mathbb R^n_+})$ implies $\big(|\nabla u_{\varepsilon}|^2+\varepsilon^2\big)^{\frac{n-2}{2}}\varepsilon^2\rightarrow 0$ uniformly on compact subset of $\overline{\mathbb R^n_+}$. Furthermore, since $u_{\varepsilon}\rightarrow u$ in $C^{1}_{\operatorname{loc}}(\overline{\mathbb R^n_+})$, we get $F^{\varepsilon}\rightarrow 0$ uniformly in the compact subset of $\overline{\mathbb R^n_+}$.

Note that the integral form of \eqref{C-PXni-eq} is as follows
\begin{eqnarray}\label{C-Pohozaev-int}
\int_{\mathbb R^n_+}Q^{\varepsilon}\eta^2\mathrm dx&=&\int_{\mathbb R^n_+}\big[\sum\limits_{i}\partial_n(X_{i}^{\varepsilon}\partial_{i}u_{\varepsilon})-n\sum\limits_{i}\partial_{i}(X_{i}^{\varepsilon}\partial_nu_{\varepsilon})\big]\eta^2\mathrm dx\nonumber\\
&=&-2\int_{\mathbb R^n_+}\big(\sum\limits_{i}X_{i}^{\varepsilon}\partial_{i}u_{\varepsilon}\eta\eta_n-n\sum\limits_{i} X_{i}^{\varepsilon}\partial_nu_{\varepsilon}\eta\eta_{i}\big)\mathrm dx\nonumber\\
&&-\int_{\partial\mathbb R^n_+}\big(\sum\limits_{i}X_{i}^{\varepsilon}\partial_{i}u_{\varepsilon}-nX_n^{\varepsilon}\partial_nu_{\varepsilon}\big)\eta^2\mathrm dx'.
\end{eqnarray}
In view of \eqref{C-boundary-iden} and \eqref{C-Pohozaev-int}, we obtain
\begin{eqnarray}\label{C-bP-iden}
&&\int_{\partial\mathbb R^n_+}e^{-u}\sum\limits_{j}E_{nj}^{\varepsilon}X_{j}^{\varepsilon}\eta^2\mathrm dx'\nonumber\\
&=&\frac{1}{n}\int_{\mathbb R^n_+}Q^{\varepsilon}\eta^2\mathrm dx+\int_{\partial\mathbb R^n_+}F^{\varepsilon}\eta^2\mathrm dx'-2\int_{\partial\mathbb R^n_+}\sum\limits_{j=1}^{n-1}X_{j}^{\varepsilon}\eta\eta_{j}\mathrm dx'\nonumber\\
&&+\frac{2}{n}\int_{\mathbb R^n_+}\big(\sum\limits_{i}X_{i}^{\varepsilon}\partial_{i}u_{\varepsilon}\eta\eta_n-n\sum\limits_{i} X_{i}^{\varepsilon}\partial_nu_{\varepsilon}\eta\eta_{i}\big)\mathrm dx.
\end{eqnarray}
Furthermore, it follows from \eqref{C-Eij-int}, \eqref{C-first-term} and \eqref{C-bP-iden} that
\begin{eqnarray}\label{C-Eij-inq}
&&\sum\limits_{i,j}\int_{\mathbb R^n_+}e^{-u}|E_{ij}^{\varepsilon}|^2\eta^2\mathrm dx\nonumber\\
&\leq& C\big[\sum\limits_{i,j}(\int_{B_{2R}^+\setminus B_R^+}e^{-u}|E_{ij}^{\varepsilon}|^2\eta^2\mathrm dx)^{\frac{1}{2}}(\int_{B_{2R}^+\setminus B_R^+}e^{-u}|X_{j}^{\varepsilon}|^2|\eta_{i}|^2\mathrm dx)^{\frac{1}{2}}\nonumber\\
&&+\int_{\mathbb R^n_+}(|P^{\varepsilon}|+|Q^{\varepsilon}|)\eta^2\mathrm dx+\int_{\partial\mathbb R^n_+}|F^{\varepsilon}|\eta^2\mathrm dx'\nonumber\\
&&+\int_{\partial\mathbb R^n_+}|X^{\varepsilon}||\nabla\eta|\eta\mathrm dx'+\int_{\mathbb R^n_+}|X^{\varepsilon}||\nabla u_{\varepsilon}||\nabla\eta|\eta\mathrm dx\big].
\end{eqnarray}
Letting $\varepsilon\rightarrow 0$, we get that
\begin{eqnarray}\label{C-limt-inq}
&&\limsup\limits_{\varepsilon\rightarrow 0}\sum\limits_{i,j}\int_{\mathbb R^n_+}e^{-u}|E_{ij}^{\varepsilon}|^2\eta^2\mathrm dx\nonumber\\
&\leq& C\big[\sum\limits_{i,j}(\int_{B_{2R}^+\setminus B_R^+}e^{-u}|E_{ij}|^2\eta^2\mathrm dx)^{\frac{1}{2}}(\int_{B_{2R}^+\setminus B_R^+}e^{-u}|X|^2|\nabla \eta|^2\mathrm dx)^{\frac{1}{2}}\nonumber\\
&&+\int_{\partial\mathbb R^n_+}|X|\eta|\nabla\eta|\mathrm dx'+\int_{\mathbb R^n_+}|X||\nabla u||\nabla\eta|\eta\mathrm dx\big].
\end{eqnarray}
Recalling Proposition \ref{1order-esti-prop} and \ref{2order-esti-prop}, we deduce
\begin{equation*}
\int_{B_{2R}^+\setminus B_{R}^+}e^{-u}|E_{ij}|^2\eta^2\mathrm dx\leq C\int_{B_{2R}^+\setminus B_{R}^+}e^{-u}[|\nabla( a(\nabla u))|^2+|\nabla u|^{2n}]\mathrm dx\leq CR^{0}=C,
\end{equation*}
\begin{equation}\label{C-int-esti2}
\int_{B_{2R}^+\setminus B_{R}^+}e^{-u}|X|^2|\nabla \eta|^2\mathrm dx\leq C\int_{B_{2R}^+\setminus B_{R}^+}e^{-u}|\nabla u|^{2n-2}|\nabla \eta|^2\mathrm dx\leq CR^{0}=C,
\end{equation}
\begin{equation}\label{C-int-esti3}
\int_{\partial\mathbb R^n_+}|X||\nabla\eta|\eta\mathrm dx'\leq C\int_{\partial\mathbb R^n_+}|\nabla u|^{n-1}|\nabla\eta|\mathrm dx'\leq\frac{C}{R},
\end{equation}
\begin{equation}\label{C-int-esti4}
\int_{\mathbb R^n_+}|X||\nabla u||\nabla\eta|\eta\mathrm dx\leq C\int_{\mathbb R^n_+}|\nabla u|^n|\nabla\eta|\mathrm dx\leq\frac{C}{R}.
\end{equation}
Combining the above with \eqref{C-limt-inq} yields
\begin{equation}\label{C-lim-inq}
\limsup\limits_{\varepsilon\rightarrow 0}\sum\limits_{i,j}\int_{\mathbb R^n_+}e^{-u}|E_{ij}^{\varepsilon}|^2\eta^2\mathrm dx\leq C+\frac{C}{R}.
\end{equation}
Since $X_{ij}^{\varepsilon}\rightharpoonup X_{ij}$ in $L^2_{\operatorname{loc}}(\overline{\mathbb R^n_+})$, $X_{i}^{\varepsilon}\rightharpoonup X_{i}$ in $L^2_{\operatorname{loc}}(\overline{\mathbb R^n_+})$, $\partial_{j}u_{\varepsilon}\rightarrow u_{j}$ in $C^{0}_{\operatorname{loc}}(\overline{\mathbb R^n_+})$ and $H^{\varepsilon}(\nabla u_{\varepsilon})\rightharpoonup H(\nabla u)$ in $L^2_{\operatorname{loc}}(\overline{\mathbb R^n_+})$, we have 
\begin{equation*}
  E_{ij}^{\varepsilon}\rightharpoonup E_{ij},\qquad\text{in}~L^2_{\operatorname{loc}}(\overline{\mathbb R^n_+}),
\end{equation*}
as $\varepsilon\rightarrow 0$. Furthermore, $e^{-\frac{u}{2}}\eta\in L^{\infty}(\overline{\mathbb R^n_+})$ implies
\begin{equation}\label{C-Eij-limit}
e^{-\frac{u}{2}}E_{ij}^{\varepsilon}\eta\rightharpoonup e^{-\frac{u}{2}}E_{ij}\eta,\qquad\text{in}~L^2(\overline{\mathbb R^n_+}),
\end{equation}
as $\varepsilon\rightarrow 0$. Therefore, by \eqref{C-lim-inq} and \eqref{C-Eij-limit}, one has
\begin{equation*}
\sum\limits_{i,j}\int_{\mathbb R^n_+}e^{-u}|E_{ij}|^2\eta^2\mathrm dx\leq C,\qquad\text{as}~R\rightarrow\infty.
\end{equation*}
Hence, $e^{-\frac{u}{2}}E_{ij}\in L^2(\mathbb R^n_+)$. This implies
\begin{equation}\label{C-lim0-eq}
\sum\limits_{i,j}\int_{B_{2R}^{+}\setminus B_{R}^{+}}e^{-u}|E_{ij}|^2\eta^2\mathrm dx\rightarrow 0,\qquad\text{as}~R\rightarrow\infty.
\end{equation}
By the inequality \eqref{C-limt-inq} with $\varepsilon\rightarrow 0$, using \eqref{C-int-esti2}-\eqref{C-int-esti4} and \eqref{C-lim0-eq} with $R\rightarrow\infty$, we obtain
\begin{equation*}
\sum\limits_{i,j}\int_{\mathbb R^n_+}e^{-u}|E_{ij}|^2\eta^2\mathrm dx\rightarrow 0, \qquad\text{as}~R\rightarrow\infty.
\end{equation*}
Hence $E_{ij}=0$ in $L^2_{\operatorname{loc}}(\overline{\mathbb R^n_+})$. It follows from $u\in C^{1,\alpha}_{\operatorname{loc}}(\overline{\mathbb R^n_+})$ and $E_{ij}=0$ that
\begin{equation*}
\partial_{j}(|\nabla u|^{n-2}u_{i})=|\nabla u|^{n-2}u_{i}u_{j}-\frac{1}{n}|\nabla u|^n\delta_{ij}\in C^{\alpha}_{\operatorname{loc}}(\overline{\mathbb R^n_+}),\qquad\forall~i,j\in\{1,2,\cdots,n\},
\end{equation*}
which implies $a(\nabla u)=|\nabla u|^{n-2}\nabla u\in C^{1,\alpha}_{\operatorname{loc}}(\overline{\mathbb R^n_+})$.
This finishes the proof.
\end{proof}

\subsection{Proof of Theorem \ref{hL-Classf-thm}}
By Proposition \ref{diff-iden-prop}, we have $E_{ij}=0$ in $L^2_{\operatorname{loc}}(\overline{\mathbb R_+^n})$, namely
\begin{equation}\label{C-diff-iden}
\partial_{j}(|\nabla u|^{n-2}u_{i})=|\nabla u|^{n-2}u_{i}u_{j}-\frac{1}{n}|\nabla u|^n\delta_{ij},\qquad \text{in}~L^2_{\operatorname{loc}}(\overline{\mathbb R_+^n}),
\end{equation}
and $a(\nabla u)=|\nabla u|^{n-2}\nabla u\in C^{1,\alpha}_{\operatorname{loc}}(\overline{\mathbb R_+^n})$.
We consider the auxiliary function
\begin{equation*}
v=e^{-\frac{1}{n-1}u},
\end{equation*} 
where $u$ is a solution of \eqref{I-hL-eq}. An immediate computation shows that $v>0$ satisfies the following equation
\begin{equation}\label{C-v-eq}
\begin{cases}
\Delta_n v= \frac{(n-1)|\nabla v|^n}{v},&\qquad \text{in}~\mathbb R_+^n,\\
|\nabla v|^{n-2}\frac{\partial v}{\partial t}=(n-1)^{1-n}, &\qquad \text{on}~\partial\mathbb R_+^n.
\end{cases}
\end{equation}
Naturally, we have $v\in C^{1,\alpha}_{\operatorname{loc}}(\overline{\mathbb R^n_+})$ and $a(\nabla v)=-\frac{1}{(n-1)^{n-1}}e^{-u}a(\nabla u)\in C^{1,\alpha}_{\operatorname{loc}}(\overline{\mathbb R^n_+})$.
By \eqref{C-diff-iden} we obtain
\begin{equation}\label{C-v-iden}
\partial_{j}(|\nabla v|^{n-2}v_{i})=\iota(x)\delta_{ij},\qquad \text{in}~L^2_{\operatorname{loc}}(\overline{\mathbb R_+^n}),
\end{equation}
where $\iota(x)=\frac{n-1}{n}\frac{|\nabla v|^n}{v}$. By the elliptic regularity theory, we have $v\in C^{2,\alpha}_{\operatorname{loc}}(\mathbb R_+^n\cap\{\nabla v\neq 0\})$, which implies $\iota(x)\in C^{1,\alpha}_{\operatorname{loc}}(\mathbb R_+^n\cap\{\nabla v\neq 0\})$. Using \eqref{C-v-iden} yields $a(\nabla v)\in C^{2,\alpha}_{\operatorname{loc}}(\mathbb R_+^n\cap\{\nabla v\neq 0\})$. Hence, fix $i\in\{1,2,\cdots,n\}$, for $j\neq i$, from \eqref{C-v-iden}, we obtain
\begin{equation*}
\partial_{j}\iota(x)=\partial_{j}\partial_{i}[a_{i}(\nabla v)]=\partial_{i}\partial_{j}[a_{i}(\nabla v)]=0
\end{equation*}
for any $x\in\mathbb R_+^n\cap\{\nabla v\neq 0\}$, namely, $\iota(x)$ is constant on each connected component of $\mathbb R_+^n\cap\{\nabla v\neq 0\}$. The asymptotic limits of $\nabla u$ in \eqref{A-gradu-esti} implies $\overline{\mathbb R_+^n}\setminus B_{R}\subset\{\nabla v\neq 0\}$ for large $R$. Consider the connected component $U$ of $\mathbb R_+^n\cap\{\nabla v\neq 0\}$ such that $\overline{\mathbb R_+^n}\setminus B_{R}\subset U$. Notice that $\iota(x)\in  C^{0,\alpha}_{\operatorname{loc}}(\overline{\mathbb R^n_+})$ by $v\in  C^{1,\alpha}_{\operatorname{loc}}(\overline{\mathbb R^n_+})$. Therefore, $a(\nabla v)=C(x-x_0)$ in $\overline{U}$ for some $C\neq0$, which implies
$$v=C_{1}+C_{2}|x-x_0|^{\frac{n}{n-1}}$$
in $\overline{U}$, for some $C_{1}\geq0,C_{2}>0,~x_0=(x_0',x_0^n)\in\mathbb R^n=\mathbb R^{n-1}\times\mathbb R$.
Using the fact that $a_n(\nabla v)=|\nabla v|^{n-1} \frac{\partial v}{\partial t}=(n-1)^{1-n}$, we have $x_0\in\mathbb R^n_{-}$ and 
\begin{equation*}
C_{2}=\frac{1}{n(-x_0^n)^{\frac{1}{n-1}}}.
\end{equation*}
Since for any $y\in\partial U$, either $y\in\partial\mathbb R^n_+$ or $a(\nabla v(y))=0$, we deduce that $\partial U\subset\partial\mathbb R_+^n$, implying $\overline{U}=\overline{\mathbb R_+^n}$. Using (\ref{C-v-eq}), we get $C_1=0$ and
\begin{equation*}
v(x)=\frac{|x-x_0|^{\frac{n}{n-1}}}{n\lambda^{\frac{1}{n-1}}},\qquad\text{in}~\overline{\mathbb R_+^n}
\end{equation*}
where $x_0=(x_0',-\lambda)\in\mathbb R^n_{-}$.
Combining with $v=e^{-\frac{1}{n-1}u}$, we obtain
\begin{equation*}
u(x)=\log v^{1-n}=\log \frac{n^{n-1}\lambda}{\big(|x'-x_0'|^2+(t+\lambda)^2\big)^{\frac{n}{2}}},\qquad\text{in}~\overline{\mathbb R_+^n}
\end{equation*}
where $x_0'\in\mathbb R^{n-1}$ and $\lambda>0$. This completes the proof of Theorem \ref{hL-Classf-thm}. 
\qed

\subsection{Complete the proof of Theorem \ref{Ot-inq-thm}} Define a functional as
\begin{equation*}
I(w)=\frac{e^{F(w)}}{G(w)},\quad w\in W_{\mu_n}(\mathbb R_+^n),
\end{equation*}
where
\begin{gather*}
F(w)= \frac{2}{n^n\sigma_{n-1}}\int_{\mathbb R_+^n}K_n(x,\nabla w)\mathrm dx+\int_{\partial\mathbb R_+^n}w\mathrm d\mu_n(x')\\ \intertext{and}
G(w)=\int_{\partial\mathbb R_+^n}e^{w}\mathrm d\mu_n(x').
\end{gather*}
If $w$ is a minimizer of $I(w)$, we have
\begin{equation*}
0=\frac{\mathrm dI(w+tv)}{\mathrm dt}\Big|_{t=0}=\frac{e^{F(w)}}{G(w)}\big[F'(w)\cdot v-\frac{1}{G(w)}G'(w)\cdot v\big],\quad\forall v\in W_{\mu_n}(\mathbb R_+^n)
\end{equation*}
that is,
\begin{eqnarray*}
F'(w)\cdot v=\frac{1}{G(w)}G'(w)\cdot v.    
\end{eqnarray*}
Thus, under the restriction $\frac{1}{G(w)}=L>0$, 
the Euler-Lagrange equation of inequality \eqref{I-Ot-inq}  is 
\begin{equation}\label{T-Ot-EL}
\begin{cases}
-\operatorname{div}\big[|X+\nabla w|^{n-2}(X+\nabla w)\big]=0, &\qquad\text{in }\mathbb R_+^n,\\
|X+\nabla w|^{n-2}(X+\nabla w)\cdot\mathbf{n}=\frac{n^{n-1}\sigma_{n-1}}{2}L\mathrm{e}^{w}\mu_n,&\qquad\text{on }\partial\mathbb R_+^n,
\end{cases}
\end{equation}
where $X=-\frac{n(x', 1+t)}{(1+t)^2 + |x'|^2}=\nabla[\log\mu_n(x)]$. Take the transformation $u(x)=\log\mu_n(x)$ $+w(x)+\log (\frac{n^{n-1}\sigma_{n-1}}{2}L)$, and then \eqref{T-Ot-EL} can be reduced to \eqref{I-hL-eq}.

Since $w\in W_{\mu_n}(\mathbb{R}_+^n)$, we have $\int_{\partial\mathbb{R}_+^n} e^w \mathrm d\mu_n<+\infty$ clearly and $\int_{\mathbb{R}_+^n} e^{\frac{n}{n-1}w} \mu_n^{\frac{n}{n-1}}\mathrm dx<+\infty$ by Proposition \ref{Pro energy In}. These show that $u$ satisfies the assumption \eqref{I-fm-assu}. Since $\log\mu_n\in W_{\operatorname{loc}}^{1,n}(\overline{\mathbb{R}_+^n})$, we know from Proposition \ref{contain-lem} that $u\in W_{\operatorname{loc}}^{1,n}(\overline{\mathbb{R}_+^n})$. That is, $u$ is a weak solution of \eqref{I-hL-eq} and satisfies all assumptions of Theorem \ref{hL-Classf-thm}.

So, $u(x)$ must be the form given in \eqref{I-hL-Solution}. Substituting it in
the transformation $u(x)=\log\mu_n(x)+w(x)+\log(\frac{n^{n-1}\sigma_{n-1}}{2}L)$, we deduce that the solution to \eqref{T-Ot-EL} must be of the form
\begin{equation}\label{T-EL-solution}
w(x)=\log\frac{\big(|x'|^2+(1+t)^2\big)^{\frac{n}{2}}\lambda}{\big(|x'-x_0'|^2+(t+\lambda)^2\big)^{\frac{n}{2}}}-\log L,
\end{equation}
where $\lambda>0$, $L>0$, $x_0'\in\mathbb R^{n-1}$ and $x=(x',t)\in\overline{\mathbb R_+^n}$. The arbitrariness of $L$ suggests that the solution to the Euler-Lagrange equation of inequality \eqref{I-Ot-inq} must be of the form \eqref{Extre-f}. 

Combining the above with Proposition \ref{Best-const}, we conclude that $w\in W_{\mu_n}(\mathbb R_+^n)$ is the extremal function of \eqref{I-Ot-inq} if and only if $w$ takes the form of \eqref{Extre-f}. By now, the proof of Theorem \ref{Ot-inq-thm} is complete. 
\qed

\appendix

\section{\textbf{Some properties of $W_{\mu_n}(\mathbb R_+^n)$}}\label{Section AP-A}

\begin{lemma}\label{equiv lem}
For any $u\in C_c^\infty(\overline{\mathbb R_+^n})$, 
there exists a positive constant $C$, depending only on $n$, such that
\begin{gather}\label{in-le-bd}
\int_{\mathbb R_+^n}|u| \mathrm d\tilde{\mu}_n(x)\leq C\big[\int_{\partial\mathbb R_+^n}|u|\mathrm d\mu_n(x')+(\int_{\mathbb R_+^n}|\nabla u|^n\mathrm dx)^{\frac 1 n}\big],\\
\int_{\partial\mathbb R_+^n}|u|\mathrm d\mu_n(x')\leq C\big[\int_{\mathbb R_+^n}|u| \mathrm d\tilde{\mu}_n(x)+(\int_{\mathbb R_+^n}|\nabla u|^n\mathrm dx)^{\frac 1 n}\big],\label{bd-le-in}
\end{gather}
where $\mathrm d\tilde{\mu}_n(x)=\mu_{n+1}(x) \mathrm dx$.
\end{lemma}

\begin{proof}
Since $1+|x'|^2+t^2\leq (1+t)^2+|x'|^2\leq 2(1+|x'|^2+t^2)$ and 
$$\int_0^{+\infty} \frac{\mathrm dt}{(1+|x'|^2+t^2)^{\frac{n+1} 2}}=\frac{1}{(1+|x'|^2)^{\frac n 2}} \int_0^{\frac\pi 2} \cos^{n-1}\theta \mathrm d\theta ,$$
we have
\begin{equation}\label{formula A-1}
\frac{C_1}{(1+|x'|^2)^{\frac n 2}}\leq \int_0^{+\infty} \mu_{n+1}(x',t)\mathrm dt \leq \frac{C_2}{(1+|x'|^2)^{\frac n 2}},
\end{equation}
where $C_2=2\sigma_n^{-1}\int_0^{\frac\pi 2} \cos^{n-1}\theta \mathrm d\theta$ and $C_1=2^{-\frac{n+1}2}C_2$.

For any $u\in C_c^\infty(\overline{\mathbb R_+^n})$, 
\begin{equation}\label{formula A-2}
|u(x',t)|=|u(x',0)|+\int_0^t\frac{\mathrm d}{\mathrm ds}|u(x',s)|\mathrm ds.
\end{equation}
Multiplying \eqref{formula A-2} by $\mu_{n+1}(x)$ 
and integrating over $\mathbb R_+^n$ yields
\begin{eqnarray}\label{AP-int-esti}
&&\int_{\mathbb R_+^n}|u|\mathrm d\tilde \mu_n(x)\nonumber\\
&\leq & C \int_{\partial\mathbb R_+^n}|u(x',0)|\mathrm d\mu_n(x') + \int_{\partial\mathbb R_+^n}\int_0^\infty\int_0^t|\nabla u(x',s)|\mu_{n+1}(x',t)\mathrm ds\mathrm dt\mathrm dx' \nonumber\\
&=& C \int_{\partial\mathbb R_+^n}|u(x',0)|\mathrm d\mu_n(x') + \int_{\partial\mathbb R_+^n}\int_0^\infty\int_s^\infty|\nabla u(x',s)|\mu_{n+1}(x',t)\mathrm dt\mathrm ds\mathrm dx'.\nonumber\\
&&
\end{eqnarray}
Since 
$$\int_s^\infty \mu_{n+1}(x',t)\mathrm dt\leq  C\int_s^\infty\frac{\mathrm dt}{\big((1+|x'|^2)^\frac1 2+t\big)^{n+1}}\leq C\big((1+|x'|^2)^\frac1 2+s\big)^{-n},$$
where $C$ depends only on $n$, it follows that
\begin{equation}\label{AP-int-esti2}
\int_{\partial\mathbb R_+^n}\int_0^\infty\int_s^\infty|\nabla u(x',s)| \mu_{n+1}(x',t)\mathrm dt\mathrm ds\mathrm dx'\leq C\big(\int_{\mathbb R_+^n}|\nabla u|^n\mathrm dx\big)^{\frac1 n},
\end{equation}
by H\"{o}lder inequality. Substituting 
\eqref{AP-int-esti2} into \eqref{AP-int-esti} yields \eqref{in-le-bd}.

Combining \eqref{formula A-1} and  \eqref{formula A-2} , we have
\begin{eqnarray*}
&&\int_{\partial\mathbb R_+^n}|u(x',0)|\mathrm d\mu_n(x')\leq \frac 1{C_1} \int_{\partial\mathbb R_+^n}|u(x',0)|\int_0^{+\infty} \mu_{n+1}(x',t)\mathrm dt\mathrm dx'\\
&\leq& C 
\big(\int_{\mathbb{R}^n_+} |u(x',t)| \mathrm d\tilde{\mu}_n(x)+ \int_{\partial\mathbb R_+^n}\int_0^\infty\int_s^\infty|\nabla u(x',s)| \mu_{n+1}(x',t)\mathrm dt\mathrm ds\mathrm dx' \big),
\end{eqnarray*}
which deduces \eqref{bd-le-in} by \eqref{AP-int-esti2}.
\end{proof}

\medskip

Define two norms as follows
\begin{gather*}
\|u\|_{In}=\int_{\mathbb R_+^n}|u| \mathrm d\tilde{\mu}_n(x)+\big(\int_{\mathbb R_+^n}|\nabla u|^n\mathrm dx\big)^{\frac 1 n},\\
\|u\|_{Bd}=\int_{\partial\mathbb R_+^n}|u|\mathrm d\mu_n(x')+\big(\int_{\mathbb R_+^n}|\nabla u|^n\mathrm dx\big)^{\frac 1 n},
\end{gather*}
and take two Banach function spaces as follows
\begin{gather*}
W_{In}=\big\{u\in L^1(\mathbb{R}_+^n, \mathrm d\tilde{\mu}_n(x)): \nabla u\in L^n(\mathbb{R}_+^n)\big\},\\
W_{Bd}=\big\{u\in L^1(\partial\mathbb{R}_+^n, \mathrm d{\mu}_n(x')): \nabla u\in L^n(\mathbb{R}_+^n)\big\}.
\end{gather*}
Combining the classical technique, we have the following completion results and will skip the proof.
\begin{proposition}\label{density-prop0}
(1)\ $W_{In}$ and $W_{Bd}$ are the completion spaces of $C_c^\infty(\overline{\mathbb{R}_+^n})$ under the norm $\|\cdot\|_{In}$ and $\|\cdot\|_{Bd}$ respectively; (2)\ $W_{In}=W_{Bd}$.
\end{proposition}

\begin{remark}\label{L1 estimate-rmk} By Proposition \ref{density-prop0}, the inequality of Lemma \ref{equiv lem} still holds for any $u\in W_{Bd}$.
Namely, for any 
$u\in W_{Bd}$, there exists a positive constant $C$, depending only on $n$, such that
\begin{equation*}
\int_{\mathbb R_+^n}|u|\mathrm d\tilde\mu_n(x)\leq C\big(\int_{\partial\mathbb R_+^n}|u|\mathrm d\mu_n(x')+\big(\int_{\mathbb R_+^n}|\nabla u|^n\mathrm dx\big)^{\frac 1 n}\big).
\end{equation*}
\end{remark}
Obviously, we have the following embedding result.

\begin{proposition}\label{density-prop}
$W_{\mu_n}$ is an embedded space of $W_{In}=W_{Bd}$ and $ W_{\mu_n}(\mathbb R_+^n)=\overline{C_c^\infty(\overline{\mathbb R_+^n})}^{\|\cdot\|_{\mu_n}}$ for $n\ge2$.
\end{proposition}

\begin{proposition}\label{Pro energy In}
If $w\in W_{\mu_n}(\mathbb{R}_+^n)$, then $\int_{\mathbb{R}_+^n} e^{\frac{n}{n-1}w} \mu_n^{\frac{n}{n-1}}\mathrm dx<+\infty.$
\end{proposition}
\begin{proof}
As $|x|\rightarrow +\infty$, it is easy to see
\begin{eqnarray*}
&&\nu_n
=O(|x|^{-\frac{n^2}{n-1}}),\qquad \big|\nabla \log(\nu_n)
\big|=O(|x|^{-1}),\\
&& \mu_n^{\frac{n}{n-1}}=O(|x|^{-\frac{n^2}{n-1}}),\qquad |\nabla\log(\mu_n)|=O(|x|^{-1}).    
\end{eqnarray*}
Noting that $\mu_n^{\frac{n}{n-1}}=\mu_n^{\frac{n+1}{n}+\frac{1}{n(n-1)}}\leq C\mu_n^{\frac{n+1}{n}}\leq C\mu_{n+1}$ and by \eqref{I-EO-inq} and \eqref{in-le-bd}, we have
\begin{eqnarray*}
&&\log\big(\int_{\mathbb{R}_+^n} e^{\frac{n}{n-1}w} \mu_n^{\frac{n}{n-1}}\mathrm dx\big) \leq C\log\big(\int_{\mathbb{R}_+^n} e^{\frac{n}{n-1}w}\mathrm d\nu_n\big)\\
&\leq& C\big(\int_{\mathbb{R}_+^n} H_n(x,\frac{n}{n-1}\nabla w)\mathrm dx+\frac{n}{n-1} \int_{\mathbb{R}_+^n} w\mathrm d\nu_n \big)\\
&\leq& C\big[\int_{\mathbb{R}_+^n} |\nabla w|^n\mathrm dx +\int_{\mathbb{R}_+^n} |\nabla w|^2 \big|\nabla \log\big((1+|x|^{\frac{n}{n-1}})^{-n}\big)\big|^{n-2}\mathrm dx\\ 
&& \qquad + \int_{\mathbb{R}_+^n} |w| \mu_n^{\frac{n}{n-1}}\mathrm dx \big]\\
&\leq& C\big(\int_{\mathbb{R}_+^n} |\nabla w|^n\mathrm dx +\int_{\mathbb{R}_+^n} |\nabla w|^2 |\nabla\log\mu_n|^{n-2}\mathrm dx +  \int_{\mathbb{R}_+^n} |w|\mathrm d\tilde{\mu}_n(x)\big)<+\infty.
\end{eqnarray*}
\end{proof}

\begin{proposition}\label{contain-lem}
$
W_{\mu_n}(\mathbb R_+^n)\subset W_{\operatorname{loc}}^{1,n}(\overline{\mathbb R_+^n}).
$
\end{proposition}

\begin{proof}
Let $u\in W_{\mu_n}(\mathbb R_+^n)$.  By Proposition \ref{density-prop}, there exists a sequence $\{u_i\}_{i\in \mathbb N}\subset C_c^{\infty}(\overline{\mathbb R_+^n})$ such that $\|u_i-u\|_{\mu_n}\rightarrow0$ as $i\rightarrow\infty$. In view of Remark \ref{L1 estimate-rmk}, we have $\|u_i-u\|_{L^1(\mathbb R_+^n,\mathrm d\tilde\mu_n)}\rightarrow0$ as $i\rightarrow\infty$, hence $u_i\rightarrow u$ in $L_{\operatorname{loc}}^{1}(\overline{\mathbb R_+^n})$.

Let $K=B_R\cap \overline{\mathbb R_+^n},~R>0$.
We claim that for any $v\in C^1(\overline{K})$, there holds
\begin{equation}\label{claim-inq}
\|v\|_{L^n(K)}\leq C_1\big(\|v\|_{L^1(K\cap\partial\mathbb R_+^n)}+\|\nabla v\|_{L^n(K)}\big),
\end{equation}
where $C_1$ is a positive constant depending only on $K$.

Next, we prove the claim 
by contradiction. Suppose the contrary of the claim is true, that is, there exists a subsequence $\{v_k\}_{k\in\mathbb N}\subset C^1(\overline{K})$ such that $\|v_k\|_{L^n(K)}=1$ and
$$\|v_k\|_{L^1(K\cap\partial\mathbb R_+^n)}+\|\nabla v_k\|_{L^n(K)}<\frac1 k.$$
Hence $\|v_k\|_{W^{1,n}(K)}\leq C$, where $C$ is a positive constant independent of $k$. We obtain that up to a subsequence, 
\begin{eqnarray*}
v_k\rightharpoonup v_0,&\quad& \text{in}\quad W^{1,n}(K),\\
v_k\rightarrow v_0,&\quad& \text{in}\quad L^n(K), \\
v_k\rightarrow v_0,&\quad& \text{in}\quad L^1(\partial K),
\end{eqnarray*}
where $\nabla v_0=0$ and $v_0|_{K\cap\partial\mathbb R_+^n}=0$, which implies $v_0=0$. This contradicts $\|v_0\|_{L^n(K)}=\lim\limits_{k\rightarrow\infty}\|v_k\|_{L^n(K)}=1$. We complete the proof of the claim.

By \eqref{claim-inq} and $u_i\in C_c^{\infty}(\overline{\mathbb R_+^n})$, 
\begin{eqnarray*}
\|u_i\|_{L^n(K)}\leq C_1\big(\|u_i\|_{L^1(K\cap\partial\mathbb R_+^n)}+\|\nabla u_i\|_{L^n(K)}\big)\leq C\|u_i\|_{\mu_n}\leq C.
\end{eqnarray*}
This implies $\|u_i\|_{W^{1,n}(K)}\leq C$, where $C$ is a positive constant independent of $i$. There exists a subsequence ( still write it as $u_i$) such that 
\begin{eqnarray*}
u_i\rightharpoonup u_0&\quad& \text{in}\quad  W^{1,n}(K),\\
u_i\rightarrow u_0&\quad& \text{in}\quad  L^n(K),
\end{eqnarray*} for some $u_0\in W^{1,n}(K)$. Since $u_i\rightarrow u$ in $L_{\operatorname{loc}}^{1}(\overline{\mathbb R_+^n})$, we deduce $u=u_0\in W^{1,n}(K)$, that is, $u\in W^{1,n}_{\operatorname{loc}}(\overline{\mathbb R_+^n})$.
\end{proof}

\section{\textbf{Proof of \eqref{A-supsolution-inq}}}\label{Section AP-B}
In this section, we give the details of the proof of \eqref{A-supsolution-inq}. All notation in this section is consistent with that in the proof of Proposition \ref{1order-esti-prop}.

Recall that $\overline{u}_{\varepsilon}$ is denoted in \eqref{A-supsolution-defi}. A straightforward computation reveals that
\begin{equation*}
\partial_{i}\overline{u}_{\varepsilon}=C_{1}^{\frac{1}{n-1}}\phi_{a_{\varepsilon},b}'\cdot\big[\frac{x_{i}}{|x|}(1-\frac{\delta t}{|x|^{\delta+1}})\big],\qquad i=1,\cdots,n-1,
\end{equation*}
\begin{equation*}
\partial_{t}\overline{u}_{\varepsilon}=C_{1}^{\frac{1}{n-1}}\phi_{a_{\varepsilon},b}'\cdot\big[\frac{t}{|x|}(1-\frac{\delta t}{|x|^{\delta+1}})+\frac{1}{|x|^{\delta}}\big],
\end{equation*}
and hence
\begin{equation*}
\nabla \overline{u}_{\varepsilon}=C_{1}^{\frac{1}{n-1}}\phi_{a_{\varepsilon},b}'\cdot\big(\frac{x'}{|x|}(1-\frac{\delta t}{|x|^{\delta+1}}),\frac{t}{|x|}(1-\frac{\delta t}{|x|^{\delta+1}})+\frac{1}{|x|^{\delta}}\big),
\end{equation*}
which implies
\begin{eqnarray*}
|\nabla \overline{u}_{\varepsilon}|^2&=&C_{1}^{\frac{2}{n-1}}\big(\phi_{a_{\varepsilon},b}'\big)^2\cdot\big[(1-\frac{\delta t}{|x|^{\delta+1}})^2+\frac{1}{|x|^{2\delta}}+\frac{2t}{|x|^{\delta+1}}(1-\frac{\delta t}{|x|^{\delta+1}})\big]\nonumber\\
&=&C_{1}^{\frac{2}{n-1}}\big(\phi_{a_{\varepsilon},b}'\big)^2\cdot\big[1+2(1-\delta)\frac{t}{|x|^{\delta+1}}+(\delta^2-2\delta)\frac{t^2}{|x|^{2\delta+2}}+\frac{1}{|x|^{2\delta}}\big]\nonumber\\
&=&C_{1}^{\frac{2}{n-1}}\big(\phi_{a_{\varepsilon},b}'\big)^2\cdot(1+A),
\end{eqnarray*}
where $A=2(1-\delta)\frac{t}{|x|^{\delta+1}}+(\delta^2-2\delta)\frac{t^2}{|x|^{2\delta+2}}+\frac{1}{|x|^{2\delta}}$. It follows that
\begin{eqnarray}\label{AP-boundary-inq1}
|\nabla \overline{u}_{\varepsilon}|^{n-2}\frac{\partial \overline{u}_{\varepsilon}}{\partial t}\big|_{t=0}&=&C_{1}(-\phi_{a_{\varepsilon},b}')^{n-2}\phi_{a_{\varepsilon},b}'(|x|)\cdot(1+\frac{1}{|x|^{2\delta}})^{\frac{n-2}{2}}\frac{1}{|x|^{\delta}}\nonumber\\
&\leq&-2C_{1}(\gamma-n)^{-1}(a_{\varepsilon}-R_{1}^{\frac{n-\gamma}{2}})\frac{1}{|x|^{n-1+\delta}}\nonumber\\
&\leq&-\frac{C_0}{|x|^{\frac{n-1+\beta}{2}}}
\end{eqnarray}
and
\begin{eqnarray}\label{AP-u-eq}
-\Delta_n\overline{u}_{\varepsilon}&=&C_{1}\operatorname{div}\big[(-\phi_{a_{\varepsilon},b}')^{n-1}\cdot(1+A)^{\frac{n-2}{2}}X\big]\nonumber\\
&=&C_{1}\big\{\big[(-\phi_{a_{\varepsilon},b}')^{n-1}\big]'\cdot(1+A)^{\frac{n}{2}}+\frac{n-2}{2}(-\phi_{a_{\varepsilon},b}')^{n-1}(1+A)^{\frac{n}{2}-2}\nabla A\cdot X\nonumber\\
&&+(-\phi_{a_{\varepsilon},b}')^{n-1}\cdot(1+A)^{\frac{n-2}{2}}\operatorname{div}X\big\},
\end{eqnarray}
where $X=\big(\frac{x'}{|x|}(1-\frac{\delta t}{|x|^{\delta+1}}),\frac{t}{|x|}(1-\frac{\delta t}{|x|^{\delta+1}})+\frac{1}{|x|^{\delta}}\big)$. Computing directly gives 
\begin{eqnarray}\label{AP-divX-esti}
\operatorname{div}X&=&\frac{n-1}{|x|}-\delta(n-\delta)\frac{t}{|x|^{\delta+2}}=\frac{n-1}{|x|}+O(|x|^{-\delta-1})\nonumber\\
&=&\frac{n-1}{|x|+\frac{t}{|x|^{\delta}}}+O(|x|^{-\delta-1}),
\end{eqnarray}
as $|x|\rightarrow+\infty$. The definition of $A$ implies that
\begin{eqnarray*}
\nabla A&=&\Big(\frac{2x'}{|x|^{\delta+2}}\big[(\delta^2-1)\frac{t}{|x|}-(\delta+1)(\delta^2-2\delta)\frac{t^2}{|x|^{\delta+2}}-\delta\frac{1}{|x|^{\delta}}\big],\nonumber\\
&&\frac{2t}{|x|^{\delta+2}}\big[(\delta^2-1)\frac{t}{|x|}-(\delta+1)(\delta^2-2\delta)\frac{t^2}{|x|^{\delta+2}}-\delta\frac{1}{|x|^{\delta}}\big]\nonumber\\
&&+\frac{2}{|x|^{\delta+1}}\big[1-\delta+(\delta^2-2\delta)\frac{t}{|x|^{\delta+1}}\big]\Big).
\end{eqnarray*}
Following these calculations, we observe that
\begin{equation}\label{AP-A-esti}
0\leq A\leq\frac{C}{|x|^{\delta}},~|\nabla A|\leq\frac{C}{|x|^{\delta+1}}\text{ and }|X|\leq C,\qquad |x|\geq R_{1}\gg1,
\end{equation}
where $C$ is a positive constant depending only on $\delta$. It follows from \eqref{A-phi-inq}, \eqref{AP-divX-esti} and \eqref{AP-A-esti} that for $|x|\geq R_{1}$
\begin{eqnarray}\label{AP-AX-esti1}
&&\frac{n-2}{2}(-\phi_{a_{\varepsilon},b}')^{n-1}(1+A)^{\frac{n}{2}-2}\nabla A\cdot X\leq\frac{C}{(|x|+\frac{t}{|x|^{\delta}})^{n-1}}|x|^{-(\delta+1)}\leq C|x|^{-n-\delta}\nonumber\\
&&
\end{eqnarray}
and
\begin{eqnarray}\label{AP-AX-esti2}
(-\phi_{a_{\varepsilon},b}')^{n-1}\cdot(1+A)^{\frac{n-2}{2}}\operatorname{div}X=\frac{n-1}{|x|+\frac{t}{|x|^{\delta}}}(-\phi_{a_{\varepsilon},b}')^{n-1}+O(|x|^{-n-\delta}),
\end{eqnarray}
as $|x|\rightarrow+\infty$, where the positive constant $C$ is independent of $\varepsilon$. By \eqref{A-phi-eq}, \eqref{AP-AX-esti2} and $1\leq (1+A)^{\frac{n}{2}}\leq 1+\frac{C}{|x|^{\delta}}$, we deduce
\begin{eqnarray}\label{AP-AX-esti3}
&&\big[(-\phi_{a_{\varepsilon},b}')^{n-1}\big]'\cdot(1+A)^{\frac{n}{2}}+(-\phi_{a_{\varepsilon},b}')^{n-1}\cdot(1+A)^{\frac{n-2}{2}}\operatorname{div}X\nonumber\\
&&=(|x|+\frac{t}{|x|^{\delta}})^{-\frac{n+\gamma}{2}}+O(|x|^{-n-\delta}),
\end{eqnarray}
as $|x|\rightarrow+\infty$. Substituting \eqref{AP-AX-esti1} and \eqref{AP-AX-esti3} into \eqref{AP-u-eq} yields that
\begin{equation}\label{AP-inter-inq}
-\Delta_n\overline{u}_{\varepsilon}=C_{1}(|x|+\frac{t}{|x|^{\delta}})^{-\frac{n+\gamma}{2}}+O(|x|^{-n-\delta})\geq C|x|^{-\frac{n+\gamma}{2}},\qquad |x|\geq R_{1},
\end{equation}
where $C$ is a positive constant independent of $\varepsilon$. Combining \eqref{AP-boundary-inq1} and \eqref{AP-inter-inq}, we obtain the validity of \eqref{A-supsolution-inq}.

\vskip 1cm
\noindent {\bf Acknowledgements}\\
This project is supported by  the National Natural Science Foundation of China (Grant No. 12471109, 12141105), Youth Innovation Team of Shaanxi Universities and the Fundamental Research Funds for the Central Universities (Grant No. GK202307001, GK202402004).

\bigskip

%


\end{document}